%% file: main.tex
\documentclass[3p,computermodern,10pt]{elsarticle}
\input{packages}

\input{commands}
\begin{document}
\begin{frontmatter}

\title{On the impact of dimensionally-consistent and physics-based inner products for POD-Galerkin and least-squares model reduction of compressible flows}
\author[a]{Eric J. Parish}
\author[b]{Francesco Rizzi}


\address[a]{Sandia National Laboratories,  Livermore, CA}
\address[b]{Nex Gen Analytics, Sheridan, WY}

\begin{abstract}
  Model reduction of the compressible Euler equations based on proper orthogonal decomposition (POD) and Galerkin orthogonality or least-squares residual minimization requires the selection of inner product spaces in which to perform projections and measure norms. The most popular choice is the vector-valued $L^2(\Omega)$ inner product space. This choice, however, yields dimensionally-inconsistent reduced-order model (ROM) formulations which often lack robustness. In this work, we try to address this weakness by studying a set of dimensionally-consistent inner products with application to the compressible Euler equations. First, we demonstrate that non-dimensional inner products have a positive impact on both POD and Galerkin/least-squares ROMs. Second, we further demonstrate that physics-based inner products based on entropy principles result in drastically more accurate and robust ROM formulations than those based on non-dimensional $L^2(\Omega)$ inner products.As test cases, we consider the following well-known problems: the one-dimensional Sod shock tube, the two-dimensional Kelvin-Helmholtz instability and two-dimensional homogeneous isotropic turbulence.
\end{abstract}

\end{frontmatter}

\input{intro}

\input{math_inner_prod}
\input{roms}
\input{results}
\begin{appendix}
\input{appendix}
\end{appendix}
\bibliographystyle{siam}
\bibliography{refs}

\end{document}

%% file: packages.tex
\usepackage{fullpage}
\usepackage{pifont}
\usepackage[usenames,dvipsnames]{xcolor}
\definecolor{airforceblue}{rgb}{0.36, 0.54, 0.66}
\definecolor{ao}{rgb}{0.0, 0.5, 0.0}
\newcommand{\cmark}{{\color{ao}\ding{51}}}%
\newcommand{\xmark}{{\color{red}\ding{55}}}%

\usepackage[ruled]{algorithm2e}
\usepackage{algpseudocode}
\usepackage[utf8]{inputenc}
\usepackage[english]{babel}
\usepackage{amsthm}

\newtheorem{remark}{Remark}[section]
\usepackage{color}
\usepackage{amsmath}
\usepackage{amsfonts}
\usepackage{mathtools}
\usepackage[titletoc,title]{appendix}
\usepackage{float}
\usepackage{graphicx} 
\usepackage{epstopdf} 
\usepackage{caption}
\usepackage{subcaption}
\usepackage{pdfpages}
\usepackage{comment}
\usepackage{soul}

\usepackage{hyperref}
\hypersetup{
    colorlinks=true,
    linkcolor=blue,
    filecolor=magenta,
    urlcolor=cyan,
}
\usepackage{tikz}
\usetikzlibrary{positioning, fit, arrows.meta, shapes}

\usetikzlibrary{
  arrows.meta, 
  fit, 
  positioning, 
}
\tikzset{
  neuron/.style={ 
    circle,draw,thick, 
    inner sep=0pt, 
    minimum size=3.5em, 
    node distance=1ex and 2em, 
  },
  group/.style={ 
    rectangle,draw,thick, 
    inner sep=0pt, 
  },
  io/.style={ 
    neuron, 
    fill=gray!15, 
  },
  conn/.style={ 
    -{Straight Barb[angle=60:2pt 3]}, 
    thick, 
  },
}

%% file: commands.tex
\newcommand{\FR}[1]{\textcolor{red}{#1}}

\newcommand{\residualConserved}{\mathcal{R}^{\qConserved}}
\newcommand{\residualEntropy}{\mathcal{R}^{\qEntropy}}

\newcommand{\RR}[1]{\mathbb{R}^{#1}}
\newcommand{\RSym}[1]{\mathbb{S}^{#1}}

\newcommand{\Flux}{\boldsymbol F}
\newcommand{\uDummy}{\boldsymbol u}
\newcommand{\vDummy}{\boldsymbol v}
\newcommand{\x}{\boldsymbol x}
\newcommand{\xDiscrete}{\mathbf{x}}
\newcommand{\nSpace}{N}
\newcommand{\rhoInf}{\rho_{\infty}}
\newcommand{\pInf}{p_{\infty}}

\newcommand{\LTwoSymbol}{L^2(\Omega)}

\newcommand{\NonDimensionalLTwoSymbol}{L^2_*(\Omega)}
\newcommand{\EntropySymbol}{L^2_{\conservativeToEntropyInf}(\Omega)}
\newcommand{\EntropyConservativeSymbol}{L^2_{\entropyToConservativeInf}(\Omega)}

\newcommand{\conservativeToEntropy}{\boldsymbol A}
\newcommand{\conservativeToEntropyInf}{\conservativeToEntropy_{\infty}}

\newcommand{\entropyToConservative}{\boldsymbol B}

\newcommand{\entropyToConservativeInf}{\entropyToConservative_{\infty}}

\newcommand{\mass}{\mathbf{M}}
\newcommand{\numStepsInWindow}{n_w}
\newcommand{\numTimeSteps}{N_t}
\newcommand{\numWindows}{N_w}
\newcommand{\genDiscreteInnerProductSpace}{\boldsymbol \ell}
\newcommand{\genDiscreteNorm}[1]{ \| #1 \|_{\genDiscreteInnerProductSpace}}

\newcommand{\snapshotMatrix}{\boldsymbol S}
\newcommand{\snapshotMatrixConserved}{\boldsymbol S^{\qConserved}}
\newcommand{\snapshotMatrixEntropy}{\boldsymbol S^{\qEntropy}}

\newcommand{\GalerkinLTwoConservedRomName}{Galerkin$_{\text{conserved}}$-$\LTwoSymbol$}
\newcommand{\GalerkinNonDimensionalLTwoConservedRomName}{Galerkin$_{\text{conserved}}$-$\NonDimensionalLTwoSymbol$}
\newcommand{\GalerkinLTwoEntropyRomName}{Galerkin$_{\text{entropy}}$-$\LTwoSymbol$}

\newcommand{\WLSLTwoConservedRomName}{WLS$_{\text{conserved}}$-$\LTwoSymbol$}
\newcommand{\WLSNonDimensionalLTwoConservedRomName}{WLS$_{\text{conserved}}$-$\NonDimensionalLTwoSymbol$}
\newcommand{\WLSEntropyLTwoConservedRomName}{WLS$_{\text{conservative}}$-$\EntropyConservativeSymbol$}
\newcommand{\WLSEntropyLTwoEntropyRomName}{WLS$_{\text{entropy}}$-$\EntropyConservativeSymbol$}

\newcommand{\GalerkinLTwoConservedRomNameShort}{Gal$_{\text{cons.}}$-$\LTwoSymbol$}
\newcommand{\GalerkinNonDimensionalLTwoConservedRomNameShort}{Gal$_{\text{cons.}}$-$\NonDimensionalLTwoSymbol$}
\newcommand{\GalerkinLTwoEntropyRomNameShort}{Gal$_{\text{ent.}}$-$\LTwoSymbol$}

\newcommand{\WLSLTwoConservedRomNameShort}{WLS$_{\text{cons.}}$-$\LTwoSymbol$}
\newcommand{\WLSNonDimensionalLTwoConservedRomNameShort}{WLS$_{\text{cons.}}$-$\NonDimensionalLTwoSymbol$}
\newcommand{\WLSEntropyLTwoConservedRomNameShort}{WLS$_{\text{cons.}}$-$\EntropyConservativeSymbol$}
\newcommand{\WLSEntropyLTwoEntropyRomNameShort}{WLS$_{\text{ent.}}$-$\EntropyConservativeSymbol$}

\newcommand{\nSnapshots}{n_S}
\newcommand{\genInnerProduct}[2]{\left(#1,#2\right)_{\genInnerProductSpace}}
\newcommand{\genInnerProductSpace}{\mathcal{I}}

\newcommand{\discreteVelocity}{\boldsymbol f}

\newcommand{\aInf}{a_{\infty}}

\newcommand{\rhouOneInf}{\rho u_{1,\infty}}
\newcommand{\rhouTwoInf}{\rho u_{2,\infty}}
\newcommand{\rhoEInf}{\rho E_{\infty}}

\newcommand{\romDim}{K}

\newcommand{\LTwoInnerProductSpace}{L^2(\physicalDomain)}
\newcommand{\EntropyInnerProductSpace}{L^2_{\conservativeToEntropyInf}(\physicalDomain)}
\newcommand{\EntropyInnerProductSpaceForConservedVariables}{L^2_{\entropyToConservativeInf}(\physicalDomain)}

\newcommand{\NonDimensionalLTwoInnerProductSpace}{L^2_*(\physicalDomain)}

\newcommand{\LTwoInnerProduct}[2]{\left( #1, #2 \right)_{\LTwoSymbol}}
\newcommand{\NonDimensionalLTwoInnerProduct}[2]{\left( #1, #2 \right)_{\NonDimensionalLTwoSymbol}}

\newcommand{\EntropyInnerProduct}[2]{\left( #1, #2 \right)_{\conservativeToEntropyInf}}
\newcommand{\EntropyInnerProductSymbol}[2]{$L_{\mathcal{S}}(\physicalDomain)$}
\newcommand{\EntropyInnerProductForConservedVariables}[2]{\left( #1, #2 \right)_{\entropyToConservativeInf}}

\newcommand{\trialBasis}{\boldsymbol \Phi}
\newcommand{\basisConserved}{\boldsymbol \Phi^{\qConserved}}

\newcommand{\basisEntropy}{\boldsymbol \Phi^{\qEntropy}}

\newcommand{\scalingMatrix}{\mathbf{D}}

\newcommand{\qConserved}{\boldsymbol U}
\newcommand{\qConservedApprox}{\tilde{\boldsymbol U}}
\newcommand{\qConservedDummy}{{\boldsymbol U}^*}

\newcommand{\qConservedDiscrete}{\mathbf{U}}
\newcommand{\qConservedDiscreteApprox}{\tilde{\mathbf{U}}}

\newcommand{\qEntropyDiscreteApprox}{\tilde{\mathbf{V}}}

\newcommand{\qEntropy}{\boldsymbol V}
\newcommand{\qEntropyApprox}{\tilde{\boldsymbol V}}
\newcommand{\qEntropyDummy}{{\boldsymbol V}^*}

\newcommand{\genStateQEntropy}{\hat{\boldsymbol V}}

\newcommand{\genStateQConserved}{\hat{\boldsymbol U}}

\newcommand{\bz}{\boldsymbol 0}
\newcommand{\physicalDomain}{\Omega}
\newcommand{\timeDomain}{[0,T]}

%% file: intro.tex
\section{Introduction}
Projection-based model reduction of fluid flows based on proper orthogonal decomposition (POD) and, e.g., (Petrov-)Galerkin projection has been a research area of great interest for the past several decades. A wealth of literature exists examining applications to incompressible flows~\cite{berkooz_turbulence_pod,iliescu_ciazzo_residual_rom,rb_3,holmes_turbstruct_pod,Lassila2014}, and more recently POD-based reduced-order models (ROMs) have been examined within the context of compressible and multiphysics compressible flows (e.g., compressible flows closed with the Reynolds--averaged Navier-Stokes (RANS) equations~\cite{carlberg_gnat}, reacting flows~\cite{HUANG2022110742}). While these works have demonstrated the potential of POD-based ROMs for more complex configurations, a common trend that has emerged is a lack of robustness.

The challenges for constructing accurate POD ROMs for compressible flows are many. Difficulties include efficiently handling non-polynomial nonlinearities~\cite{eim,deim,everson_sirovich_gappy,carlberg_gnat}, developing effective sampling techniques~\cite{qdeim_drmac,adeim_peherstorfer}, properly accounting for boundary conditions~\cite{KaBa10,GuPeSh07}, and a slowly decaying Kolmogorov $n$-width for convection-dominated problems~\cite{LeCa20,Pe20}. We emphasize that all of these issues are significant and require attention. In this work we focus on a fundamental challenge that underpins any ROM formulation: the selection of the inner product space(s). We place a particular focus on dimensional consistency of inner product spaces and their compatibility for use in model reduction of both \textit{dimensional} and \textit{non-dimensional} problem configurations and solvers.

The choice of inner product used in POD and in the (Petrov)-Galerkin projection plays a critical role on ROM accuracy and stability, and it is well known\footnote{but perhaps not well appreciated, or commonly undertaken!} that significantly more care must be taken when choosing an inner product for compressible flows (as opposed to incompressible flows). This aspect is emphasized in Rowley et al.~\cite{ROWLEY2004115}, where the authors highlight that, in an incompressible flow with state variables $\mathbf{u}_1$, $\mathbf{u}_2$, and $\mathbf{u}_3$, a standard vector-valued $L^2(\Omega)$ inner product induces a norm that measures the total kinetic energy of the system. The same cannot be said for a compressible flow with, e.g., conserved variables, where the standard vector-valued $L^2(\Omega)$ inner product induces a norm that not only does not measure any physically meaningful quantity, but is also \textit{dimensionally inconsistent}, i.e., it adds together quantities that have different units. This challenge is amplified for more complex configurations, such as reacting flows, where different state variables span numerous orders of magnitude and computational fluid dynamics solvers are often written to solve the dimensional form of the governing equations. As a result it is unclear which inner product is most appropriate for compressible flows. We further emphasize that the choice of inner product is important for both POD \textit{and} the projection step used to construct the ROM, and it is possible that different inner products could be used in each step.

Several studies have examined the impact of the choice of inner product on the ROM performance, and we aim to provide a brief review here. For starters, Rowley et al.~\cite{ROWLEY2004115} proposed a new, dimensionally-consistent, energy-based inner product for use in POD and Galerkin projections for the isothermal Euler equations. The approach is shown to be effective, but is not directly applicable to the non-isothermal Euler equations. In a similar spirit, Refs.~\cite{BARONE20091932,BaSaTh08,KaBa10,GiPhXa12} developed and investigated symmetrizing inner products and stabilizing projections for the linearized Euler equations discretized in primitive variables; recent work has adapted this inner product within the context of control~\cite{ReWe20}. Although not emphasized, we point out here that the inner products developed in these works are dimensionally consistent. Galerkin ROMs based on these stabilizing inner products have been shown to be more effective than traditional Galerkin methods. Follow on work by Tezaur (f.k.a Kalashnikova) and Barone put forth a framework for developing Galerkin ROMs of the compressible Navier-Stokes equations discretized in entropy variables via a vector-valued $L^2(\Omega)$ inner product~\cite{KaBa11}. This choice of variables is motivated by the fact that (1) it is a symmetrizing transformation and (2) a Galerkin projection process in the vector-valued $L^2(\Omega)$ inner product with entropy variables results in a statement of the second law of thermodynamics, which is akin to energy for incompressible flows. Although again not emphasized in Ref.~\cite{KaBa11}, this formulation results in ROMs that are dimensionally consistent. Results were only presented for scalar conservation laws and the choice of inner product for POD was not discussed. On a similar note, recent work has developed entropy-stable hyper-reduced ROMs~\cite{CHAN2020109789} based on conserved variables with the standard vector-valued $L^2(\Omega)$ inner product and re-projection of the states after mapping to and back from entropy variables. The choice of inner product employed for POD was not explicitly discussed in~\cite{CHAN2020109789}, but it appears a standard $L^2(\Omega)$ inner product was employed. Lastly, while we have focused the above review to the Euler equations, similar ideas exist for other physics. Ref.~\cite{Kaptanoglu_2021}, for instance, develops POD-Galerkin ROMs for magnetohydrodynamics (MHD) based on an ``MHD inner product". Similarly, Ref.~\cite{sockwellThesis} develops POD-Galerkin ROMs of the rotating shallow water equations based on an energy inner product space.

Another body of work that can be closely tied to the choice of inner product is the development of ROMs based on least-squares residual minimization\cite{legresley_1,legresley_2,legresley_3,bui_resmin_steady,rovas_thesis,carlberg_thesis,bui_thesis,bui_unsteady,carlberg_thesis,carlberg_gnat,carlberg_lspg,carlberg_lspg_v_galerkin,PaCa21}. Instead of performing Galerkin projection, these techniques compute solutions that lie within a reduced subspace that minimize the residual of the system of interest. Like the aforementioned entropy-based ROMs, least-squares techniques symmetrize the system and enhance stability. We additionally note that these techniques are often interpreted as Petrov-Galerkin projections with a time-varying test basis. Typically, least-squares ROMs are formulated to minimize the residual in a standard $\LTwoSymbol$-norm\footnote{Strictly speaking, most least-squares have been formulated at the fully discrete level and minimize an $\ell^2$ norm.}. While this choice of norm has resulted in ROMs that are more robust than their Galerkin counterparts, application of the standard $\LTwoSymbol$-norm to compressible flows results in a ROM formulation that is dimensionally inconsistent (e.g., the residual of the momentum equation has different units than the residual of the energy equation, and the standard $\LTwoSymbol$ norm simply adds these two quantities together). This has been increasingly realized in recent years. Washabaugh ans Huang et al.~\cite{HUANG2022110742}, for example, formulate least-squares ROMs that minimizes a scaled residual. Provided that the scaling is done consistently, this process can result in a ROM formulation that is dimensionally consistent. We additionally note that Refs.~\cite{HUANG2022110742} and~\cite{BlRiHo21} provide a similar scaling for POD. While this scaling process results in dimensionally-consistent ROMs, their resulting performance will be, of course, tied to the choice of scaling; such a drawback was pointed out in Rowley et al.~\cite{ROWLEY2004115}.

This work aims to provide a simple but comprehensive overview of these aforementioned issues by studying a set of inner products for performing POD and building Galerkin and least-squares ROMs within the context of the compressible Euler equations. In our study we focus on the dimensional consistency of existing formulations and highlight the advantages of physics-based inner products. Using dimensional arguments, we additionally extend existing ideas in entropy-based ROMs by proposing a new inner product in which to perform POD for entropy variables as well as a new ``entropy-least-squares" ROM based on this entropy inner product. In summary, the contributions of this work are as follows:
\begin{enumerate}
\item A concise framework for building various dimensionally-consistent ROMs through the language of inner products.
\item A new inner product space in which to perform POD based on entropy principles and dimensional arguments.
\item A new least-squares formulation based on an entropy inner product space. This formulation results in solutions that are more accurate, both qualitatively and quantitatively, than solutions obtained from Galerkin and least-squares ROMs obtained in the traditional $\LTwoSymbol$ inner product.
\item Numerical experiments via the open-source software \textsf{pressio-demoapps} demonstrating the impact of the choice of inner product/norm.
\end{enumerate}
Lastly, we emphasize that we do not consider hyper-reduction in this work. Our focus is not to create the most effective ROM, but rather to directly highlight the impact of the inner product/norm on the ROM results.

This paper proceeds as follows.
First, in Section~\ref{sec:math} we outline the compressible Euler equations. Next, in Section~\ref{sec:podInnerProduct} we introduce POD, highlight challenges for systems with multiple states, and introduce a candidate set of inner products in which to perform POD. Next, Section~\ref{sec:roms} outlines several ROM formulations obtained through Galerkin projection and least-squares residual minimization in different inner products. Next Sections~\ref{sec:results}-\ref{sec:numerical_summary} present numerical experiments demonstrating the importance of the choice of inner product. In these experiments we build and deploy ROMs on dimensional and non-dimensional versions of the same problem to evaluate the performance and consistency of the various methods. Conclusions are then provided in Section~\ref{sec:summary}.

%% file: math_inner_prod.tex
\section{The compressible Euler equations}\label{sec:math}

In this work we consider model reduction of the compressible Euler equations in one and two-dimensions. In what follows we provide the formulation for the two-dimensional setting.

The compressible Euler equations on the two-dimensional spatial domain $\physicalDomain$ are given for $t \in [0,T]$ as
\begin{equation}\label{eq:compressible_euler}
\frac{\partial \qConserved}{\partial t} + \nabla \cdot \Flux(\qConserved)  =\bz,
\end{equation}
where $\qConserved: [0,T] \times \Omega \rightarrow \RR{4}$ comprise the density, $\x_1$ and $\x_2$ momentum, and total energy. The inviscid fluxes are given by $\Flux = \left[ \Flux_1 , \Flux_2\right]$ where $\Flux_1(\cdot) \in \RR{4}$ is the inviscid flux in the $\x_1$ direction and $\Flux_2(\cdot) \in \RR{4}$ is the inviscid flux in the $\x_2$ direction. The state and inviscid fluxes are given by
$$
\qConserved = \begin{Bmatrix}
\rho \\ \rho u_1 \\ \rho u_2 \\ \rho E \end{Bmatrix}, \qquad \Flux_1(\qConserved) = \begin{Bmatrix} \rho u_1 \\ \rho u_1^2 +      p \\ \rho u_1 u_2 \\ u_1(E + p) \end{Bmatrix},
\qquad \Flux_2(\qConserved) = \begin{Bmatrix} \rho u_2 \\ \rho u_1 u_2  \\ \rho u_2^2 + p \\ u_2(E + p) \end{Bmatrix},
$$
where we have assumed a calorically perfect gas such that the thermodynamic pressure is given by
$$p = \left(\gamma - 1\right)\left( \rho E - \frac{1}{2} \rho \left( u_1^2 +  u_2^2 \right) \right),$$
with $\gamma = 1.4$ being the heat-capacity ratio. We assume the existence of initial and boundary conditions such that Eq.~\eqref{eq:compressible_euler} is well-posed.

The Euler equations are typically formulated in terms of conservative variables, as in Eq.~\eqref{eq:compressible_euler}. An alternative path relevant to this work is to express the Euler equations in terms of \textit{entropy variables}. Specifically, leveraging the formulation developed by Hughes~\cite{HUGHES1986223}, we introduce a set of entropy variables as
$$\qEntropy = \begin{Bmatrix} \frac{-s}{\gamma - 1} + \frac{\gamma + 1}{\gamma - 1} - \frac{\rho E}{p} \\
\frac{\rho u_1}{p} \\
\frac{\rho u_2}{p} \\
-\frac{\rho}{p}
  \end{Bmatrix}$$
with $s = \ln(p) - \gamma \log(\rho)$. The Euler equations can be written in terms of the entropy variables as
\begin{equation}\label{eq:compressible_euler_entropy}
\conservativeToEntropy \frac{\partial \qEntropy}{\partial t} + \frac{\partial \Flux_1}{\partial \qConserved} \conservativeToEntropy \frac{\partial \qEntropy}{\partial \x_1} +
\frac{\partial \Flux_2}{\partial \qConserved} \conservativeToEntropy \frac{\partial \qEntropy}{\partial \x_2} = \bz,
\end{equation}
where $\conservativeToEntropy \equiv \frac{\partial \qConserved}{\partial \qEntropy} \in \RSym{4}$ is the Jacobian of the variable transformation and $\RSym{N}$ is the space of $N\times N$ symmetric matrices. The analytic form of the Jacobian $\conservativeToEntropy$ is provided in \ref{sec:appendix}. Critically, $\conservativeToEntropy$, $ \frac{\partial \Flux_1}{\partial \qConserved}  \conservativeToEntropy$, and $ \frac{\partial \Flux_2}{\partial \qConserved} \conservativeToEntropy$ are all symmetric matrices, and thus Eq.~\eqref{eq:compressible_euler_entropy} is ``symmetrized". We emphasize that, while the Euler equations expressed in terms of entropy variables~\eqref{eq:compressible_euler_entropy} is analytically equivalent to the Euler equations expressed in terms of conserved variables~\eqref{eq:compressible_euler}, the symmetric nature of the entropy formulation leads to numerical formulations that enjoy, e.g., enhanced stability properties~\cite{HARTEN1983151,HUGHES1986223}.

\begin{remark}
In this work we will consider several ROM formulations that employ entropy variables. In our numerical experiments, these ROMs are constructed from a flow solver that spatially discretizes the conservative from~\eqref{eq:compressible_euler} rather than the entropy form~\eqref{eq:compressible_euler_entropy}. We consider this case as most flow solvers do not discretize in entropy variables. Consequences of this choice are a loss of symmetry and entropy stability at the discrete level.\footnote{With special care these properties can be retained, but we do not consider this case here.}
\end{remark}

\section{Model reduction via POD}\label{sec:podInnerProduct}
Projection-based reduced-order models generate approximate solutions to the system of interest by first restricting the state to live within a low-dimensional trial space. Various approaches exist for constructing this trial space (e.g., the reduced-basis method, space--time methods) and here we focus on the classic spatial-only POD-based approach~\cite{holmes_turbstruct_pod,berkooz_turbulence_pod}.
POD, which is closely related to principle component analysis and the Karhunen Loeve expansion, is a data-compression method that finds a $K$-dimensional orthonormal basis that minimizes the projection error of a dataset in an inner product space $\genInnerProductSpace$ equipped with the inner product $\genInnerProduct{\cdot}{\cdot}$. For a scalar-valued system with a state variable $u(\x,t)$, POD results in the state being approximated as\footnote{We note that some bodies of work consider an affine subspace that, e.g., satisfies initial conditions or boundary conditions~\cite{GuPeSh07}; the ideas presented here generalize to this setting}
$$u(\x,t) = \sum_{i=1}^\romDim \phi_i(\x) \hat{u}_i(t),$$
where $\phi_i : \Omega \rightarrow \RR{}$, $i=1,\ldots,\romDim$ are the POD basis vectors and $\hat{u}_i : [0,T]\rightarrow \RR{}$, $i=1,\ldots,\romDim$ are generalized coordinates.

The situation is more complicated for vector-valued systems with a state variable $\mathbf{u}(\x,t) : \Omega \times [0,T] \rightarrow \RR{M}$ where $M>1$ is the number of state variables (i.e., $M=4$ for the Euler equations in two dimensions). Using the taxonomy of~\cite{ROWLEY2004115}, two approaches exist for applying POD to vector-valued systems like the Euler equations:
\begin{enumerate}
\item Scalar-valued POD: A separate POD basis is constructed for each state variable, i.e., $M$ separate POD problems are solved, one for each for variable.
The advantage of this approach is that each state variable is treated separately and as a result the process is dimensionally consistent. The disadvantages of this approach are that it cannot capture correlations between the different state variables and results in larger ROMs.

\item Vector-valued POD: A global basis is constructed for all the state variables. In this approach one POD problem is solved and the result is used to construct a basis for all the variables. Advantages of this approach are that it can capture correlations between the different state variables and it leads to smaller ROMs. The primary disadvantage of the approach is that it results in dimensionally-inconsistent formulations if special care is not taken when defining the inner product.
\end{enumerate}
In this work we exclusively consider vector-valued POD. The reasons for this are several, including (1) that vector-valued POD appears to be the more common approach, (2) past works have \textit{not} demonstrated that scalar-valued POD outperforms vector-valued POD~\cite{ROWLEY2004115,RoCoMu01,Ro02}, and (3) scalar-valued POD results in ROMs that are of a significantly higher dimension than vector-valued POD.

We now outline vector-valued POD within the context of the compressible Euler equations in two-dimensions. Specifically, given a dataset of solution ``snapshots" $\snapshotMatrix : \physicalDomain \rightarrow \RR{4 \times \nSnapshots}$, where $\nSnapshots$ is the number of samples, vector-valued POD solves the minimization problem
\begin{equation}\label{eq:vector-valued-pod}
\trialBasis = \underset{\{ \trialBasis^* : \Omega \rightarrow \RR{4 \times K} | \genInnerProduct{ \trialBasis^*}{ \trialBasis^*} = \mathbf{I} \} } {\text{arg min} } \sum_{i=1}^{\nSnapshots}\| \snapshotMatrix_i -  \trialBasis^*  \genInnerProduct{ \trialBasis^* }{\snapshotMatrix_i} \|_{\genInnerProductSpace},
\end{equation}
where $\| \cdot \| \equiv \sqrt{\genInnerProduct{\cdot}{\cdot}}$ and, to simplify the notation, we use $\snapshotMatrix_i$ to denote the $i$th column of the snapshot matrix, i.e., the $i$th sample. Assuming that the POD process is employed on snapshots of the conservative variables, the conservative variables are then approximated as
$$\qConserved(\x,t) \approx \trialBasis(\x) \genStateQConserved(t),$$
where $\trialBasis: \physicalDomain \rightarrow \RR{4 \times \romDim}$ is the ROM basis and $\genStateQConserved: [0,T] \rightarrow \RR{\romDim}$ are the generalized coordinates. We note that a similar expression exists if POD was performed on a different variable set, e.g., entropy variables.

The accuracy of the POD basis and the resulting ROM clearly depends on the choice of the inner product space $(\cdot,\cdot)_{\genInnerProductSpace}$. This is particularly true for vector-valued POD, where the minimization problem~\eqref{eq:vector-valued-pod} combines variables with different units. The remainder of this section outlines various inner products within the context of the compressible Euler equations to study this issue.

\subsection{Vector-valued $L^2(\Omega)$ inner product in conserved variables}
We first study the classic vector-valued $L^2(\Omega)$ inner product with conserved variables.
For $\uDummy,\vDummy : \physicalDomain \rightarrow \RR{4}$ we define the vector-valued $L^2(\physicalDomain)$ inner product as
$$\LTwoInnerProduct{ \uDummy}{\vDummy} = \int_{\Omega} \uDummy(\x)^T \vDummy(\x) d\x \equiv \sum_{i=1}^4  \int_{\Omega} \uDummy_i(\x) \vDummy_i(\x) d\x .$$
This inner product is \textit{not} dimensionally consistent for the compressible Euler equations with solution vectors defined in terms of the conserved variables as
$$\LTwoInnerProduct{ \qConserved}{\qConserved} \equiv  \int_{\physicalDomain} \rho(\x) \rho(\x) d\x + \int_{\physicalDomain} \rho u_1(\x) \rho u_1(\x) d\x + \int_{\physicalDomain} \rho u_2(\x) \rho u_2(\x) d\x + \int_{\physicalDomain} \rho E(\x) \rho E(\x) d\x,$$
which does not make sense on dimensional grounds.
Despite this inconsistency, this inner product has been widely employed for model reduction of compressible flows (particularly in the Galerkin/least-squares projection process; see, for example, the residual-minimization formulation in the seminal work~\cite{carlberg_gnat}). The vector-valued $L^2(\Omega)$ inner product results in a trial basis defined as the solution to the minimization problem
\begin{equation}\label{eq:podLTwo}
\basisConserved = \underset{\{ \trialBasis^* : \Omega \rightarrow \RR{4 \times K} | \LTwoInnerProduct{ \trialBasis^*}{ \trialBasis^*} = \mathbf{I} \} } {\text{arg min} } \sum_{i=1}^{\nSnapshots} \big|\big| \snapshotMatrixConserved_i -  \trialBasis^*  \LTwoInnerProduct{ \trialBasis^* }{\snapshotMatrixConserved_i} \big| \big|_{\LTwoInnerProductSpace}
\end{equation}
where $\snapshotMatrixConserved : \physicalDomain \rightarrow \RR{4 \times \nSnapshots}$ is the sample matrix of the conserved state variables. 

%
\subsection{Non-dimensional vector-valued $L^2$ inner product in conserved variables}
The next inner product we consider is similar to the $L^2(\physicalDomain)$ inner product, but this time includes a scaling function that non-dimensionalizes the variables. Specifically, we define a non-dimensional $L^2(\Omega)$ inner product for $\uDummy,\vDummy : \physicalDomain \rightarrow \RR{4}$ as
\begin{equation}\label{eq:nonDimInnerProduct}
\NonDimensionalLTwoInnerProduct{ \uDummy}{\vDummy} = \int_{\physicalDomain} \uDummy(\x)^T \scalingMatrix(\x) \vDummy(\x) d\x
\end{equation}
where $\scalingMatrix: \physicalDomain \rightarrow \RR{4 \times 4}$ is a  diagonal matrix given by
$$ \scalingMatrix  = \begin{bmatrix} \frac{1}{\rhoInf^2} & 0 & 0 & 0 \\ 0 & \frac{1}{(\rhouOneInf)^2} & 0 & 0 \\ 0 & 0 & \frac{1}{(\rhouTwoInf)^2} & 0 \\ 0 & 0 & 0 & \frac{1}{(\rhoEInf)^2} \end{bmatrix}.$$
In the above $\rhoInf , \rhouOneInf, \rhouTwoInf$, and $ \rhoEInf: \physicalDomain \rightarrow \RR{}$ are the reference states for the density, momentum, and energy, and may or may not depend on the spatial location. This inner product is \textit{dimensionally consistent} for the compressible Euler equations with solution vectors defined in terms of the conserved variables as
\begin{multline}
\NonDimensionalLTwoInnerProduct{ \qConserved}{\qConserved} \equiv  \int_{\physicalDomain}\frac{1}{\rhoInf(\x)^2} \rho(\x) \rho(\x) d\x + \int_{\physicalDomain} \frac{1}{(\rhouOneInf)^2}  \rho u_1(\x) \rho u_1(\x) d\x + \\ \int_{\physicalDomain} \frac{1}{(\rhouTwoInf(\x))^2}  \rho u_2(\x) \rho u_2(\x) d\x + \int_{\physicalDomain} \frac{1}{(\rhoEInf)^2}  \rho E(\x) \rho E(\x) d\x.
\end{multline}
This non-dimensional $L^2(\Omega)$ inner product results in a trial basis defined as the solution to the minimization problem
\begin{equation}\label{eq:podNonDimensionalLTwo}
\basisConserved = \underset{\{ \trialBasis^* : \Omega \rightarrow \RR{4 \times K} | \NonDimensionalLTwoInnerProduct{ \trialBasis^*}{ \trialBasis^*} = \mathbf{I} \} } {\text{arg min} } \sum_{i=1}^{\nSnapshots}\| \snapshotMatrixConserved_i -  \trialBasis^*  \NonDimensionalLTwoInnerProduct{ \trialBasis^* }{\snapshotMatrixConserved_i} \|_{\NonDimensionalLTwoInnerProductSpace}
\end{equation}
where $\snapshotMatrixConserved : \physicalDomain \rightarrow \RR{4 \times \nSnapshots}$ is the sample matrix of the conserved state variables.

We emphasize that scaled vector-valued POD processes like the one outlined above have been pursued in the literature. Ref.~\cite{HUANG2022110742}, for instance, performs POD on the snapshots where each state variable is scaled by its $L^2(\Omega)$ norm on the training set. Ref.~\cite{BlRiHo21} pursues a similar approach, but instead scales by the maximum value of each state variable.

\begin{remark} \textbf{Selection of the reference state.}
The non-dimensional inner product requires specification of the reference state. Numerous choices exist, like those described above. To the best of the authors' knowledge, no approach has been demonstrated to be ``superior".
\end{remark}

\subsection{Energy inner product in transformed variables}
An interesting possibility for creating dimensionally-consistent ROMs is based off of the ``energy" inner products proposed in~\cite{ROWLEY2004115}. Unfortunately, the energy inner products proposed in Ref.~\cite{ROWLEY2004115} are only applicable to the isentropic Euler equations, while the energy inner products proposed in Ref.~\cite{FiTeBa15} are not formal inner products, mathematically speaking (positivity is not guaranteed). As a result, we don't consider these energy inner products.

\subsection{Entropy inner product in entropy variables}\label{sec:entropy_pod}
The next inner product we consider is an ``entropy-based" inner product that acts on solution vectors defined in terms of entropy  variables. We start by noting that the product of the entropy variables and the time derivative of the conserved variables results in~\cite{HUGHES1986223}\footnote{We use the $\sim$ notation as numerous definitions of the entropy exist, and different definitions may have constants of proportionality.}
$$\qEntropy^T \frac{\partial}{\partial t} \qConserved \sim \frac{\partial}{\partial t} \rho s.$$
We observe that this product is \textit{dimensionally consistent} and results in a physically meaningful quantity.
Via the chain rule we can write the above as
$$ \qEntropy^T \frac{\partial }{\partial t} \qConserved  = \qEntropy^T \frac{\partial \qConserved}{\partial \qEntropy} \frac{\partial }{ \partial t} \qEntropy.$$
On dimensional grounds we propose an inner product space for use with entropy variables that is based of an inner product of the type $\qEntropy^T \frac{\partial \qConserved}{\partial \qEntropy} \qEntropy$. We note the following important points:
\begin{itemize}
\item The matrix $ \frac{\partial \qConserved}{\partial \qEntropy}$ is symmetric positive definite so can be used to define an inner product.
\item To create an inner product space that is of use for POD, the weighting matrix $ \frac{\partial \qConserved}{\partial \qEntropy}$ needs to be constant (e.g., defined for some background state).
\end{itemize}
Thus, to create our entropy inner product space, we start by defining the matrix
$$\conservativeToEntropyInf(\x) \equiv \frac{\partial \qConserved}{\partial \qEntropy}\left(\qEntropy_{\infty}(\x)\right) $$
where $\qEntropy_{\infty}$ is a reference state. Leveraging the fact that $\conservativeToEntropyInf : \physicalDomain \mapsto \RSym{4}$ is a symmetric positive definite matrix, we can define the inner product
\begin{equation}\label{eq:entropyInnerProduct}
\EntropyInnerProduct{ \uDummy}{\vDummy} = \int_{\physicalDomain} \uDummy(\x)^T \conservativeToEntropyInf(\x) \vDummy(\x) d\x.
\end{equation}
 This inner product is naturally \textit{dimensionally consistent} for the compressible Euler equations with solution vectors defined in terms of the entropy variables, and as a result does not require any additional scaling when deployed on dimensional problems.
This entropy-based inner product results in a trial basis for the entropy variables defined as the solution to the minimization problem
\begin{equation}\label{eq:podEntropy}
\basisEntropy = \underset{\{ \trialBasis^* : \Omega \rightarrow \RR{4 \times K} | \EntropyInnerProduct{ \trialBasis^*}{ \trialBasis^*} = \mathbf{I} \} } {\text{arg min} } \sum_{i=1}^{\nSnapshots}\| \snapshotMatrixEntropy_i -  \trialBasis^*  \EntropyInnerProduct{ \trialBasis^* }{\snapshotMatrixEntropy_i} \|_{\EntropyInnerProductSpace}
\end{equation}
where $\snapshotMatrixEntropy : \physicalDomain \rightarrow \RR{4 \times \nSnapshots}$ is the sample matrix of the entropy variables.

\subsection{Entropy inner product in conserved variables}\label{sec:entropy_pod_conserved}
The last inner product we consider is an entropy-based inner product that acts on solution vectors defined in terms of the conserved variables. Similar to the entropy inner product space presented in the previous section, we derive this inner product based on dimensional arguments. We start by recalling the dimensional consistency of $\qEntropy^T \frac{\partial \qConserved}{\partial \qEntropy} \qEntropy$, and use this dimensional consistency to create an analogous inner product for use with the conserved variables. On dimensional grounds this results in the inner product $\qConserved^T  \frac{\partial \qEntropy}{\partial \qConserved} \qConserved$, which has the same units as $\qEntropy^T \frac{\partial \qConserved}{\partial \qEntropy} \qEntropy$. Similarly to our entropy inner product, we make the following comments:
\begin{itemize}
\item The matrix $\frac{\partial \qEntropy}{\partial \qConserved}$ is positive definite and thus can be used to define an inner product and norm.
\item To create an inner product space that is of use for POD, the weighting matrix $ \frac{\partial \qEntropy}{\partial \qConserved}$ needs to be constant.
\end{itemize}
Thus, to create this entropy inner product space we start by defining the matrix
$$\entropyToConservativeInf(\x) \equiv \frac{\partial \qEntropy}{\partial \qConserved}\left(\qConserved_{\infty}(\x)\right) = \conservativeToEntropyInf^{-1}(\x)$$
where $\qConserved_{\infty}$ is a reference state. Leveraging the fact that $\entropyToConservativeInf : \physicalDomain \mapsto \RSym{4}$ is a symmetric positive definite matrix, we can define the inner product
\begin{equation}\label{eq:entropyInnerProductConserved}
\EntropyInnerProductForConservedVariables{ \uDummy}{\vDummy} = \int_{\physicalDomain} \uDummy(\x)^T \entropyToConservativeInf(\x) \vDummy(\x) d\x.
\end{equation}
 This inner product is \textit{dimensionally consistent} for the compressible Euler equations with solution vectors defined in terms of the conserved variables, and as a result does not require any additional scaling when deployed on dimensional problems. This entropy-based inner product results in a trial basis for the conserved variables defined as the solution to the minimization problem
\begin{equation}\label{eq:podEntropyConserved}
\basisConserved = \underset{\{ \trialBasis^* : \Omega \rightarrow \RR{4 \times K} | \EntropyInnerProductForConservedVariables{ \trialBasis^*}{ \trialBasis^*} = \mathbf{I} \} } {\text{arg min} } \sum_{i=1}^{\nSnapshots}\| \snapshotMatrixConserved_i -  \trialBasis^*  \EntropyInnerProductForConservedVariables{ \trialBasis^* }{\snapshotMatrixConserved_i} \|_{\EntropyInnerProductSpaceForConservedVariables}
\end{equation}
where $\snapshotMatrixEntropy : \physicalDomain \rightarrow \RR{4 \times \nSnapshots}$ is the sample matrix of the entropy variables.

\subsection{Summary of inner products}
In summary, we considered four potential inner products: the classic $\LTwoSymbol$ inner product, a non-dimensional $\NonDimensionalLTwoSymbol$ inner product, and two entropy-based inner products: one for use with solution vectors defined in terms of entropy variables and the other for use with solution vectors defined in terms of conservative variables. The classic $\LTwoSymbol$ inner product is dimensionally inconsistent when used with solution vectors defined in terms of conservative variables, while the non-dimensional $\NonDimensionalLTwoSymbol$ inner product, $\EntropySymbol$ inner product, and $\EntropyConservativeSymbol$ inner product are dimensionally consistent when used with solution vectors defined in terms of the conservative variables, entropy variables, and conservative variables, respectively. Table~\ref{tab:InnerProductSummary} summarizes these properties.
\begin{table}
\begin{center}
\begin{tabular}{ |c|c|c|c| c| }
\hline
& $\LTwoSymbol$ & $\NonDimensionalLTwoSymbol$ & $\EntropyInnerProductSpace$ & $\EntropyInnerProductSpaceForConservedVariables$  \\ \hline
Dimensionally consistent & {\color{red}\xmark} & {\color{ao}\cmark} & {\color{ao}\cmark} & {\color{ao}\cmark}  \\
Solution variables & $\qConserved$ & $\qConserved$ & $\qEntropy$ & $\qConserved$ \\
\hline
\end{tabular}
\caption{Summary of potential inner products for POD}
\label{tab:InnerProductSummary}
\end{center}
\end{table}

\begin{remark}
\textbf{A comment on orthogonality for Galerkin ROMs.}
POD results in bases that are orthonormal in the inner product in which they were constructed. This orthonormality is typically leveraged in a Galerkin projection process to simplify the construction of the ROM. If the POD bases are constructed in an inner product which differs from the inner product used in Galerkin projection, then the typical simplifications stemming from orthonormality are lost. We emphasize here that this is not a significant issue as (1) ROM performance relies on the subspace described by the basis, not orthonormality (e.g., consider two ROMs: one with an orthonormal basis and one with a non-orthonormal basis. So long as the bases have the same range, Galerkin projection will result ROMs that yield the same solution), and (2) the basis can always be re-orthogonalized in the desired inner product via, e.g., a QR decomposition.
\end{remark}

%% file: roms.tex
\section{Reduced-order models}\label{sec:roms}
Thus far we have outlined the construction of a reduced basis for the state variables through POD in different inner products. In addition to construction of this reduced-basis, a ROM requires a projection process applied to the governing equations. This section outlines several such approaches. First, Section~\ref{sec:Galerkin} outlines the construction of the Galerkin ROM through various inner products. Next, Section~\ref{sec:least_squares} outlines several least-squares ROMs that are obtained from residual minimization in norms induced by different inner product spaces. 
\subsection{Galerkin-based ROMs}\label{sec:Galerkin}
Galerkin ROMs rely on a residual orthogonalization process where the residual of the full-order model is restricted to be orthogonal to the trial space in some inner product. Similar to POD, the performance of the Galerkin ROM depends on the definition of the inner product. We now outline the construction of several Galerkin-type ROMs via the inner products defined above.

\subsubsection{Galerkin ROM}\label{sec:galerkin}
The Galerkin ROM is arguably the most popular ROM formulation and, in general, results in accurate and stable solutions for symmetric systems; in this case Galerkin projection can be interpreted to be optimal in a system-specific norm. In the present context the Galerkin ROM with solution vectors defined in terms of the conservative variables is given by
\begin{equation}\label{eq:galerkin}
\LTwoInnerProduct{ \basisConserved_i }{ \basisConserved \frac{\partial}{\partial t} \genStateQConserved} + \LTwoInnerProduct{ \basisConserved_i }{ \nabla \cdot \Flux( \basisConserved \genStateQConserved) } = \bz
\end{equation}
for $i = 1,\ldots,\romDim$. We will refer to this ROM formulation as the \GalerkinLTwoConservedRomName\ ROM. As the compressible Euler equations are not symmetric with respect to the conservative variables, we expect the Galerkin ROM to, in general, perform poorly. Further, this Galerkin ROM is \text{not} dimensionally consistent for the same reasons outlined above. We further emphasize that, even if the POD bases are defined via a dimensionally-consistent inner product, this ROM formulation is still dimensionally inconsistent. We expect the Galerkin ROM to be the least accurate method we explore.

\subsubsection{Non-dimensional Galerkin reduced-order model}
To remedy the dimensional inconsistency of the Galerkin ROM outlined above, we can instead define a non-dimensional Galerkin ROM with solution vectors defined in terms of the conservative variables as
\begin{equation}\label{eq:nonDimensionalGalerkin}
\NonDimensionalLTwoInnerProduct{ \basisConserved_i }{ \basisConserved \frac{\partial}{\partial t} \genStateQConserved} + \NonDimensionalLTwoInnerProduct{ \basisConserved_i }{ \nabla \cdot \Flux( \basisConserved \genStateQConserved) } = \bz
\end{equation}
for $i=1,\ldots,\romDim$. We will refer to this ROM formulation as the \GalerkinNonDimensionalLTwoConservedRomName\ ROM. Unlike Eq.~\eqref{eq:galerkin}, Eq.~\eqref{eq:nonDimensionalGalerkin} \textit{is} dimensionally consistent. Despite being dimensionally consistent, however, we still expect this ROM to perform poorly as the Euler equations are not symmetric with respect to the conserved variables.

\subsubsection{Galerkin reduced-order model in entropy variables}
We can develop a ROM addressing both of the deficiencies of the above Galerkin ROMs, namely dimensional consistency and a lack of symmetry, by creating a Galerkin ROM with solution vectors defined in terms of the entropy variables as
\begin{equation}\label{eq:entropyGalerkin}
\LTwoInnerProduct{ \basisEntropy_i }{  \frac{\partial \qConserved}{\partial \qEntropy} \basisEntropy \frac{\partial}{\partial t} \genStateQEntropy} + \LTwoInnerProduct{ \basisEntropy_i }{ \nabla \cdot \Flux\left( \qConserved\left( \basisEntropy  \genStateQEntropy \right) \right) } = \bz
\end{equation}
for $i=1,\ldots,\romDim$. We will refer to this ROM as the \GalerkinLTwoEntropyRomName\ ROM. This ROM is both dimensionally consistent \textit{and} symmetric with respect to the entropy variables. As a result we expect this ROM to be more accurate and more stable than both the ROMs outlined above.

\subsection{Least-squares reduced-order models}\label{sec:least_squares}
Least-squares ROMs comprise a popular alternative to Galerkin-based ROMs. Unlike Galerkin ROMs, which rely on residual orthogonalization, least-squares ROMs rely on residual minimization, i.e., they compute solutions that lie within a low-dimensional trial subspace that minimize the residual of the system of interest in some norm. A variety of least-squares ROM formulations exist, and here we formulate our least-squares ROMs within the context of the windowed least-squares (WLS) approach~\cite{PaCa21}, which addresses addresses several deficiencies of the more standard least-squares Petrov-Galerkin (LSPG) approach.

To describe WLS ROMs, we first define the residual of the compressible Euler equations expressed in terms of the conserved variables as 
\begin{equation}\label{eq:residualConserved}
\residualConserved : \qConserved \mapsto \frac{\partial}{\partial t} \qConserved + \nabla \cdot \Flux (\qConserved).
\end{equation}
Analogously, we define a residual for the compressible Euler equations in terms of entropy variables as
\begin{equation}\label{eq:residualEntropy}
\residualEntropy : \qEntropy \mapsto \frac{\partial}{\partial t} \qConserved(\qEntropy) + \nabla \cdot \Flux \left(\qConserved \left(\qEntropy \right) \right).
\end{equation}

\begin{figure}[!t]
\begin{center}
\begin{subfigure}[t]{0.85\textwidth}
\includegraphics[trim={0 5cm 0 4cm},clip,width=1.0\linewidth]{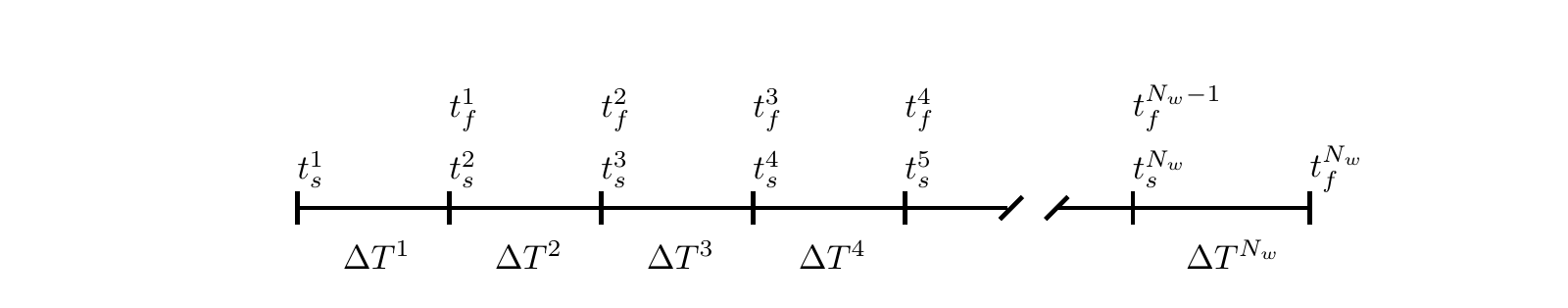}
\end{subfigure}\\
\caption{Depiction of the decomposition of the time domain into windows for WLS.}
\label{fig:wls_fig}
\end{center}
\end{figure}

WLS ROMs operate by partitioning the time domain into a set of $\numWindows$ non-overlapping time windows $[t_s^n,t_f^n] \subset \timeDomain$ of length $\Delta T^n = t_f^n - t_s^n$, $n=1,\ldots,\numWindows$ such that $t_s^1 = 0$ and $t_f^{\numWindows} = T$ and $t_s^{n+1}= t_f^n$, $n=1,\ldots,\numWindows - 1$; see Figure~\ref{fig:wls_fig}. 
WLS then sequentially seeks approximate solutions that lie within the trial subspace that minimize the residual over each window in some norm defined by the inner product space $\mathcal{I}$, i.e., WLS solves the sequence of optimization problems for $n=1,\ldots,\numWindows$
\begin{equation}\label{eq:wlsGeneralProblem}
\qConservedApprox(t) =   \underset{\qConservedDummy(t) \in \text{Range}(\basisConserved) }{\text{arg min} }\int_{t_s^n}^{t_f^n} \genInnerProduct{ \residualConserved\left(\qConservedDummy \right) }{\residualConserved\left(\qConservedDummy\right)} dt
\end{equation}
subject to the initial condition constraints $\qConservedApprox(t_s^n) = \qConservedApprox(t_f^{n-1})$ for $n=2,\ldots,\numWindows$ and $\qConservedApprox(t_s^1) = \qConserved_0$, where $\qConserved_0$ is, e.g., projection of the initial conditions onto the trial subspace. The analogous ROM formulation for a solution vector defined in terms of entropy variables is given by
 \begin{equation}\label{eq:wlsGeneralProblemEntropy}
\qEntropyApprox(t) =   \underset{\qEntropyDummy(t) \in \text{Range}(\basisEntropy) }{\text{arg min} }\int_{t_s^n}^{t_f^n} \genInnerProduct{ \residualEntropy\left(\qEntropyDummy \right) }{\residualEntropy\left(\qEntropyDummy\right)} dt 
\end{equation}
subject to analogous initial condition constraints. 

Various approaches exist to solve the sequence of WLS minimization problems, including \textit{optimize-then-discretize} approaches that involve solving the Euler--Lagrange equations associated with the residual minimization principle, and \textit{discretize-then-optimize} approaches that involve discretizing the objective function in space and time and solving the resulting nonlinear least-squares problem with a discrete least-squares solver (e.g., Gauss--Newton).\footnote{In this case WLS is essentially LSPG, but with the residual defined over multiple time steps.} We consider the latter here due to its simplicity and ease of implementation. In particular, we will consider a discretize-then-optimize scheme where the temporal derivative in Eq.~\eqref{eq:wlsGeneralProblem} (or Eq.~\eqref{eq:wlsGeneralProblemEntropy}) is discretized with the Crank-Nicolson method and the spatial derivatives are discretized with the finite volume method.

Similar to POD, the choice of inner product dictates the accuracy of WLS ROMs. The remainder of this section will outline several possible inner products for WLS ROMs within the context of the Euler equations to investigate this issue.

\subsubsection{$\LTwoSymbol$ inner product}
The most straight-forward inner product to employ in the optimization problems~\eqref{eq:wlsGeneralProblem} and~\eqref{eq:wlsGeneralProblemEntropy} is the standard vector-valued $\LTwoSymbol$ inner product. In this case our ROM formulation with a basis for the conserved variables becomes 
\begin{equation}\label{eq:wlsLTwoProblem}
\qConservedApprox(t) =   \underset{\qConservedDummy(t) \in \text{Range}(\basisConserved) }{\text{arg min} }\int_{t_s^n}^{t_f^n} \LTwoInnerProduct{ \residualConserved\left(\qConservedDummy \right) }{\residualConserved\left(\qConservedDummy\right)} dt.
\end{equation}
A similar definition exists if we discretize in entropy variables.
We will refer to WLS ROMs that employ the above inner product and discretize the conserved variables as \WLSLTwoConservedRomName\ ROMs.  

We emphasize that, like the classic $\LTwoSymbol$ inner product employed in POD, this choice of inner product is \textit{dimensionally inconsistent} for the residual minimization problems~\eqref{eq:wlsGeneralProblem} and~\eqref{eq:wlsGeneralProblemEntropy} as it results in a norm that, e.g., adds the residual of the energy equation to the residual of the momentum equations. Regardless, this inner product as been widely employed in least-squares model reduction (see again Ref.~\cite{carlberg_gnat}). Due to their symmetric nature and tie to an underlying residual minimization problem, least-squares ROMs based on the $\LTwoSymbol$ inner product have, in general, been observed to be more robust than their Galerkin counterparts. 

\subsubsection{Non-dimensional $\NonDimensionalLTwoSymbol$ inner product}
A more consistent approach is to employ a non-dimensional $\NonDimensionalLTwoSymbol$ inner product like that defined in Eq.~\eqref{eq:nonDimInnerProduct}. With this inner product, and discretizing in conserved variables, our ROM formulation becomes
 \begin{equation}\label{eq:wlsNonDimLTwoProblem}
\qConservedApprox(t) =   \underset{\qConservedDummy(t) \in \text{Range}(\basisConserved) }{\text{arg min} }\int_{t_s^n}^{t_f^n} \NonDimensionalLTwoInnerProduct{ \residualConserved\left(\qConservedDummy \right) }{\residualConserved\left(\qConservedDummy\right)} dt. 
\end{equation}
A similar definition again exists if we discretize in entropy variables. We will refer to WLS ROMs that employ the above inner product and discretize the conserved variables as \WLSNonDimensionalLTwoConservedRomName\ ROMs. 

This choice of inner product is dimensionally consistent for the residual minimization problems~\eqref{eq:wlsGeneralProblem} and~\eqref{eq:wlsGeneralProblemEntropy}. We emphasize that non-dimensional inner products have been considered in the literature, see, e.g., Refs.~\cite{HUANG2022110742,BlRiHo21}, and typically have resulted in more robust ROMs than those obtained with direct $\LTwoSymbol$ residual minimization. Despite being dimensionally consistent, this inner product is arguably not a ``physics-based" norm as it simply adds scaled variables of different units. 

\subsubsection{Entropy inner product}
Lastly, we consider WLS ROMs based on entropy inner products. Noting that the residual of the Euler equations has units of $\frac{d}{dt} \qConserved$, we employ the inner product $\EntropyConservativeSymbol$ described in Section~\ref{sec:entropy_pod_conserved}. With this inner product, and discretizing in conserved variables, our ROM formulation becomes
 \begin{equation}\label{eq:wlsEntropyConservedProblem}
\qConservedApprox(t) =   \underset{\qConservedDummy(t) \in \text{Range}(\basisConserved) }{\text{arg min} }\int_{t_s^n}^{t_f^n} \EntropyInnerProductForConservedVariables{ \residualConserved\left(\qConservedDummy \right) }{\residualConserved\left(\qConservedDummy\right)} dt. 
\end{equation}
We will refer to WLS ROMs that employ the $\EntropyConservativeSymbol$ inner product and discretize the conserved variables as \WLSEntropyLTwoConservedRomName\ ROMs. 
A similar definition exists if we discretize in entropy variables, in which case we obtain the ROM formulation
 \begin{equation}\label{eq:wlsEntropyEntropyProblem}
\qEntropyApprox(t) =   \underset{\qEntropyDummy(t) \in \text{Range}(\basisEntropy) }{\text{arg min} }\int_{t_s^n}^{t_f^n} \EntropyInnerProductForConservedVariables{ \residualEntropy\left(\qEntropyDummy \right) }{\residualEntropy\left(\qEntropyDummy\right)} dt. 
\end{equation}
We will refer to ROMs based on the above formulation as \WLSEntropyLTwoEntropyRomName\ ROMs. This choice of inner product is dimensionally consistent for the residual minimization problems~\eqref{eq:wlsEntropyConservedProblem} and~\eqref{eq:wlsEntropyEntropyProblem}, and as a result does not require any additional scaling if deployed on a dimensional problem.

\begin{remark}
We chose to construct our ``entropy norm" based on a fixed background state $\qConserved_{\infty}$. We note that an alternative option would be to left multiply our residual by the Cholesky decomposition of $\frac{\partial \qEntropy}{\partial \qConserved}$ evaluated about the current state and minimize the residual in the standard $\LTwoSymbol$-norm, i.e., we could solve
$$
\qConservedApprox(t) =   \underset{\qConservedDummy(t) \in \text{Range}(\basisConserved) }{\text{arg min} }\int_{t_s^n}^{t_f^n} \LTwoInnerProduct{ \sqrt{ \frac{\partial \qEntropy}{\partial \qConserved} (\qConservedDummy)} \residualConserved\left(\qConservedDummy \right) }{ \sqrt{ \frac{\partial \qEntropy}{\partial \qConserved}(\qConservedDummy) } \residualConserved\left(\qConservedDummy\right)}dt
$$
where we use $\sqrt{\frac{\partial \qEntropy}{\partial \qConserved}}$ to denote the Cholesky decomposition of $\frac{\partial \qEntropy}{\partial \qConserved}$ such that $\sqrt{ \frac{\partial \qEntropy}{\partial \qConserved}}^T  \sqrt{ \frac{\partial \qEntropy}{\partial \qConserved}} = \frac{\partial \qEntropy}{\partial \qConserved}$. In this case our ROM formulation could be interpreted as a ``preconditioned" ROM~\cite{LiFiCa22}. We choose a fixed background state as (1) it directly relates to an inner product space (e.g., that could be used in POD) and (2) it simplifies the numerical implementation, but note that employing the ``preconditioned" version may have advantages. 
\end{remark}
\subsection{Summary of considered ROMs}
In summary we outlined the construction of 6 types of ROMs. The first three ROMs comprised Galerkin-type methods that were obtained by restricting the residual of the governing equations to be orthogonal to the trial space in a standard $\LTwoSymbol$ inner product, a non-dimensional $\LTwoSymbol$ inner product, and in a standard $\LTwoSymbol$ inner product but with a POD basis for entropy variables. The second three ROMs considered were least-squares ROMs that were obtained via residual minimization in a norm induced by the standard $\LTwoSymbol$ inner product, a non-dimensional $\NonDimensionalLTwoSymbol$ inner product, and an ``entropy" $\EntropyConservativeSymbol$ inner product. All of the ROMs considered are dimensionally consistent except for the Galerkin and WLS ROMs employing the $\LTwoSymbol$ inner product. Only the ``entropy" ROMs employ ``physics-based" inner products. Tables~\ref{tab:romSummaryGalerkin} and~\ref{tab:romSummaryWls} summarize these properties. Lastly, we note that additional ROMs could be constructed from the described inner products that we do not consider here. For instance, we could create a Galerkin ROM discretized in conserved variables based on the $\EntropyConservativeSymbol$ inner product; this ROM would be dimensionally consistent and may have better stability properties than the standard Galerkin ROM.
 
\begin{table}
\begin{center}
\begin{tabular}{ |c|c|c|c|}
\hline
& \GalerkinLTwoConservedRomNameShort &  \GalerkinNonDimensionalLTwoConservedRomNameShort & \GalerkinLTwoEntropyRomNameShort  \\ \hline
Dimensionally consistent & \xmark & \cmark & \cmark \\ \hline
Physics-based & \xmark & \xmark & \cmark  \\
\hline
\end{tabular}
\caption{Summary of Galerkin ROMs considered.}
\label{tab:romSummaryGalerkin}
\end{center}
\end{table}

\begin{table}
\begin{center}
\begin{tabular}{ |c| c|c|c|c | }
\hline
& \WLSLTwoConservedRomNameShort & \WLSNonDimensionalLTwoConservedRomNameShort & \WLSEntropyLTwoConservedRomNameShort & \WLSEntropyLTwoEntropyRomNameShort \\ \hline
Dimensionally consistent & \xmark & \cmark & \cmark & \cmark \\ \hline
Physics-based &  \xmark & \xmark & \cmark & \cmark\\
\hline
\end{tabular}
\caption{Summary of WLS ROMs considered.}
\label{tab:romSummaryWls}
\end{center}
\end{table}

\section{Current results for numerical analyses of ROMs for the compressible Euler equations}
Before proceeding, we briefly summarize available results (to the best of the authors' knowledge) for numerical analyses of the above ROMs applied to the compressible Euler equations. We note that analyses for the Galerkin ROMs deployed on simplified versions of the Euler equations (e.g., the isothermal Euler equations, the linearized Euler equations) exist in the literature~\cite{KaBa10,ROWLEY2004115}, but here we restrict our discussion to the fully nonlinear Euler equations. 

\begin{itemize}
\item \textbf{\GalerkinLTwoConservedRomName.}
There are presently no theoretical results guaranteeing stability of the Galerkin ROM with a basis for the conserved variables applied to the compressible Euler equations (numerical results have demonstrated that the Galerkin ROM is often unstable). There are additionally no practical, sharp \textit{a priori} or \textit{a posteriori} error bounds that exist. We note that analyses for the Galerkin ROM have been undertaken for a generic dynamical system in Ref.~\cite{carlberg_lspg_v_galerkin}; these analyses rely on the availability of difficult-to-compute Lipshitz constants and sharpness was not addressed. 

\item \textbf{\GalerkinNonDimensionalLTwoConservedRomName.}
The non-dimensional Galerkin ROM, to the best of the authors knowledge, has not been considered in the literature. As a result there are no available analyses; we fully expect the non-dimensional Galerkin ROM to suffer from the same shortcomings of the standard Galerkin ROM.

\item \textbf{\GalerkinLTwoEntropyRomName.}
The Galerkin ROM with a basis for the entropy variables can leverage (in part) past analyses by e.g.,  Hughes, Tadmor, and Harten~\cite{HUGHES1986223,tadmor1,tadmor2}, that was performed in the context of finite elements. The \GalerkinLTwoEntropyRomName\ ROM has been analyzed within the ROM context in Ref.~\cite{KaBa11}, where Tezaur and Barone demonstrate that the Galerkin ROM discretized in entropy variables is ``entropy stable", meaning that it can be shown to satisfy the second law of thermodynamics (the Clausius-Duhem inequality). In order to be entropy-stable in practice the \GalerkinLTwoEntropyRomName\ ROM must employ an entropy stable discretization of the governing equations (e.g., via entropy conserving fluxes~\cite{IsRo09}). We remark that, even if the \GalerkinLTwoEntropyRomName\ ROM is entropy stable, it is not totally clear if this guarantees numerical stability. Further, no error bounds have been derived for the \GalerkinLTwoEntropyRomName\ ROM. 

\item \textbf{\WLSLTwoConservedRomName.}
Analyses for WLS ROMs have been carried out for a generic dynamical system in Ref.~\cite{PaCa21}, but similar to the Galerkin ROM these analyses rely on difficult-to-compute Lipshitz constants and may not be sharp. These analyses indicated that minimizing the residual over larger windows resulted in stable and more accurate ROMs, but in general it is not clear if stability can be guaranteed. A similar ROM formulation that preceded WLS is the least-squares Petrov Galerkin (LSPG) ROM, and analyses of LSPG applied to a generic dynamical system were again presented in Ref.~\cite{carlberg_lspg_v_galerkin}. Stability for LSPG is not guaranteed, and the error bounds suffer from the same issues of sharpness and difficult-to-compute Lipshitz constants.  

\item \textbf{\WLSNonDimensionalLTwoConservedRomName, \WLSEntropyLTwoConservedRomName, and \WLSEntropyLTwoEntropyRomName.}
The WLS ROM in these norms has not been considered in the literature, and as a result there are no available analyses. 
\end{itemize}

In summary, none of the ROMs described in this work are equipped with guarantees of stability when they are applied to the compressible Euler equations, nor are they equipped with \textit{a priori} or \textit{a posteriori} error bounds. We believe this should be an area for future work. 

%% file: results.tex
\section{Description of numerical experiments}\label{sec:results}
We now present numerical experiments studying the impact of the choice of inner product on the POD bases and the resulting ROM formulation. We consider three numerical experiments: the Sod shock tube in one dimension, the Kelvin Helmholtz instability in two dimensions, and homogeneous isotropic turbulence in two dimensions. For each problem we create a ``non-dimensional" and ``dimensional" configuration. These configurations are designed to operate at the same Mach number and non-dimensional time step, and are defined over the same non-dimensional time horizon. Methods that are dimensionally consistent should produce identical results on the two configurations, while methods that are not dimensionally consistent should produce different results on the two configurations. All ROMs considered are reproductive, and our aim is only to study the impact of the inner product.

All three problems are simulated in Python using \href{https://pressio.github.io/pressio-demoapps}{\textsf{pressio-demoapps}} (PDA).
PDA is an open-source project developed at Sandia National Laboratory that is aimed at providing benchmark ROM cases, and is part of the \href{https://pressio.github.io/}{\textsf{Pressio}} ecosystem. 
Features of PDA include a cell-centered finite volume-based solver equipped with several flux schemes, exact Jacobians, sample mesh capabilities, and benchmark cases in one, two, and three dimensions. We emphasize that PDA's Euler solver discretizes the governing equations in their conserved form with conserved variables. For ROMs employing entropy variables, we simply employ the required change of variables before evaluating the spatial fluxes in PDA.
We now describe the specific cases studied in this work.

\subsection{Sod shock tube}
The first problem considered is the Sod shock tube in one dimension, which is a typical benchmark problem for the compressible Euler equations.
The problem is governed by the compressible Euler equations on the 1D domain $\physicalDomain \in [-0.5,0.5]$.
For the spatial discretization, we use a finite volume scheme with 500 uniformly sized cells and Weno5 flux reconstruction.
For the temporal discretization, we employ a low-storage fourth-order Runge-Kutta (RK4) scheme. The initial conditions are
$$ \left[ \rho, \; u_1, \; p \right] =
\begin{cases}
\left[ \rhoInf, \; 0, \; \gamma \pInf \right] & x < 0 \\
\left[ \frac{1}{8} \rhoInf , \;  0 , \frac{\gamma}{10} \pInf \right] & x \ge 0 \end{cases}.$$
The PDA webpage for this case can be found \href{https://pressio.github.io/pressio-demoapps/euler_1d_sod.html}{here}.

\begin{figure}[!h]
\begin{center}
\begin{subfigure}[t]{0.32\textwidth}
\includegraphics[trim={0 0 0 0},clip,width=1.0\linewidth]{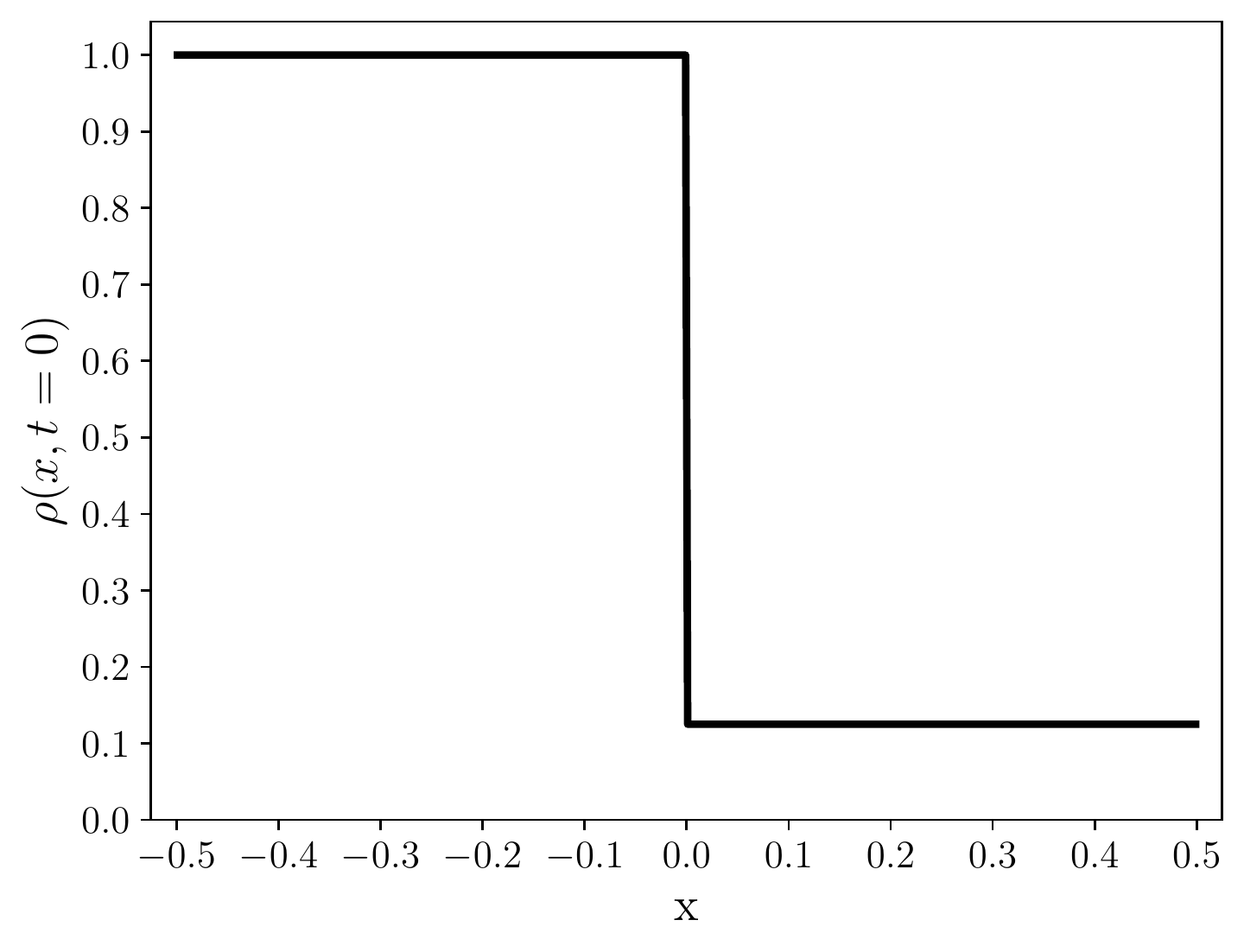}
\end{subfigure}
\begin{subfigure}[t]{0.32\textwidth}
\includegraphics[trim={0 0 0 0},clip,width=1.0\linewidth]{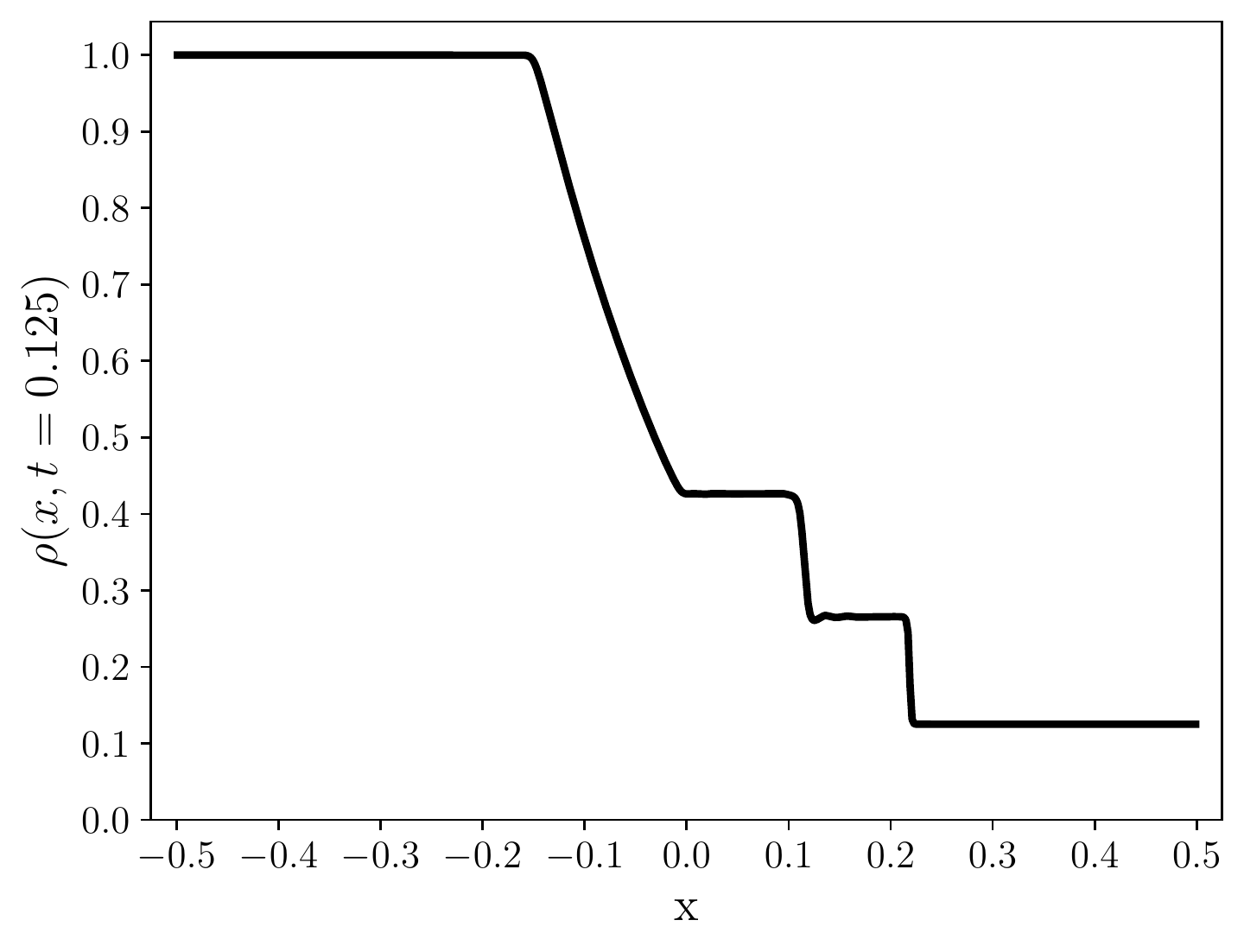}
\end{subfigure}
\begin{subfigure}[t]{0.32\textwidth}
\includegraphics[trim={0 0 0 0},clip,width=1.0\linewidth]{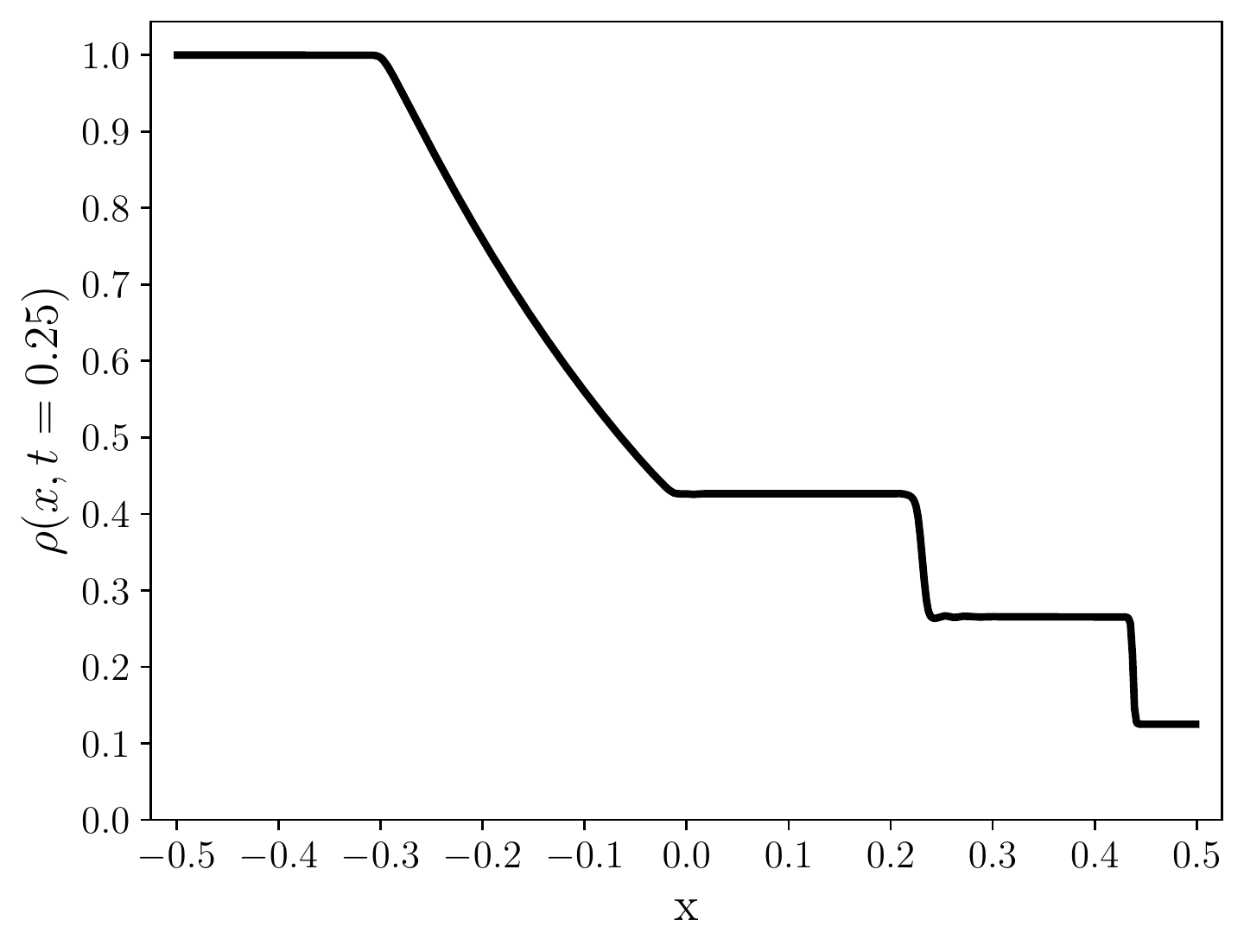}
\end{subfigure}
\caption{Representative plots of the density field at $t=0$ (left),
  $t=0.125$ (middle) and $t=0.25$ (right) obtained for the ``non-dimensional'' Sod 1D FOM.}
\label{fig:fom_sod1d_figs}
\end{center}
\end{figure}

We perform two FOM simulations: a ``non-dimensional" FOM and a dimensional FOM.
The non-dimensional FOM employs $\rhoInf = 1$ and $\pInf = \frac{1}{\gamma}$, while the dimensional FOM employs $\rhoInf = 1.225$ and $\pInf = 101325$.
The FOMs are simulated for $t \in [0,0.25 / a_{\infty}]$, where $a_{\infty} = \sqrt{\gamma \pInf / \rhoInf}$ with a time step of ${\Delta t} = \frac{0.25 \Delta x}{a_{\infty}}$ where $\Delta x = \frac{1}{500}$. We collect solution snapshots every time step for construction of the POD bases. Figure~\ref{fig:fom_sod1d_figs} shows representative plots of the density at three different times simulated using the ``non-dimensional'' FOM.

\subsection{Kelvin Helmholtz instability}
The second example we consider is the Kelvin Helmholtz instability. The problem is governed by the compressible Euler equations on the periodic domain $\physicalDomain \in [-5,5] \times [-5,5]$ and is run for 2.5 flow through times. For spatial discretization we again employ a finite volume scheme with Weno5\footnote{We note that in the Weno reconstruction, a small number is typically added to the denominator of the nonlinear weights. We found that the standard value of $10^{-6}$ led to inconsistencies between our dimensional and non-dimensional configurations for the second two numerical examples. For these two cases we employed a value of $10^{-20}$.} flux reconstruction on a discrete domain comprising 256x256 uniformly spaced volumes. Time integration is performed with a low-storage RK4 scheme. The initial conditions are:
$$ \left[ \rho, \; u_1, \; u_2, \; p \right] =   \begin{cases} \left[ 2 \rhoInf, \; \frac{1}{2} \aInf, \; 0, \; \frac{7}{2} \pInf  \right] & \x \in \Omega_1 \\
\left[ \rhoInf, \; -\frac{1}{2} \aInf, \; 0,\; \frac{7}{2} \pInf \right] & \x \in \Omega_2 \end{cases}$$
where $\Omega_2 =  [-5,5] \times [-2 + \cos( 0.8 \pi x) , 2 + \cos(0.8 \pi x)]$ and $\Omega_1 = \Omega - \Omega_2 $.
The PDA webpage for this case can be found \href{https://pressio.github.io/pressio-demoapps/euler_2d_kelvin_helmholtz.html}{here}.

\begin{figure}[!h]
\begin{center}
\begin{subfigure}[t]{0.319\textwidth}
\includegraphics[trim={0 0 4cm 0},clip,width=1.0\linewidth]{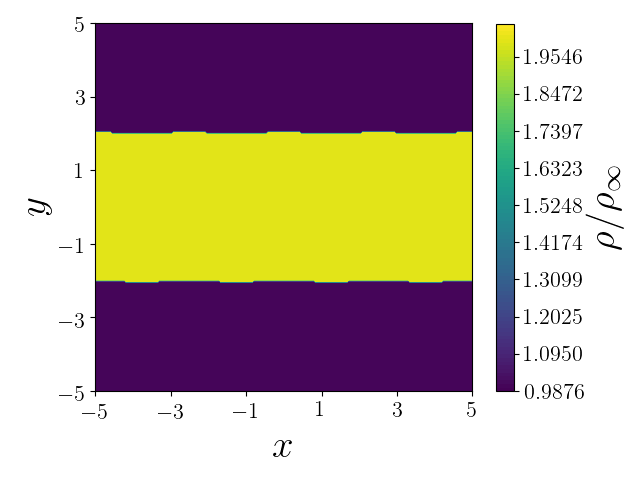}
\caption{$t=0.0$}
\end{subfigure}
\begin{subfigure}[t]{0.28\textwidth}
\includegraphics[trim={1.5cm 0 4cm 0},clip,width=1.0\linewidth]{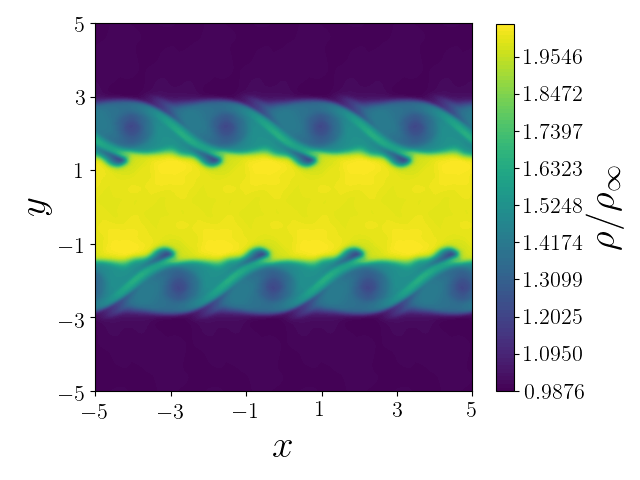}
\caption{$t=24.4$}
\end{subfigure}
\begin{subfigure}[t]{0.383\textwidth}
\includegraphics[trim={1.5cm 0 0 0},clip,width=1.0\linewidth]{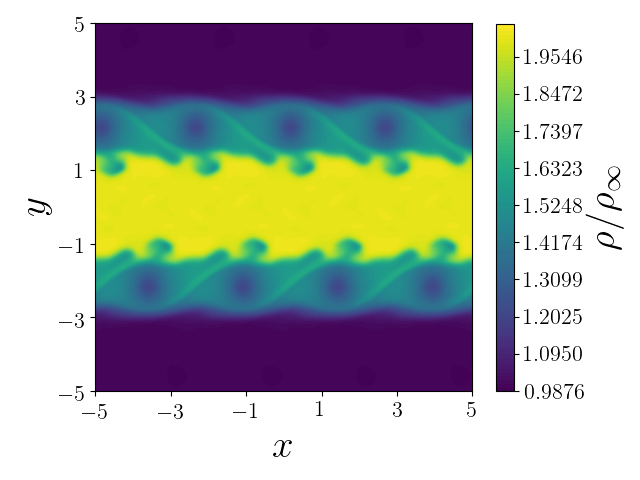}
\caption{$t=48.8$}
\end{subfigure}
\caption{Representative plots of the density at $t=0.0$ (left),
  $t=24.4$ (middle) and $t=48.8$ (right) obtained for the ``non-dimensional'' 2D Kelvin Helmholtz FOM.}
\label{fig:fom_kh2d_figs}
\end{center}
\end{figure}

We again perform two FOM simulations: a ``non-dimensional" FOM with $\rhoInf = 1, \; \pInf = \frac{1}{\gamma}, \; \aInf = 1$, and a ``dimensional" FOM with $\rhoInf = 1.225,\; \pInf = 101325$ and $\aInf =  \sqrt{ \frac{\gamma \pInf }{\rhoInf} } $. The FOMs are simulated for $t \in [0,50/a_{\infty}]$ at a time step of $\Delta t = \frac{0.25 \Delta x}{a_{\infty}}$ with $\Delta x = \frac{10}{256}$. We export solution snapshots every five time steps for construction of the POD bases. Figure~\ref{fig:fom_kh2d_figs} shows representative plots of the density at three different times simulated using the ``non-dimensional'' FOM.

\subsection{Two-dimensional homogeneous turbulence at an infinite Reynolds number}
The final problem setup we consider comprises homogeneous turbulence at an infinite Reynolds number in two dimensions. The problem is again governed by the compressible Euler equations on the periodic domain $\physicalDomain \in [-5,5] \times [-5,5]$. For spatial discretization we employ a finite volume scheme with Weno5 reconstruction on a discrete domain comprising 512x512 uniformly spaced volumes. Time integration is performed with a low-storage RK4 scheme.
\FR{\footnote{This problem is not currently available in PDA, but it will be included in the next release.}}
The problem is initialized in the frequency domain based on the kinetic energy spectrum
$$ e(\mathbf{k}) = \frac{50}{k_p} a u_0^2 \left( \frac{\| \mathbf{k}\|}{k_p} \right)^{2s + 1} \exp\left( -\left(s + \frac{1}{2} \right) \left( \frac{\| \mathbf{k} \|}{k_p} \right)^2 \right) $$
where $\mathbf{k} = \left[ k_1 , k_2 \right]$ is the wave number, $\| \mathbf{k} \| = \sqrt{k_1^2 + k_2^2}$ is the wave number magnitude, and $u_0 = 25 \aInf$, $k_p = 25$, and $s = 3$. The initial density and pressure are set to be uniform with values of $\rhoInf$ and $\pInf$, respectively. The problem is initialized in physical space by first prescribing a vorticity field with the desired energy spectrum and then constructing the $x_1$ and $x_2$ velocity fields based on the vorticity field; see Ref.~\cite{Ch97} for details.

\begin{figure}[!h]
\begin{center}
\begin{subfigure}[t]{0.319\textwidth}
\includegraphics[trim={0 0 4cm 0},clip,width=1.0\linewidth]{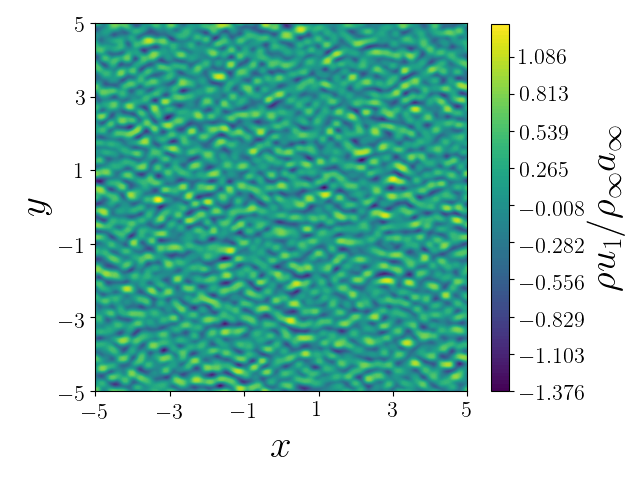}
\caption{$t=0.0$}
\end{subfigure}
\begin{subfigure}[t]{0.28\textwidth}
\includegraphics[trim={1.5cm 0 4cm 0},clip,width=1.0\linewidth]{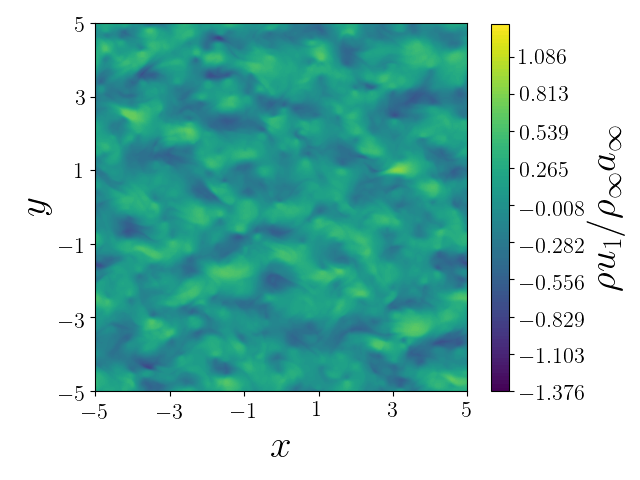}
\caption{$t=9.8$}
\end{subfigure}
\begin{subfigure}[t]{0.383\textwidth}
\includegraphics[trim={1.5cm 0 0 0},clip,width=1.0\linewidth]{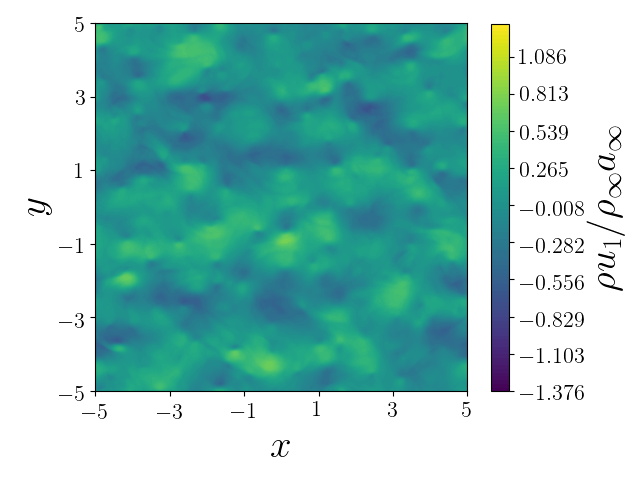}
\caption{$t=19.5$}
\end{subfigure}
\caption{Representative plots of the $x_1$ momentum at $t=0.0$ (left),
  $t=9.8$ (middle) and $t=19.5$ (right) obtained for the ``non-dimensional'' 2D homogeneous turbulence FOM.}
\label{fig:fom_hit2d_figs}
\end{center}
\end{figure}

We again perform two FOM simulations: a ``non-dimensional" FOM with $\rhoInf = 1, \; \pInf = \frac{1}{\gamma}, \; \aInf = 1$, and a ``dimensional" FOM with $\rhoInf = 1.225,\; \pInf = 101325$ and $\aInf =  \sqrt{ \frac{\gamma \pInf }{\rhoInf} }$.
The FOMs are simulated for $t \in [0,20/a_{\infty}]$ at a time step of $\Delta t = \frac{0.25 \Delta x}{a_{\infty}}$ with $\Delta x = \frac{10}{512}$.
We export solution snapshots every five time steps for construction of the POD bases. Figure~\ref{fig:fom_hit2d_figs} shows representative plots of the $x_1$ momentum at three different times simulated using the ``non-dimensional'' FOM.

\subsection{Reduced-order models studied}
In our experiments we study the performance of 7 ROM configurations using the various inner products described above. Tables~\ref{tab:galerkin_rom_summary} and~\ref{tab:wls_rom_summary} summarize these ROMs and identify the variables that the basis represents, the POD problem that is solved to obtain the bases, and the formulation used to construct the ROM. All ROMs employ the same finite volume discretization used by the FOMs. All Galerkin ROMs employ the same RK4 time stepping scheme used by the FOMs, while all WLS ROMs employ the Crank Nicolson method for time discretization.\footnote{The error incurred by temporal discretization is small for the CFL numbers we consider in the experiments, and the difference between the Crank Nicolson and RK4 time schemes is minimal. We employ Crank Nicolson for convenience of evaluating the residual in the least-squares minimization problem arising in WLS.} The time step for all ROMs is set to be the same as the FOM time step. For all WLS ROMs we employ a uniform window size. For the 1D Sod problem the window size is set to be $\Delta T = 10\Delta t$. For the 2D Kelvin Helmholtz and homogeneous turbulence problems the window size is set to be $\Delta T = 2 \Delta t$. The nonlinear least-squares problem resulting from the WLS formulation is solved using \textsf{scipy}'s \href{https://docs.scipy.org/doc/scipy/reference/generated/scipy.optimize.least_squares.html}{\textsf{scipy.optimize.least\_squares}} solver with default settings. We emphasize again that ROM run time is not an aspect we are considering in this work; we are focused on the impact of the inner product space on the ROM performance. Lastly, for each problem we study ROMs of varying basis dimensions. For the 1D Sod problem we study ROMs of basis dimensions $K = 5,10,15,20,25,30,35,40,45,50$, while for the 2D Kelvin Helmholtz and homogeneous turbulence problems we study ROMs of dimensions $K=25,50,100,$ and $150$.

We now discuss construction of the inner product spaces. The above experiments were all constructed such that the non-dimensional configurations have a reference speed of sound $a_{\infty} = 1$ based on a thermodynamic state of $\rho_{\infty} = 1, p_{\infty} = \frac{1}{\gamma}$. Similarly, the ``dimensional" configurations were constructed to have a reference speed of sound of $a_{\infty} = 340.3$ based on a thermodynamic state of $\rho_{\infty} = 1.225$ and $p_{\infty} = 101325$. We use these reference variables for the non-dimensional inner product~\eqref{eq:nonDimInnerProduct}. For the non-dimensional configuration we set $\rho_{\infty} = 1$, $\rho_{\infty} u_{1,\infty} = \rho_{\infty} a_{\infty} $, $\rho_{\infty} u_{2,\infty} = \rho_{\infty} a_{\infty}$, and $\rho E_{\infty} = \rho_{\infty} a_{\infty}^2$ with $a_{\infty} = 1$. For the dimensional configuration we set $\rho_{\infty} = 1.225$, $\rho u_{1,\infty} = \rhoInf a_{\infty} $, $\rho u_{2,\infty} =\rhoInf a_{\infty}$, and $\rho E_{\infty} = \rho_{\infty} a_{\infty}^2$ with $a_{\infty} = 340.3$. This selection of reference states is not unique, and it is possible the other choices may give better results. We additionally note that this choice of reference states makes the $\NonDimensionalLTwoSymbol$ inner product equivalent to the $\LTwoSymbol$ inner product when deployed on the non-dimensional configuration, which is intentional. Lastly, for the entropy inner products~\eqref{eq:entropyInnerProduct} and~\eqref{eq:entropyInnerProductConserved} we define the background states $\qEntropy_{\infty}$ and $\qConserved_{\infty}$ as their mean over the training snapshots\footnote{We note that this choice of reference state \textit{cannot} be used to select the values of, e.g., $\rho_{\infty} u_{1,\infty}$ for use in the $\NonDimensionalLTwoSymbol$ inner product as it can result in the inner product not being positive definite if, e.g., the mean value of $\rho u_{1}$ is negative. While there are other logical choices, such as employing the norm of the snapshot training data or the mean absolute value, this selection is ultimately subjective and it is not clear if one approach is favored. The entropy inner product is not subject to these issues as positivity is guaranteed, even if the mean value is negative.} used to compute the basis (e.g., no spatial dependence). Again, it is possible that a different choice of background state could yield better results.

\begin{table}[H]
\begin{centering}
\begin{tabular}{l c c c c }
\hline
& \GalerkinLTwoConservedRomName &  \GalerkinNonDimensionalLTwoConservedRomName & \GalerkinLTwoEntropyRomName \\
\hline
Basis variables & $\qConserved$ & $\qConserved$ & $\qEntropy$ &  \\
POD method & Eq.~\eqref{eq:podLTwo} &  Eq.~\eqref{eq:podNonDimensionalLTwo} & Eq.~\eqref{eq:podEntropy}  \\
ROM method &  Eq.~\eqref{eq:galerkin} & Eq.~\eqref{eq:nonDimensionalGalerkin} & Eq.~\eqref{eq:entropyGalerkin}  \\
\hline
\end{tabular}
\caption{Summary of Galerkin ROMs investigated.}
\label{tab:galerkin_rom_summary}
\end{centering}
\end{table}

\begin{table}[H]
\begin{centering}
\begin{tabular}{l c c c c c}
\hline
& \WLSLTwoConservedRomName& \WLSNonDimensionalLTwoConservedRomName & \WLSEntropyLTwoConservedRomName & \WLSEntropyLTwoEntropyRomName \\
\hline
Basis variables & $\qConserved$ & $\qConserved$ & $\qConserved$ & $\qEntropy$ \\
POD method &  Eq.~\eqref{eq:podLTwo} & Eq.~\eqref{eq:podNonDimensionalLTwo} & Eq.~\eqref{eq:podNonDimensionalLTwo} & Eq.~\eqref{eq:podEntropy} \\
ROM method &  Eq.~\eqref{eq:wlsLTwoProblem} & Eq.~\eqref{eq:wlsNonDimLTwoProblem} & Eq.~\eqref{eq:wlsEntropyConservedProblem} & Eq.~\eqref{eq:wlsEntropyEntropyProblem} \\
\hline
\end{tabular}
\caption{Summary of WLS ROMs investigated.}
\label{tab:wls_rom_summary}
\end{centering}
\end{table}

\subsection{Error metrics}
We employ the relative, time-integrated mean-squared error as our error metric for the following experiments.\footnote{We note that we investigated other error metrics, e.g., the $\ell^{\infty}$ error, and found that they resulted in similar conclusions.} We compute separate errors for each conserved solution variable, and thus have four (three in one-dimension) errors that we monitor. For, e.g., the density field, this error metric is defined by
$$e_{\rho}= \frac{ \int_0^T  \int_{\Omega}\left( \tilde{\rho}(\x,t) - \rho(\x,t) \right)^2 d\x dt }{\int_0^T  \int_{\Omega} \rho(\x,t)^2 d\x dt}$$
where $\tilde{\rho}(\x,t)$ is the ROM approximation to the density field. We define analogous errors for the $x_1$-momentum, $x_2$-momentum, and energy fields as $e_{\rho u_1}$, $e_{\rho u_2}$, and $e_{\rho E}$, respectively.

\subsection{Stability}
In our numerical experiments we will use the terms ``stable" and ``unstable". Numerous definitions of stability exist in the literature. In our work, we will refer to a ROM as ``unstable" if it results in a solution containing a NaN, and thus the solver cannot continue. Similarly, we will refer to a ROM as ``stable" if it results in a real-valued solution.

\section{Non-dimensional vs. dimensional $L^2(\Omega)$ inner products}\label{sec:results2}
We first explore the impact of constructing POD-Galerkin and POD-WLS ROMs based on standard non-dimensional and dimensional $L^2(\Omega)$ inner products; we will consider the physics-based ``entropy" inner products in subsequent sections. We emphasize that we consider this setting first as Galerkin and least-squares ROMs based on the standard $L^2(\Omega)$ inner products are by far the most studied ROMs.

\begin{figure}[!t]
\begin{center}
\begin{subfigure}[t]{0.65\textwidth}
\includegraphics[trim={0 0 0 0},clip,width=1.0\linewidth]{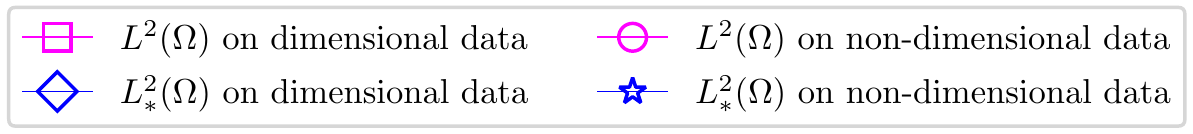}
\end{subfigure}\\
%
\begin{subfigure}[t]{0.33\textwidth}
\includegraphics[trim={0 1.1cm 0 0},clip,width=1.0\linewidth]{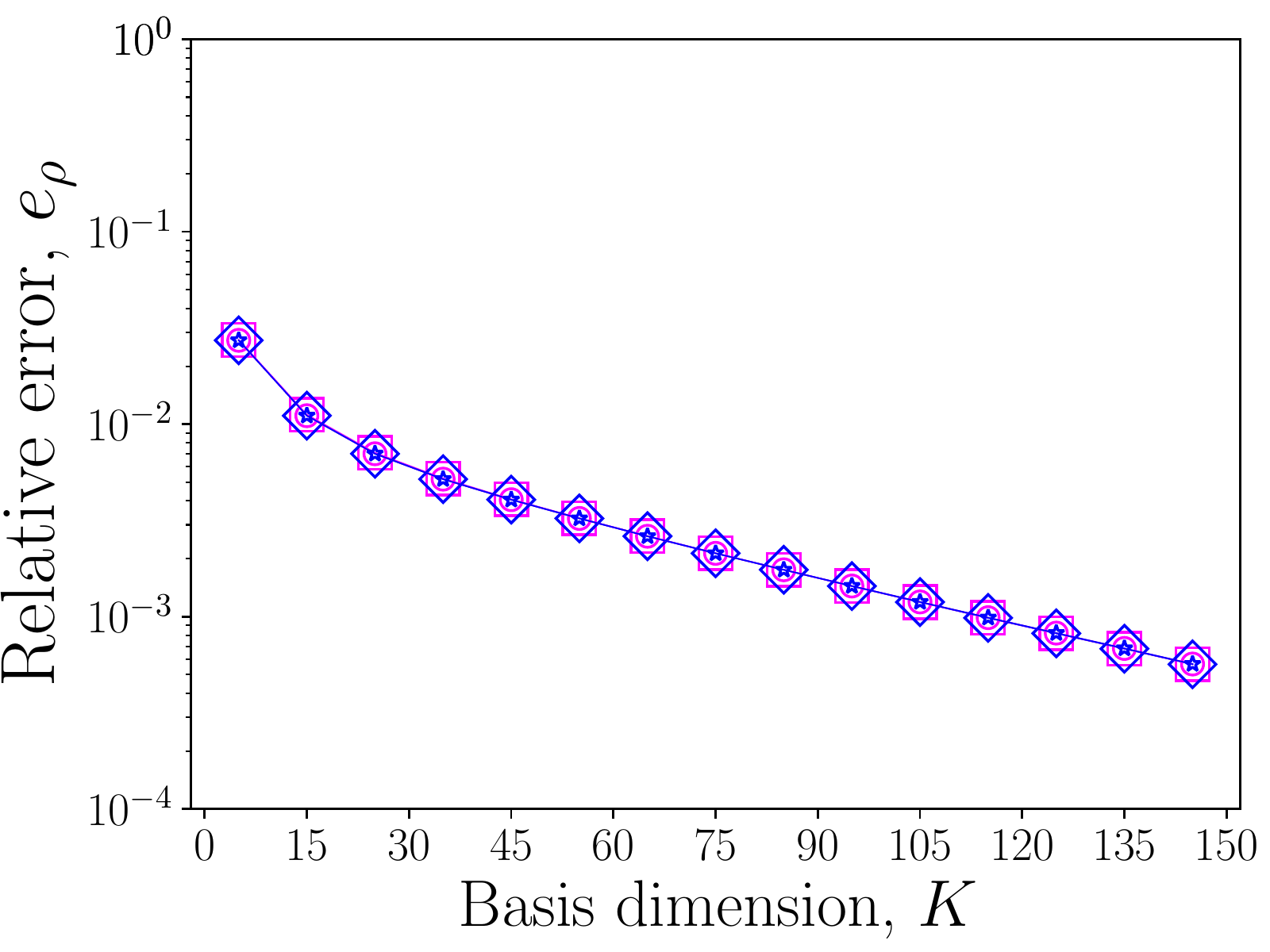}
\end{subfigure}
\begin{subfigure}[t]{0.31\textwidth}
\includegraphics[trim={0cm 0.5cm 0 0},clip,width=1.0\linewidth]{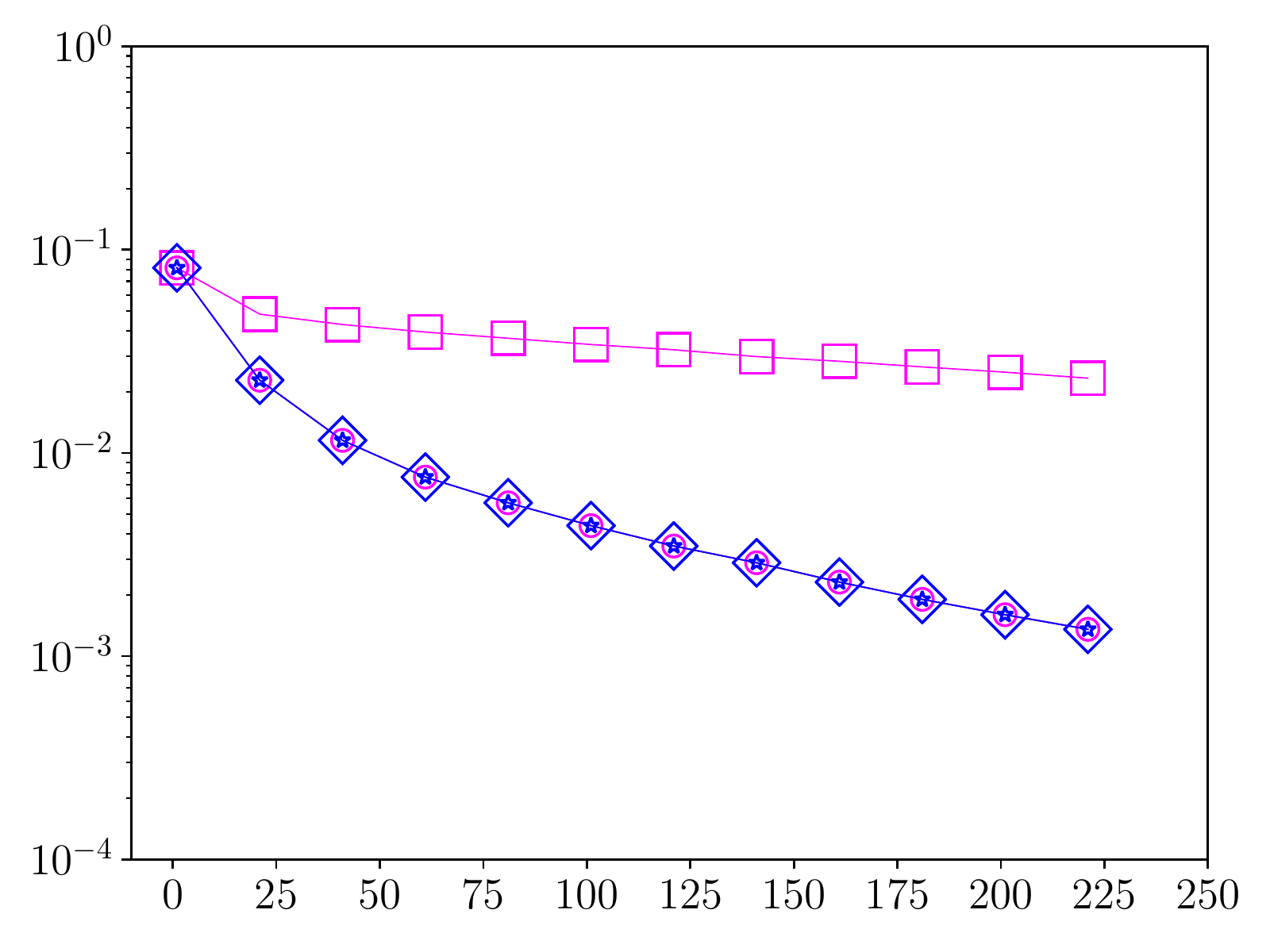}
\end{subfigure}
\begin{subfigure}[t]{0.31\textwidth}
\includegraphics[trim={0cm 0.5cm 0 0},clip,width=1.0\linewidth]{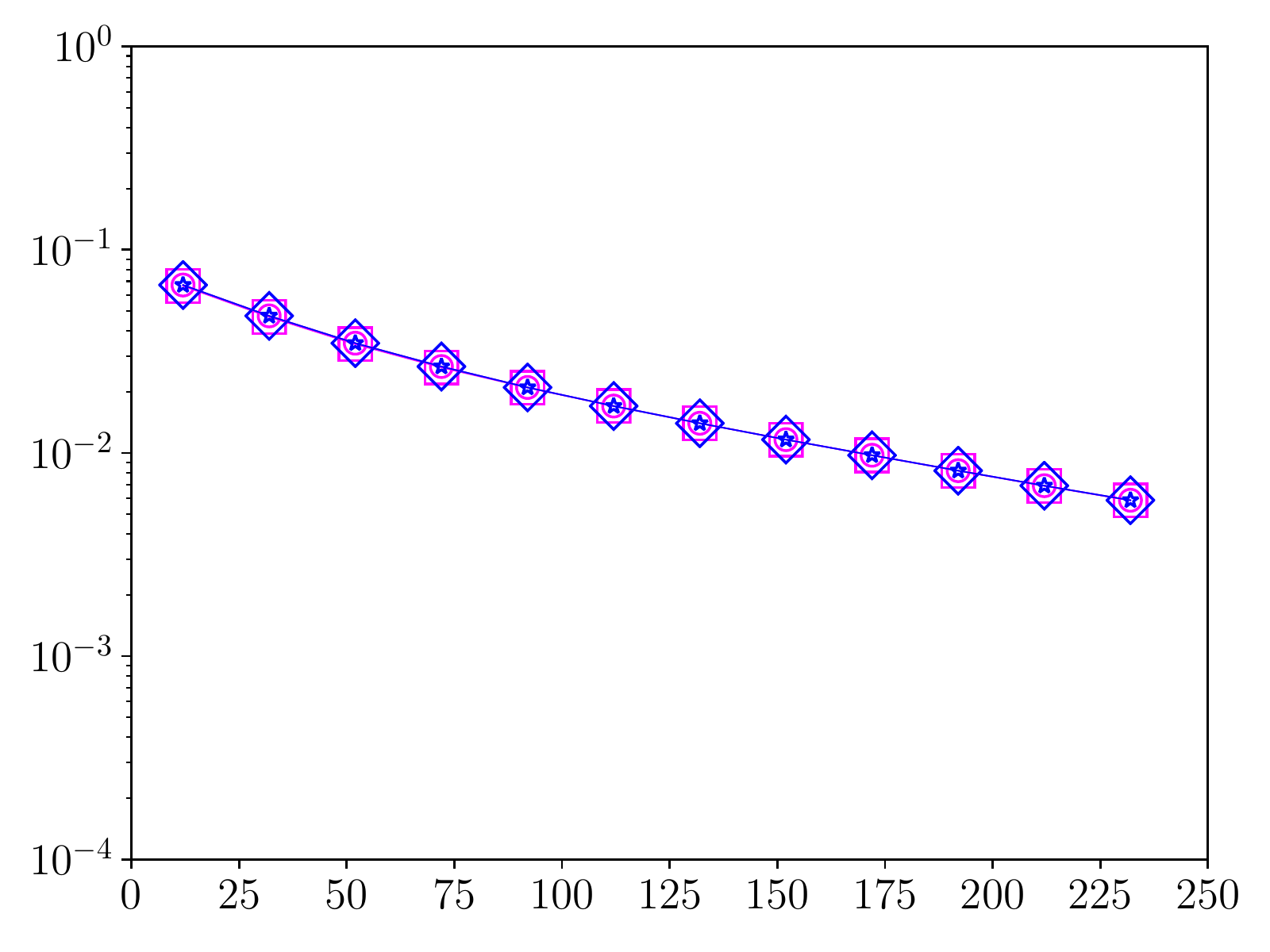}
\end{subfigure}
%
\begin{subfigure}[t]{0.33\textwidth}
\includegraphics[trim={0 1.1cm 0 0},clip,width=1.0\linewidth]{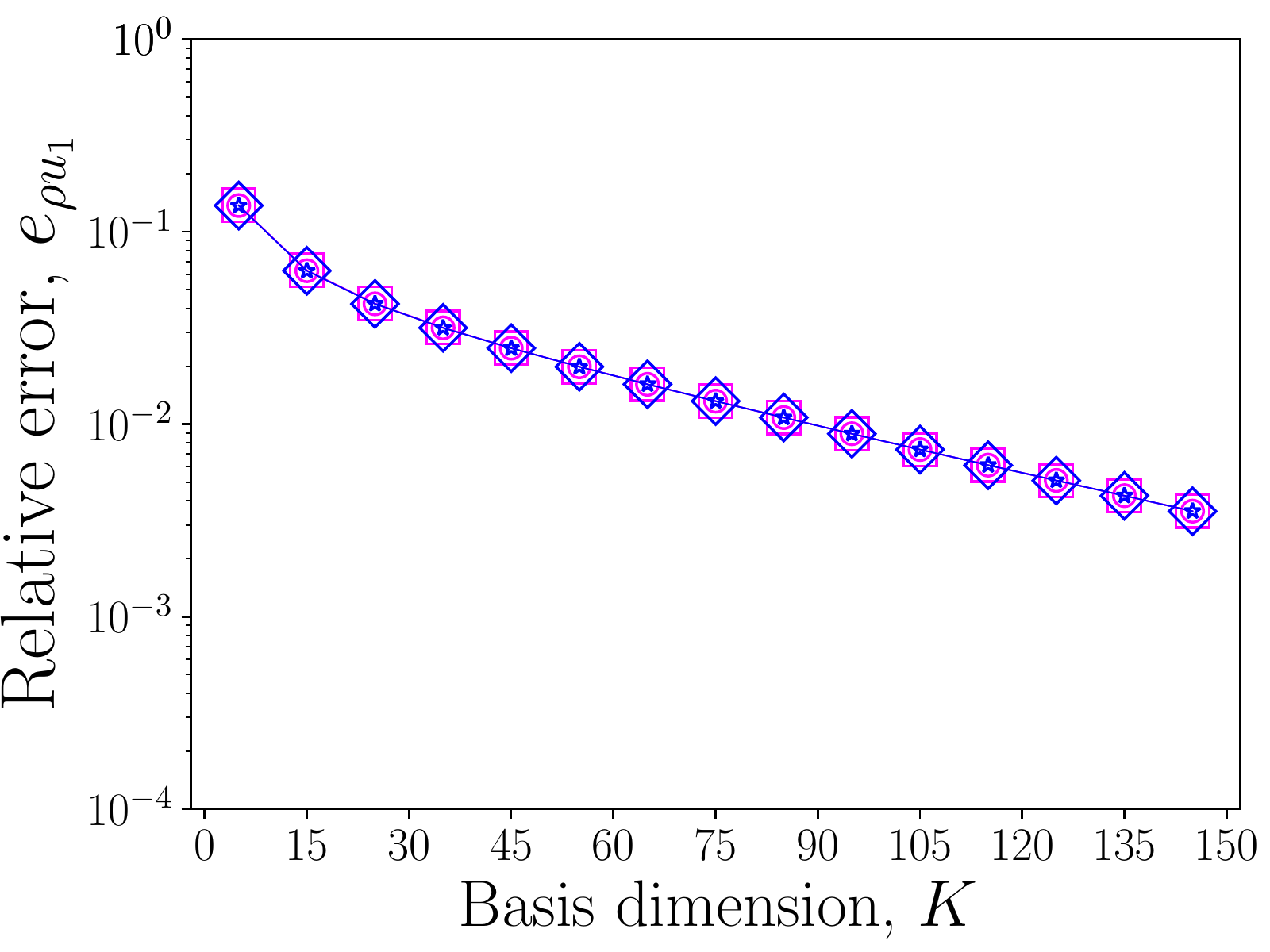}
\end{subfigure}
\begin{subfigure}[t]{0.31\textwidth}
\includegraphics[trim={0cm 0.5cm 0 0},clip,width=1.0\linewidth]{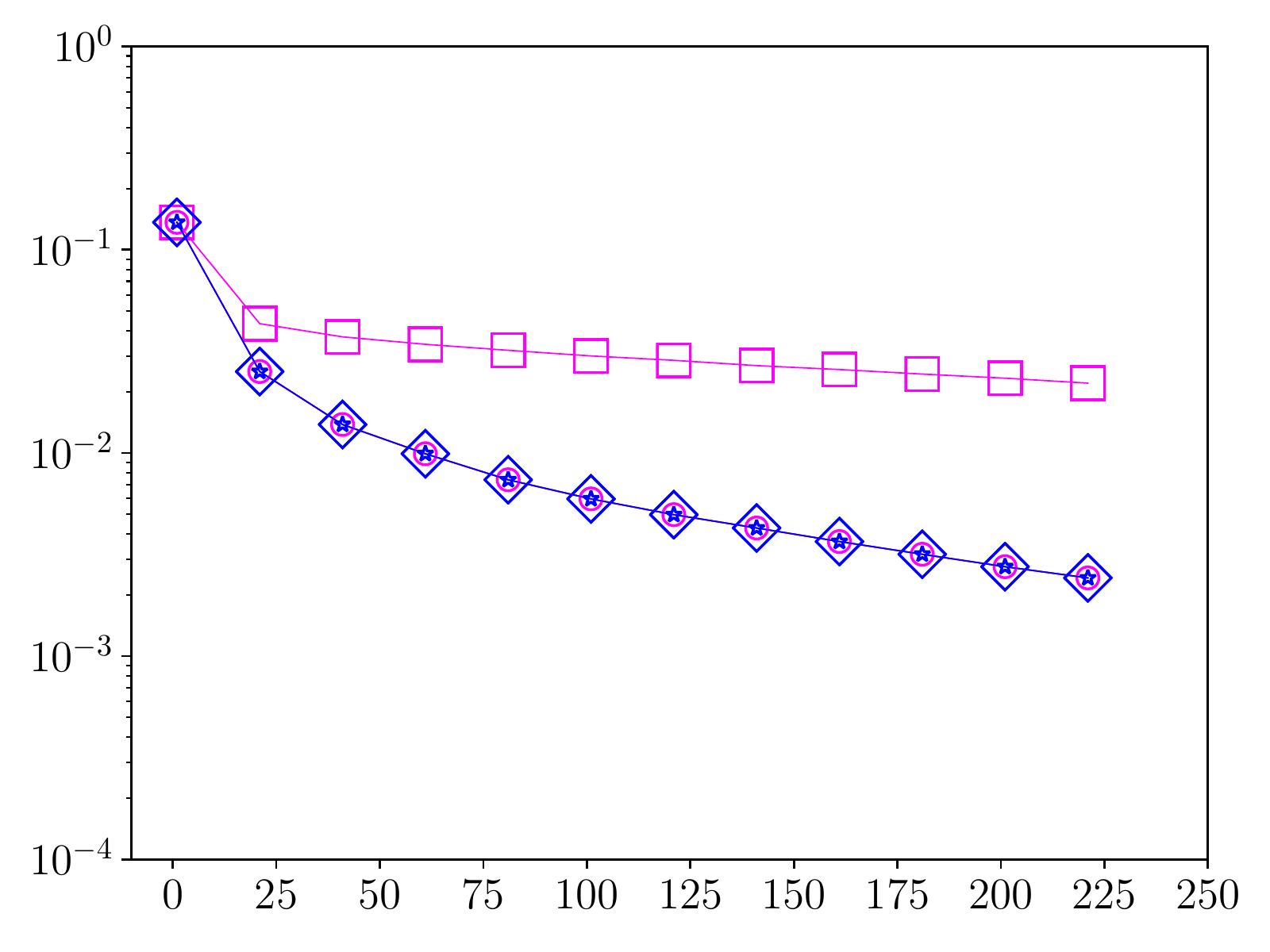}
\end{subfigure}
\begin{subfigure}[t]{0.31\textwidth}
\includegraphics[trim={0cm 0.5cm 0 0},clip,width=1.0\linewidth]{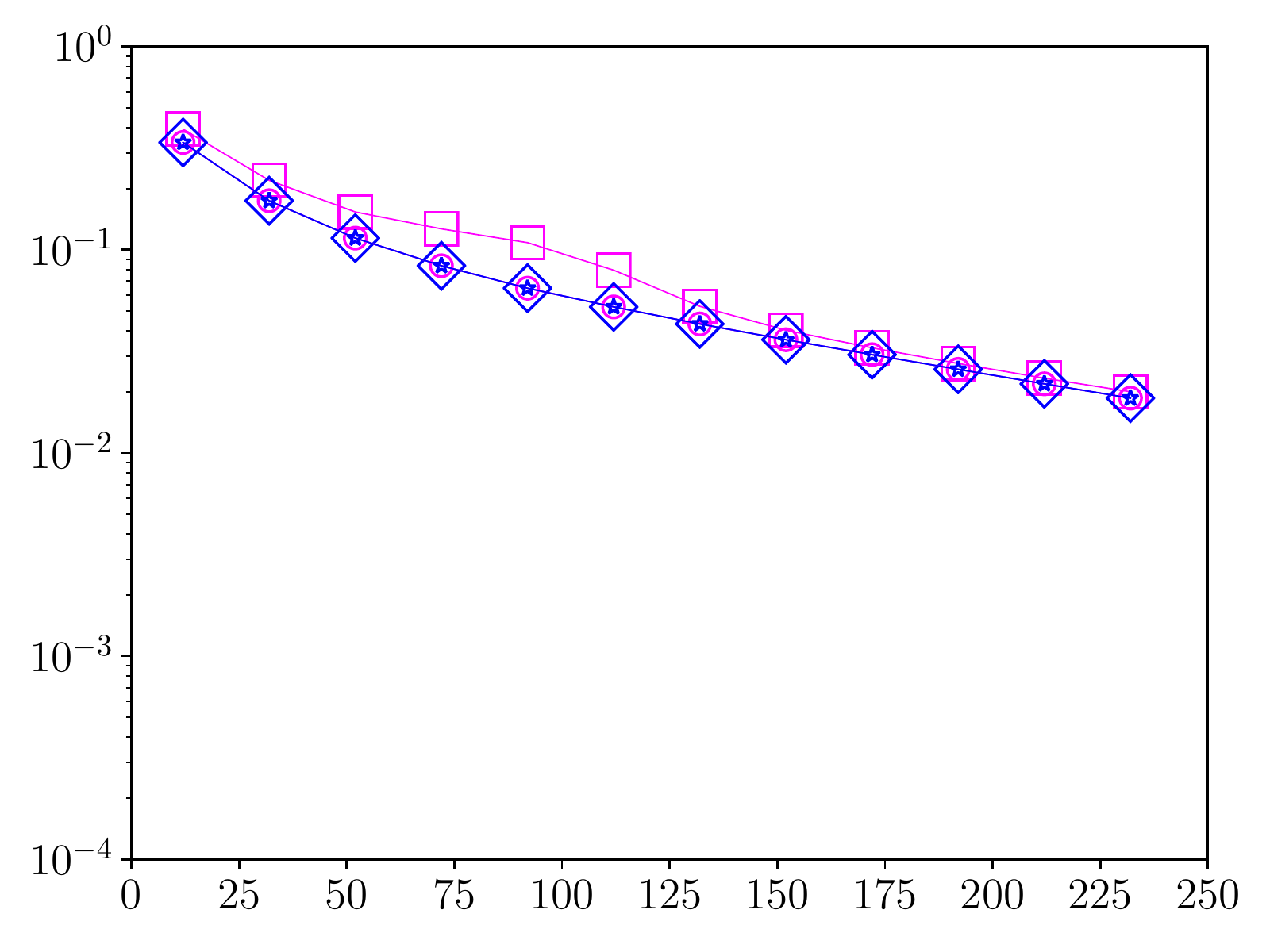}
\end{subfigure}
%
\begin{subfigure}[t]{0.33\textwidth}
\includegraphics[trim={0 1.1cm 0 0},clip,width=1.0\linewidth]{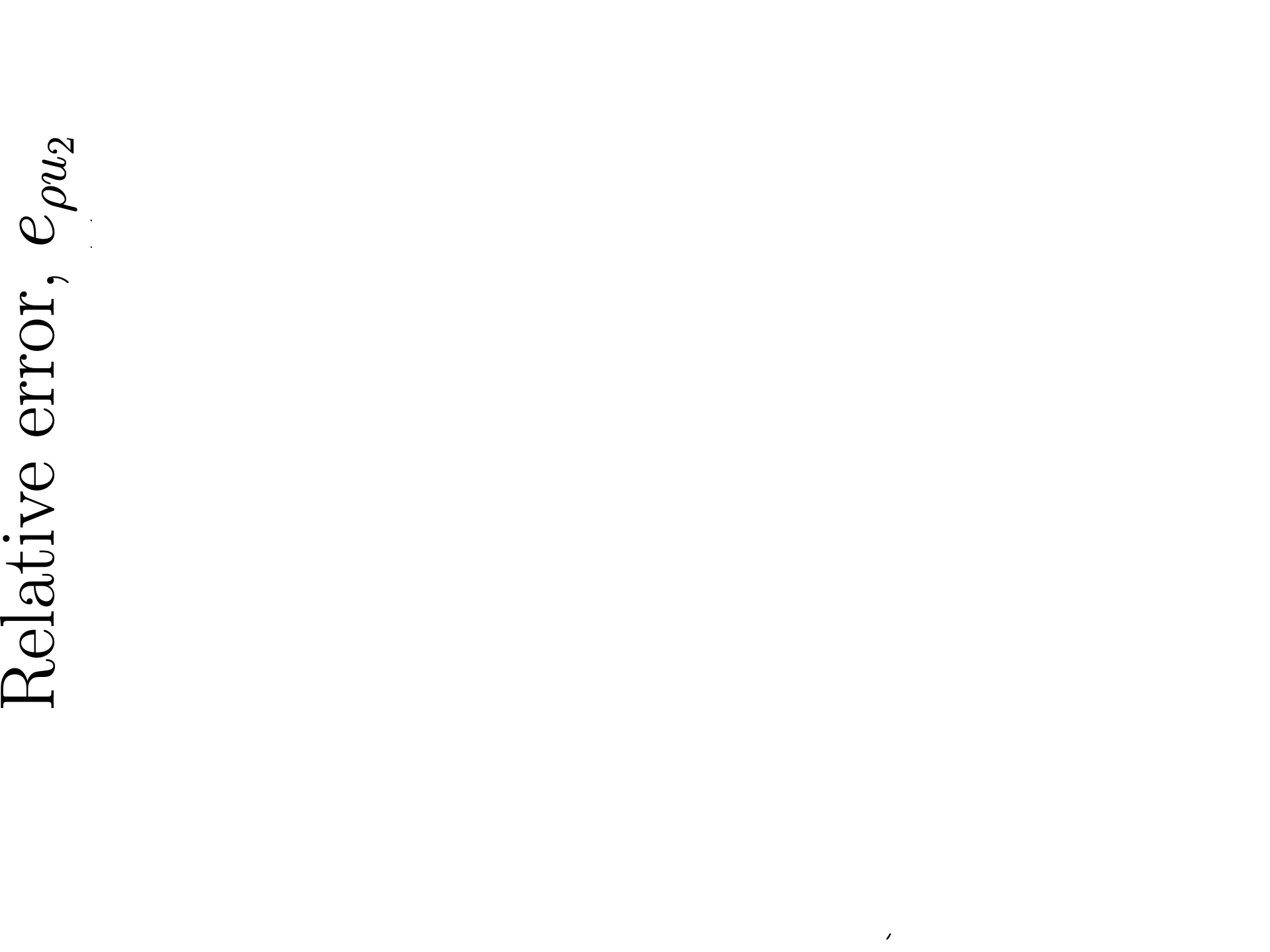}
\end{subfigure}
\begin{subfigure}[t]{0.31\textwidth}
\includegraphics[trim={0cm 0.5cm 0 0},clip,width=1.0\linewidth]{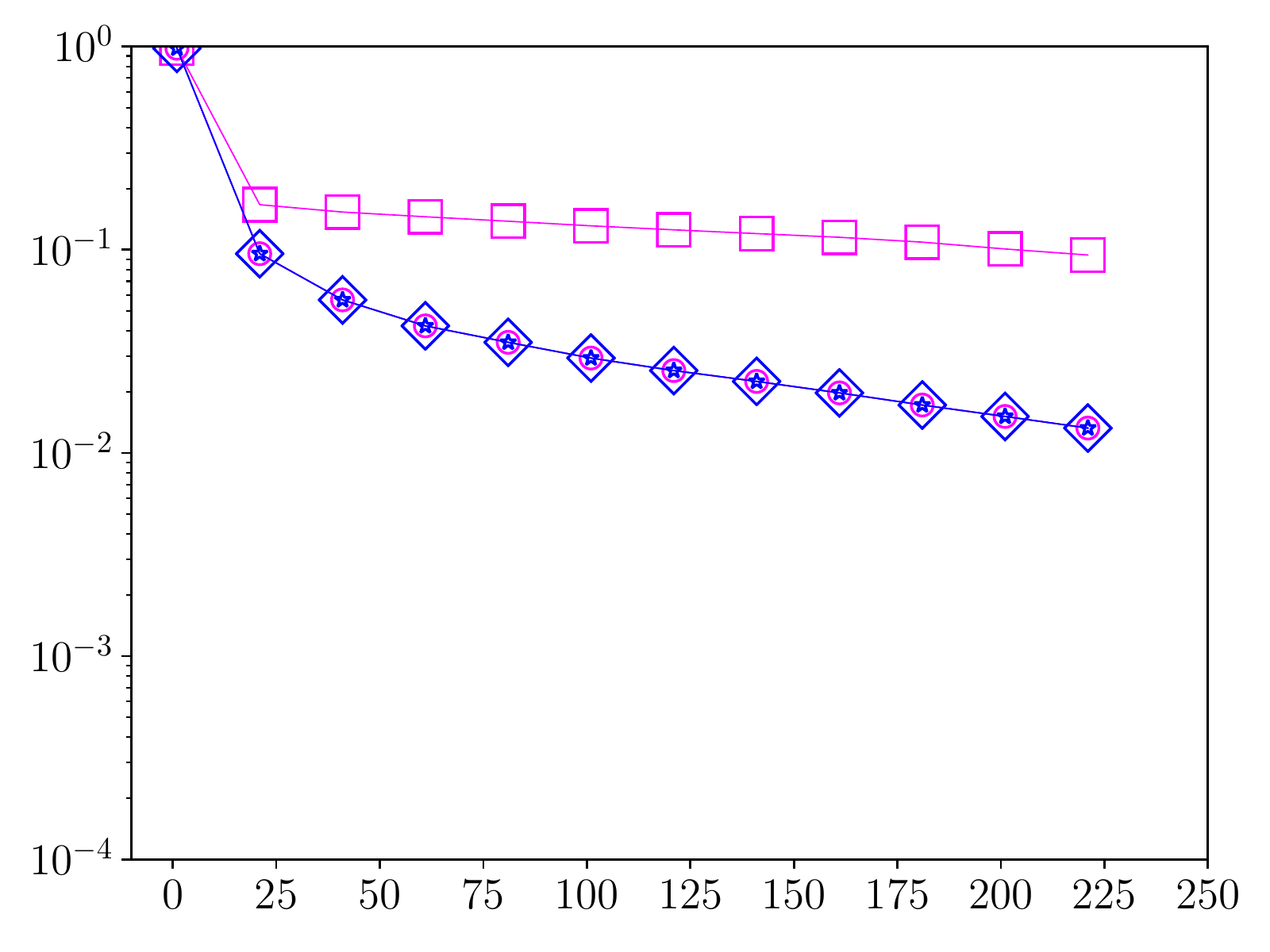}
\end{subfigure}
\begin{subfigure}[t]{0.31\textwidth}
\includegraphics[trim={0cm 0.5cm 0 0},clip,width=1.0\linewidth]{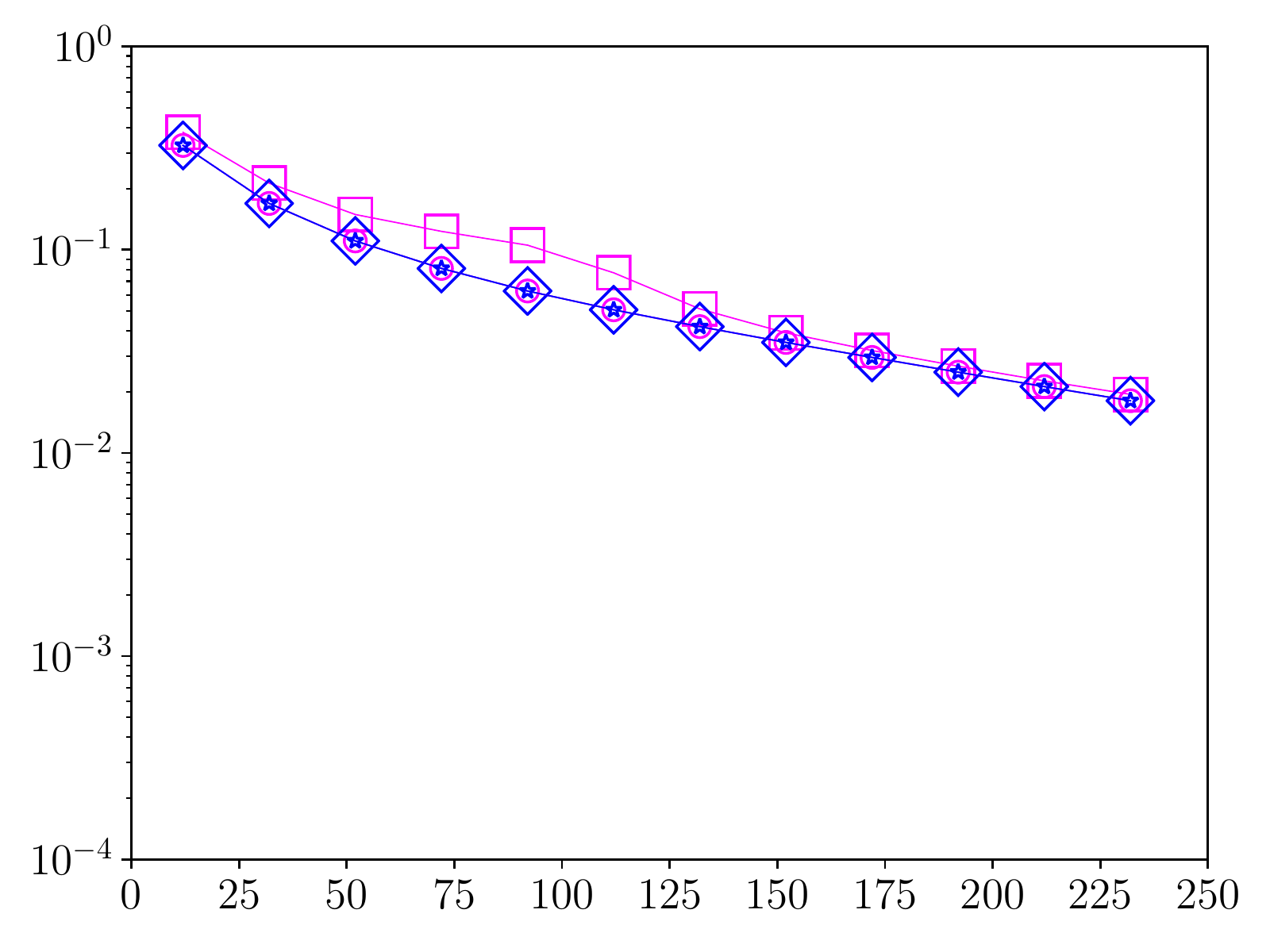}
\end{subfigure}
%
\begin{subfigure}[t]{0.31\textwidth}
\includegraphics[trim={0 0 0 -0.1cm},clip,width=1.0\linewidth]{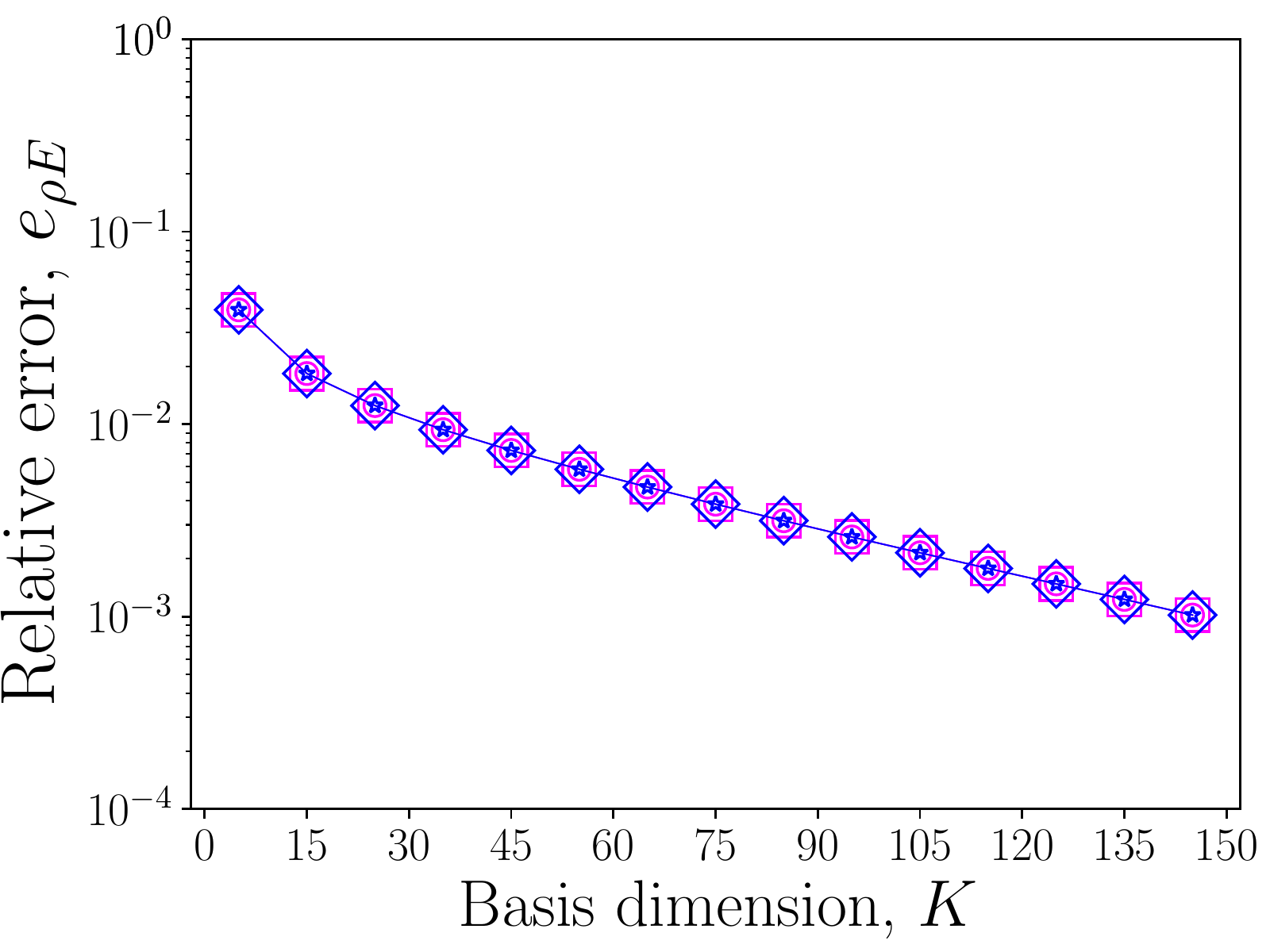}
\end{subfigure}
\begin{subfigure}[t]{0.31\textwidth}
\includegraphics[trim={0cm 0 0.4cm 0},clip,width=1.0\linewidth]{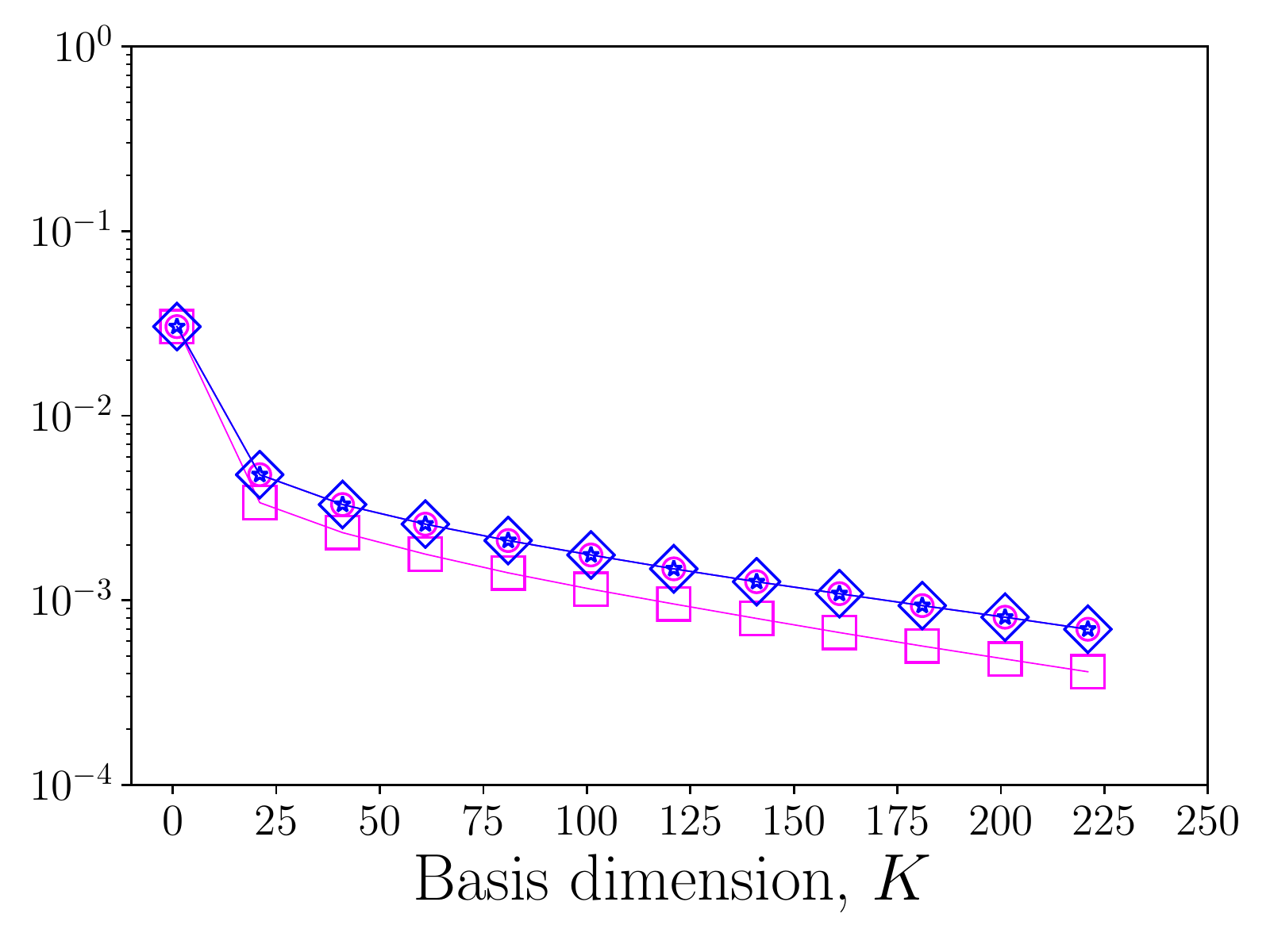}
\end{subfigure}
\begin{subfigure}[t]{0.31\textwidth}
\includegraphics[trim={0cm 0 0.4cm 0},clip,width=1.0\linewidth]{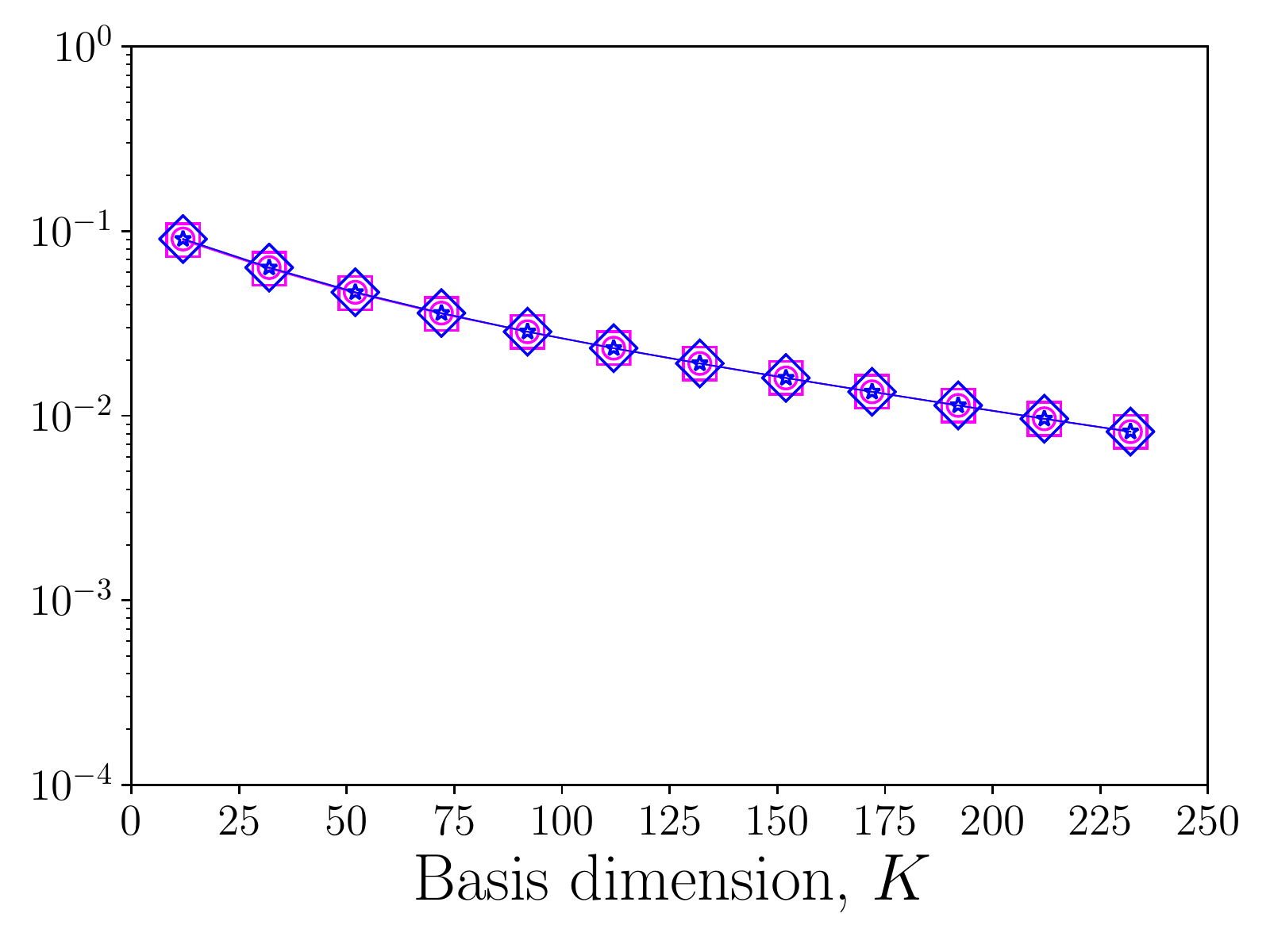}
\end{subfigure}
\caption{Relative errors associated with projection of the snapshot data on to the trial subspace for the Sod 1D problem (left column), Kelvin Helmholtz problem (center column),
and homogeneous turbulence problem (right column). Results are shown for the mass (top row),
$x_1$ momentum (second row), $x_2$ momentum (where applicable, third row),
and energy (last row) state variables as a function of basis dimension.
Results are shown for both the dimensional and non dimensional data.}
\label{fig:pod_l2_fig}
\end{center}
\end{figure}

We first study POD in the various inner products. Figure~\ref{fig:pod_l2_fig} presents the error associated with projection of the FOM training data onto the trial subspace as a function of basis dimension for the various inner products examined in Section~\ref{sec:podInnerProduct}. We note that projection is done in the same inner product as which the corresponding basis is constructed. Results are presented for the inner products applied to the dimensional and non-dimensional FOM configurations. In all cases we observe that the $L^2(\Omega)$ inner product results in different bases when applied to the dimensional data vs. non-dimensional data. This difference is very clear for the Kelvin Helmholtz problem (Figure~\ref{fig:pod_l2_fig}, center column) as well as for the $x_1$ and $x_2$ momentum fields for the homogeneous turbulence problem (Figure~\ref{fig:pod_l2_fig}, right column). The discrepancy is very minor for the Sod problem (Figure~\ref{fig:pod_l2_fig}, left column) and is not visibly obvious, but is still present. For the Kelvin Helmholtz problem we observe that, when compared to the non-dimensional $L^2_*(\Omega)$ inner product, the $L^2(\Omega)$ inner product applied to the dimensional data results in bases that have a poor approximation power for the density, $x_1$, and $x_2$ momentum fields, but a good approximation power for the total energy field. This result is intuitive as the total energy snapshots are orders of magnitude larger than the density and velocity snapshots and hence dominate the $L^2(\Omega)$ inner product. We next observe that, as expected, the non-dimensional $L^2_*(\Omega)$ inner product results in bases with equivalent approximation power when applied to the dimensional and non-dimensional data. Overall, we observe that the non-dimensional $L_*^2(\Omega)$ inner product is slightly ``better" for bases construction.

\begin{figure}
\begin{center}
\begin{subfigure}[t]{0.8\textwidth}
\includegraphics[trim={0 7cm 0 0},clip,width=1.0\linewidth]{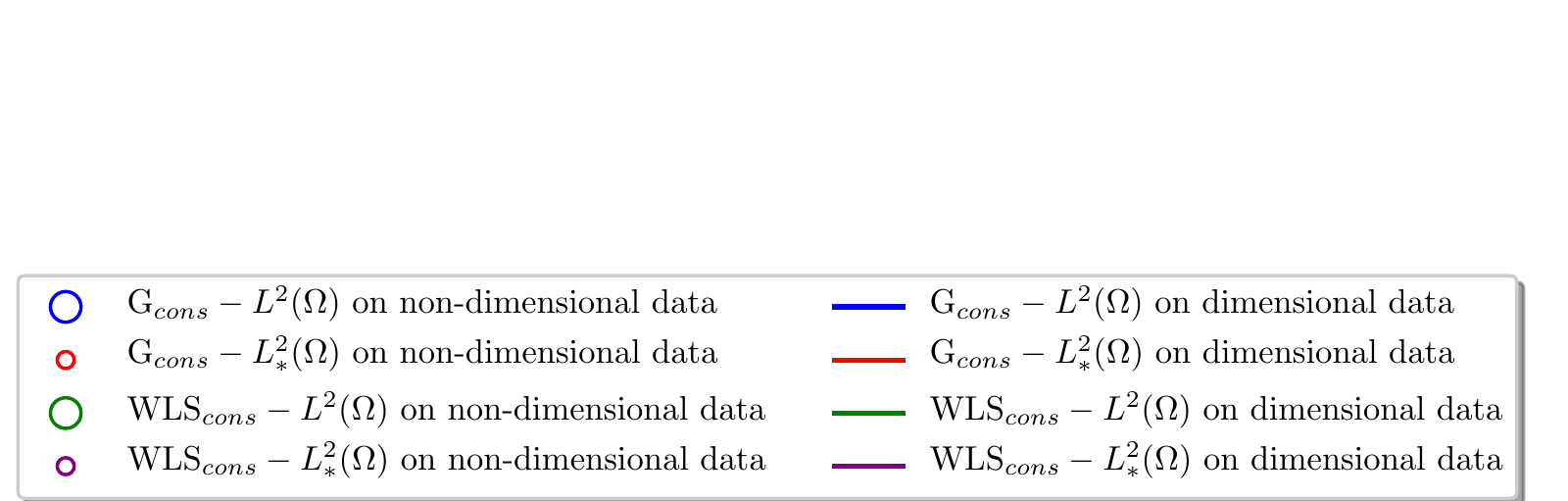}
\end{subfigure}

\begin{subfigure}[t]{0.48\textwidth}
\includegraphics[clip,width=1.0\linewidth]{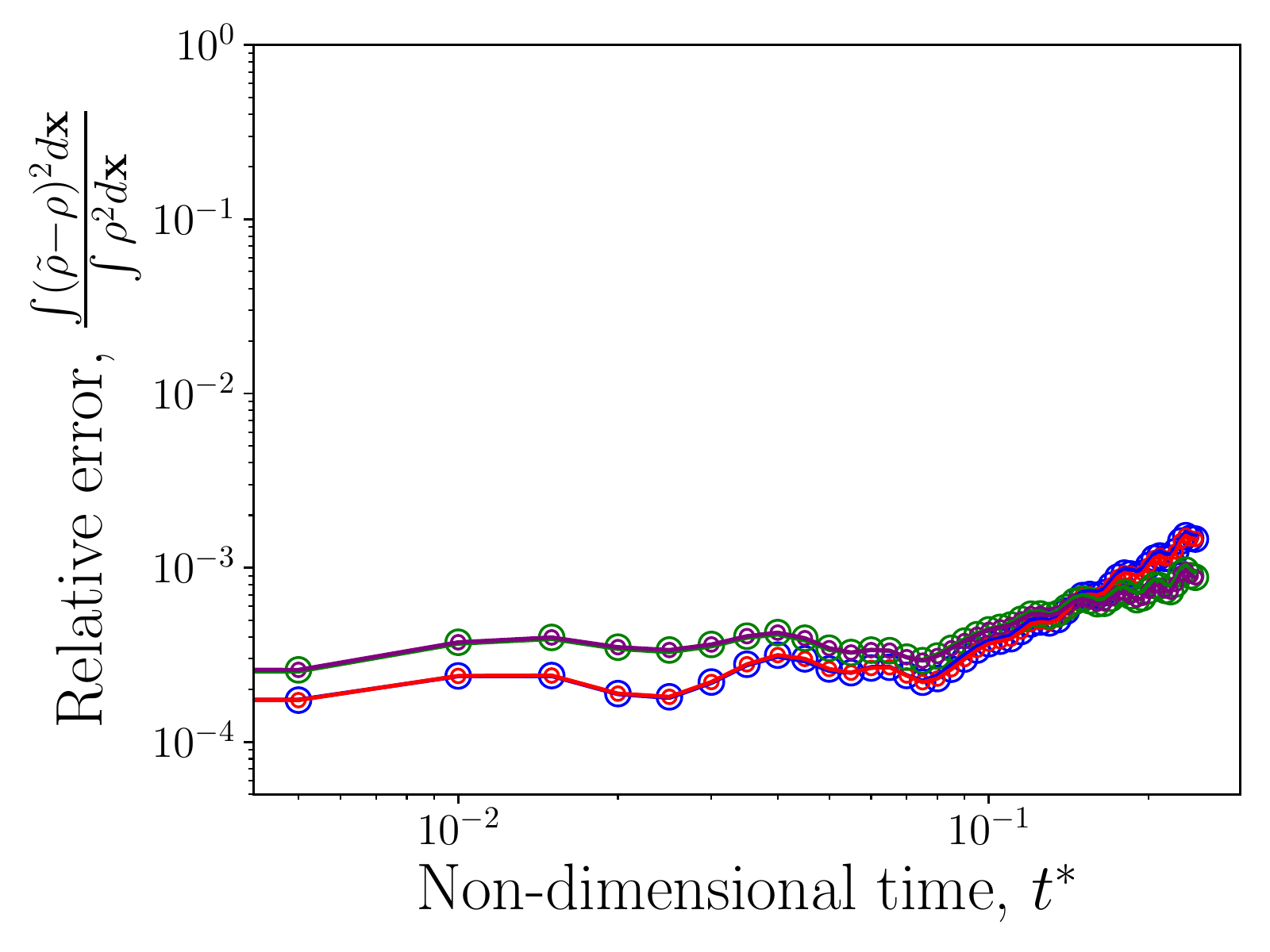}
\caption{Sod 1D problem}
\label{fig:l2_rom_summary_a}
\end{subfigure}
\begin{subfigure}[t]{0.48\textwidth}
\includegraphics[clip,width=1.0\linewidth]{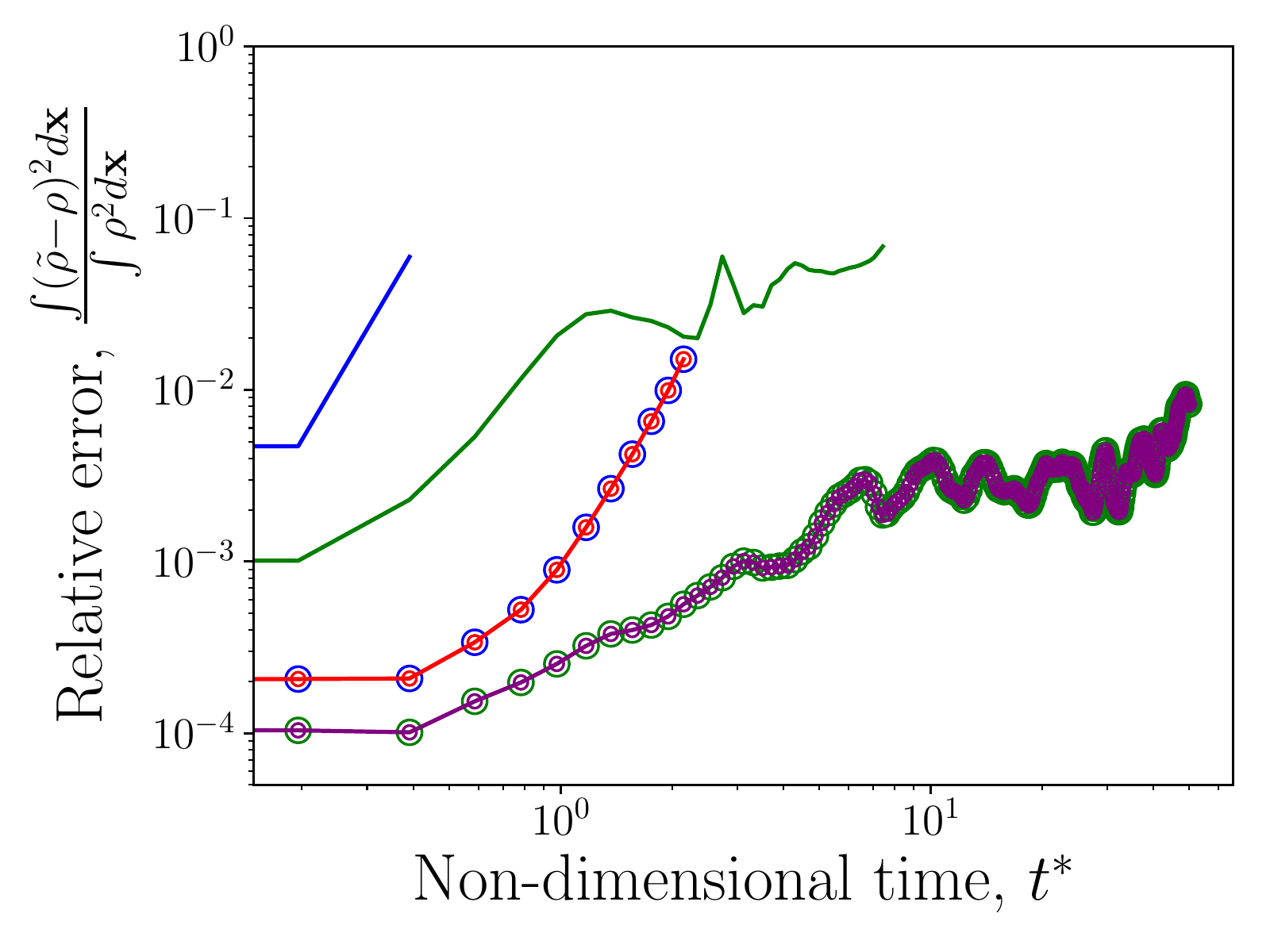}
\caption{Kelvin Helmholtz problem}
\label{fig:l2_rom_summary_b}
\end{subfigure}
\begin{subfigure}[t]{0.48\textwidth}
\includegraphics[clip,width=1.0\linewidth]{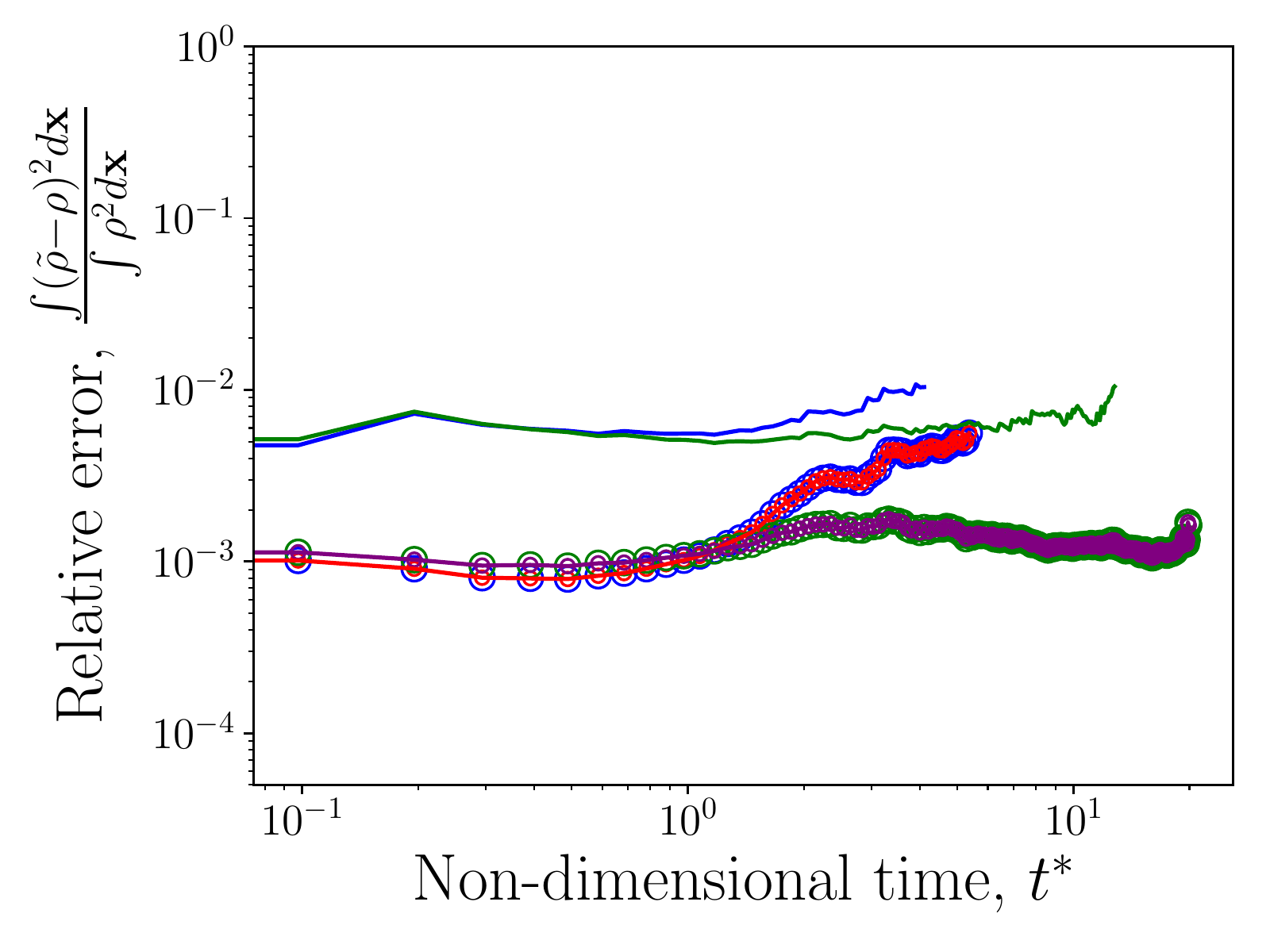}
\caption{Homogeneous turbulence problem}
\label{fig:l2_rom_summary_c}
\end{subfigure}
\begin{subfigure}[t]{0.48\textwidth}
\includegraphics[clip,width=1.0\linewidth]{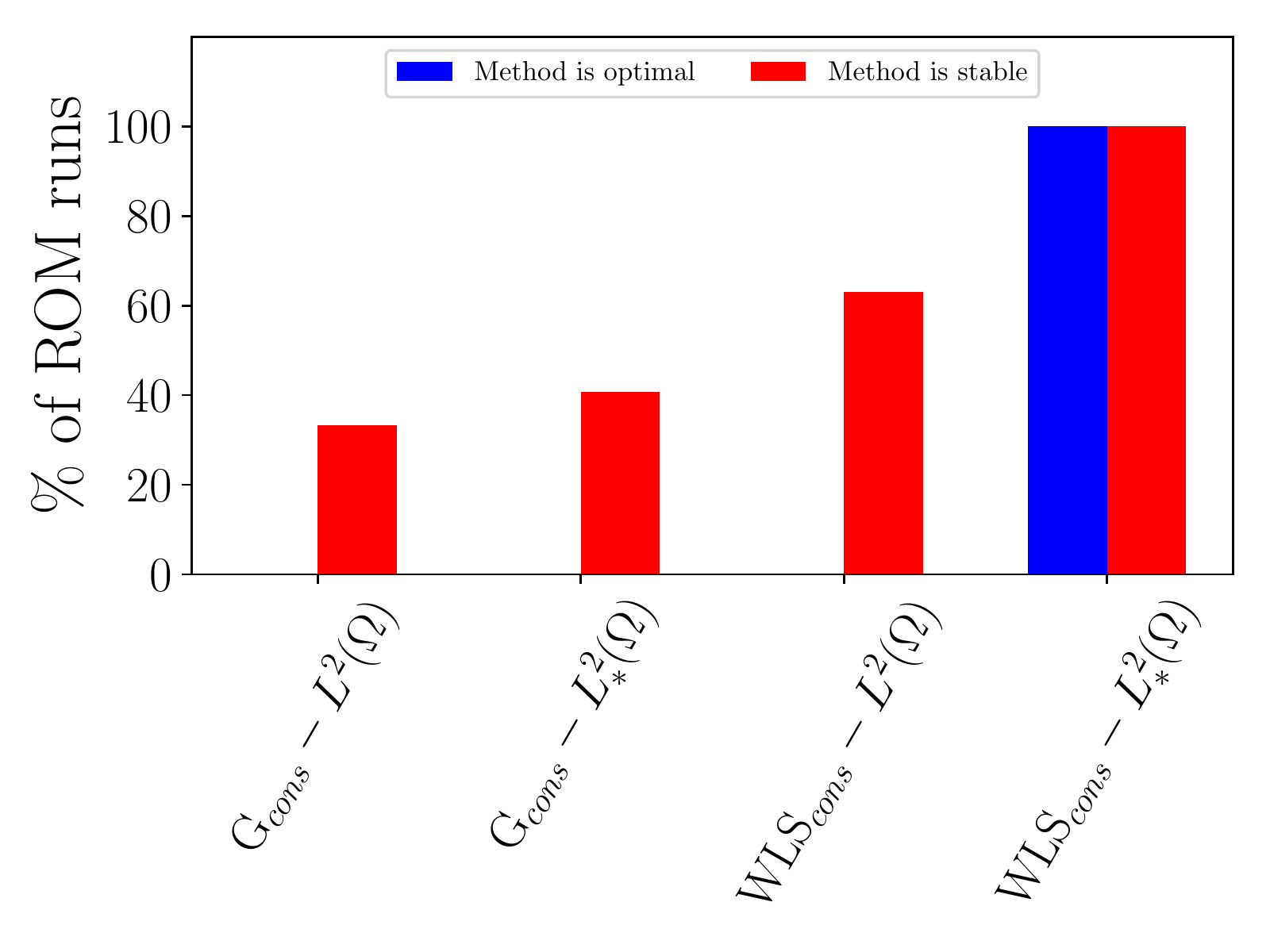}
\caption{ROM summary}
\label{fig:l2_rom_summary_d}
\end{subfigure}
\caption{Comparison of ROM results obtained with dimensional $\LTwoSymbol$ and non-dimensional $\NonDimensionalLTwoSymbol$ inner products. Figures~\ref{fig:l2_rom_summary_a}-\ref{fig:l2_rom_summary_c} depict errors in the density field as a function of time for the Sod 1D, Kelvin Helmholtz, and homogeneous turbulence problems, respectively. Figure~\ref{fig:l2_rom_summary_d} summarizes the percentage of time a ROM formulation was optimal and the percentage of time a ROM formulation was stable.}
\label{fig:l2_rom_summary}
\end{center}
\end{figure}

Next we consider the results of the \GalerkinLTwoConservedRomName, \GalerkinNonDimensionalLTwoConservedRomName, \WLSLTwoConservedRomName, and \WLSNonDimensionalLTwoConservedRomName\ ROMs. Figure~\ref{fig:l2_rom_summary} provides a summary of the performance of these ROMs on the three exemplars. First, Figures~\ref{fig:l2_rom_summary_a}-\ref{fig:l2_rom_summary_c} depict the error for the density field predicted by the various ROMs at dimensions of $K=50, 50,$ and $100$ as a function of time for the Sod 1D, Kelvin Helmholtz, and homogeneous turbulence problems, respectively. Results are shown for both the dimensional and non-dimensional configurations. Next, Figure~\ref{fig:l2_rom_summary_d} summarizes the overall ROM performance by tabulating (1) the number of times a ROM formulation gave solutions with the lowest errors for a state variable (among the \GalerkinLTwoConservedRomName, \GalerkinNonDimensionalLTwoConservedRomName, \WLSLTwoConservedRomName, and \WLSNonDimensionalLTwoConservedRomName\ ROMs) and (2) the percentage of the time a ROM formulation was stable. The results in Figure~\ref{fig:l2_rom_summary_d} are computed by aggregating across all ROM runs on the dimensional configurations. We make the following observations:
\begin{itemize}

\item (Figures~\ref{fig:l2_rom_summary_a}-\ref{fig:l2_rom_summary_c}.) The \GalerkinNonDimensionalLTwoConservedRomName\ and \WLSNonDimensionalLTwoConservedRomName\ ROMs, which are based on non-dimensional inner products, yield the same results when deployed on the dimensional and non-dimensional configurations. This result is expected and demonstrates that these methods are dimensionally consistent.

\item (Figures~\ref{fig:l2_rom_summary_a}-\ref{fig:l2_rom_summary_c}.) The \GalerkinLTwoConservedRomName\ and \WLSLTwoConservedRomName\ ROMs, which are \textit{not} based on dimensionally-consistent inner products, yield \textit{different} results for the dimensional and non-dimensional configurations. This result is very minor for the Sod 1D problem, but is very evident for the Kelvin Helmholtz and homogeneous turbulence examples. This result is again expected, and demonstrates that these ROM formulations are dimensionally inconsistent.

\item (Figures~\ref{fig:l2_rom_summary_a}-\ref{fig:l2_rom_summary_d}.) When applied to dimensional problems, the ROMs based on the non-dimensional $\NonDimensionalLTwoSymbol$ inner product consistently outperform their counterparts that are based on the standard dimensional $\LTwoSymbol$ inner product.

\item (Figure~\ref{fig:l2_rom_summary_d}.) The \GalerkinLTwoConservedRomName\ and \GalerkinNonDimensionalLTwoConservedRomName\ ROMs are rarely stable; the \GalerkinLTwoConservedRomName\ ROM was stable less than 40\% of the time, while the \GalerkinNonDimensionalLTwoConservedRomName\ ROM was stable approximately 40\% of the time.

\item (Figure~\ref{fig:l2_rom_summary_d}.) The \WLSNonDimensionalLTwoConservedRomName\ ROM was always stable and always resulted in solutions with the lowest errors for each conserved variable.
\end{itemize}

This section demonstrated that constructing ROMs on dimensionally-consistent inner products results in more stable ROMs that produce more accurate solutions. As a result we will no longer focus on the \GalerkinLTwoConservedRomName\ and \WLSLTwoConservedRomName\ ROMs in the remainder of this work as they are dimensionally inconsistent (which makes them theoretically unappealing) and, as just seen, result in inaccurate/unstable ROMs.

\section{Non-dimensional inner products vs physics-based inner product}\label{sec:results3}
The previous section demonstrated the importance of having a dimensionally-consistent inner product when working with dimensional problems. However, even though it is dimensionally consistent, the $L^2_*(\Omega)$ inner product is arguably not a ``physics-based" inner product as it simply results in a (weighted) sum of the different state variables and carries no physical meaning. In this section, we investigate the performance of ROMs built on the physics-based ``entropy" inner products introduced in Section~\ref{sec:podInnerProduct} and compare them to ROMs built using the non-dimensional $L^2_*(\Omega)$ inner products. We emphasize that these entropy inner products are dimensionally consistent.

\begin{figure}[!t]
\begin{center}
\begin{subfigure}[t]{0.65\textwidth}
\includegraphics[trim={0 0 0 0},clip,width=1.0\linewidth]{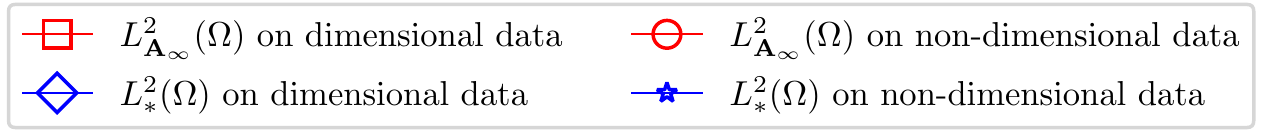}
\end{subfigure}\\
%
\begin{subfigure}[t]{0.33\textwidth}
\includegraphics[trim={0 1.1cm 0 0},clip,width=1.0\linewidth]{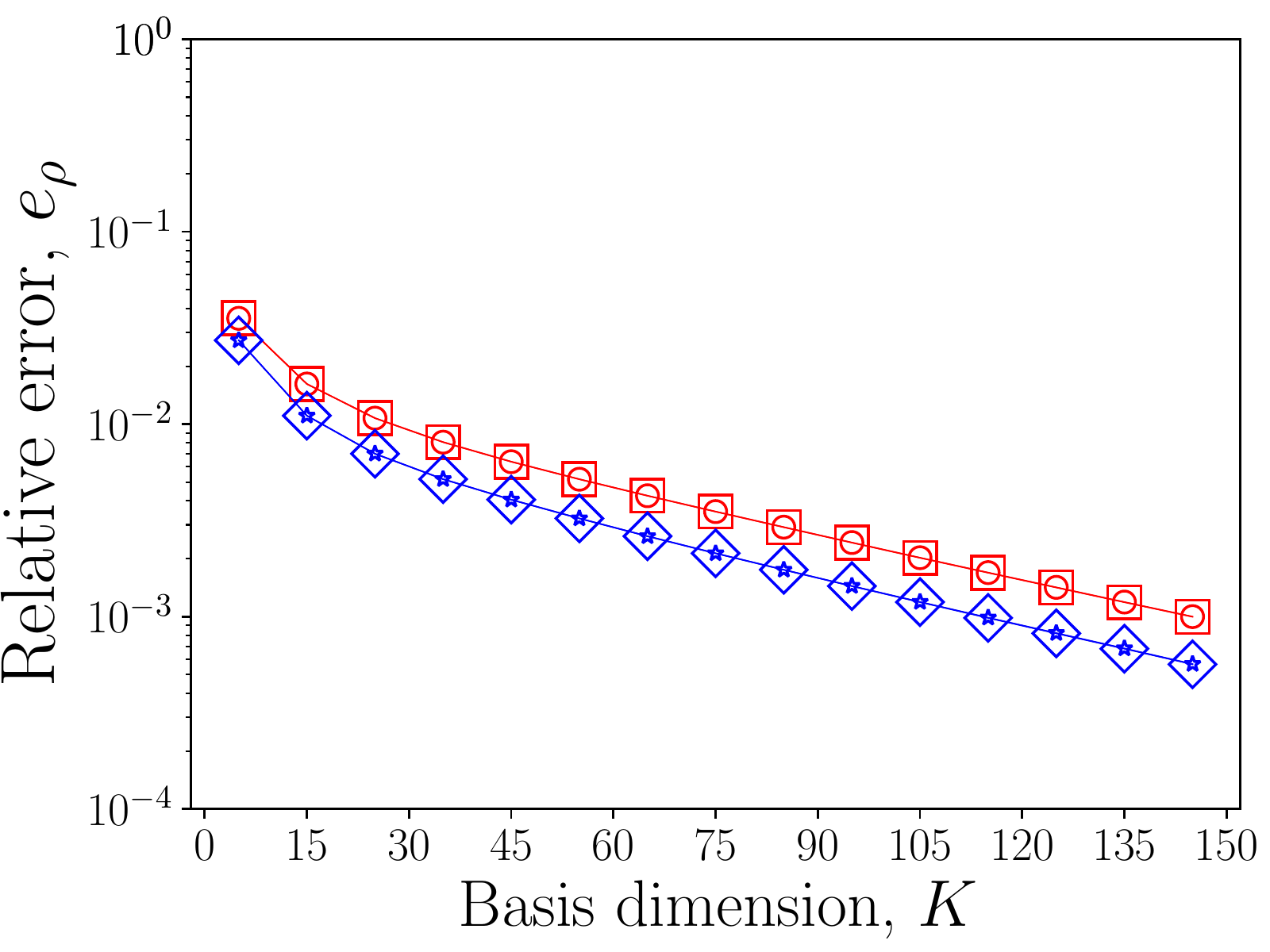}
\end{subfigure}
\begin{subfigure}[t]{0.31\textwidth}
\includegraphics[trim={0cm 0.5cm 0 0},clip,width=1.0\linewidth]{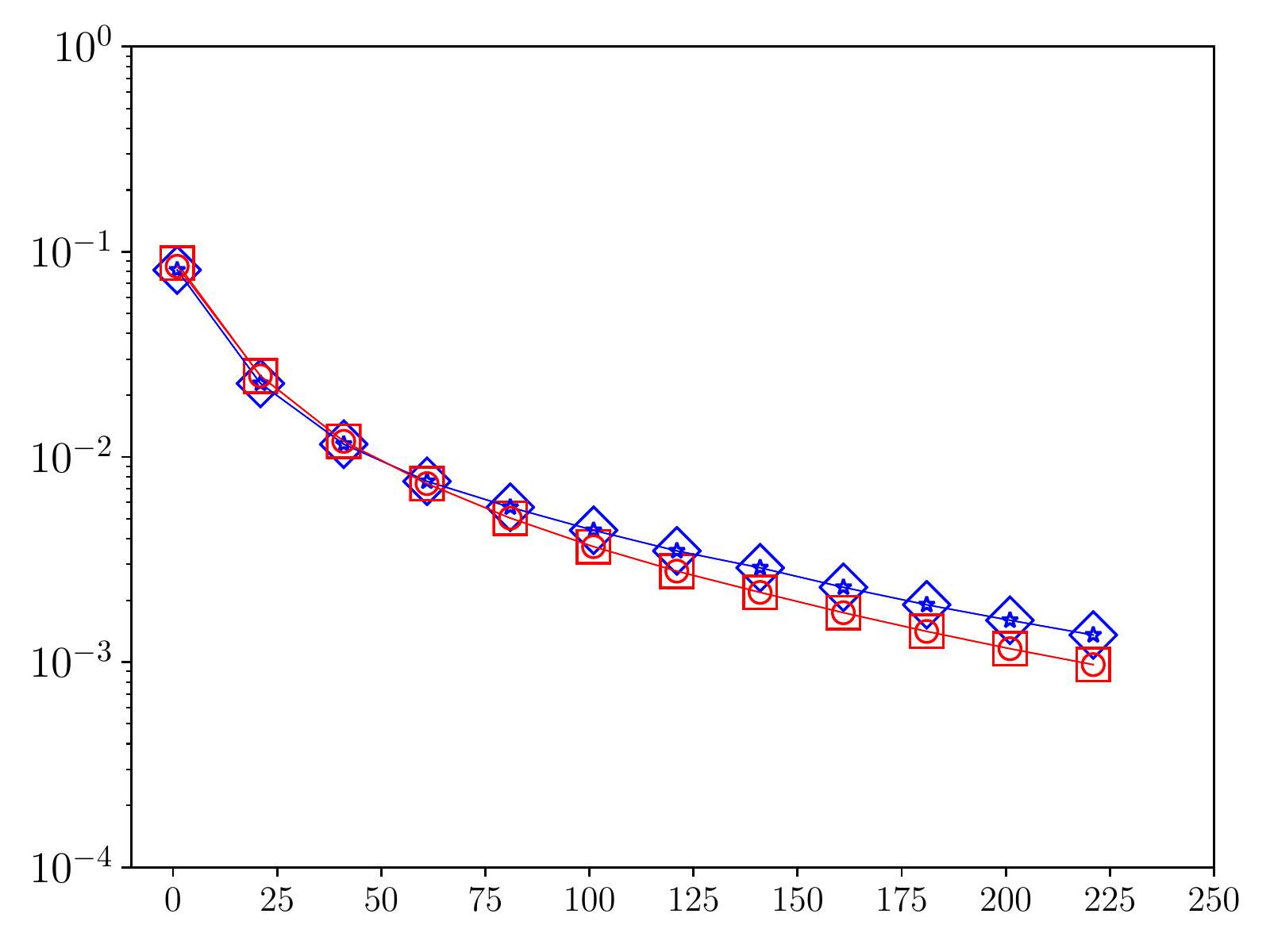}
\end{subfigure}
\begin{subfigure}[t]{0.31\textwidth}
\includegraphics[trim={0cm 0.5cm 0 0},clip,width=1.0\linewidth]{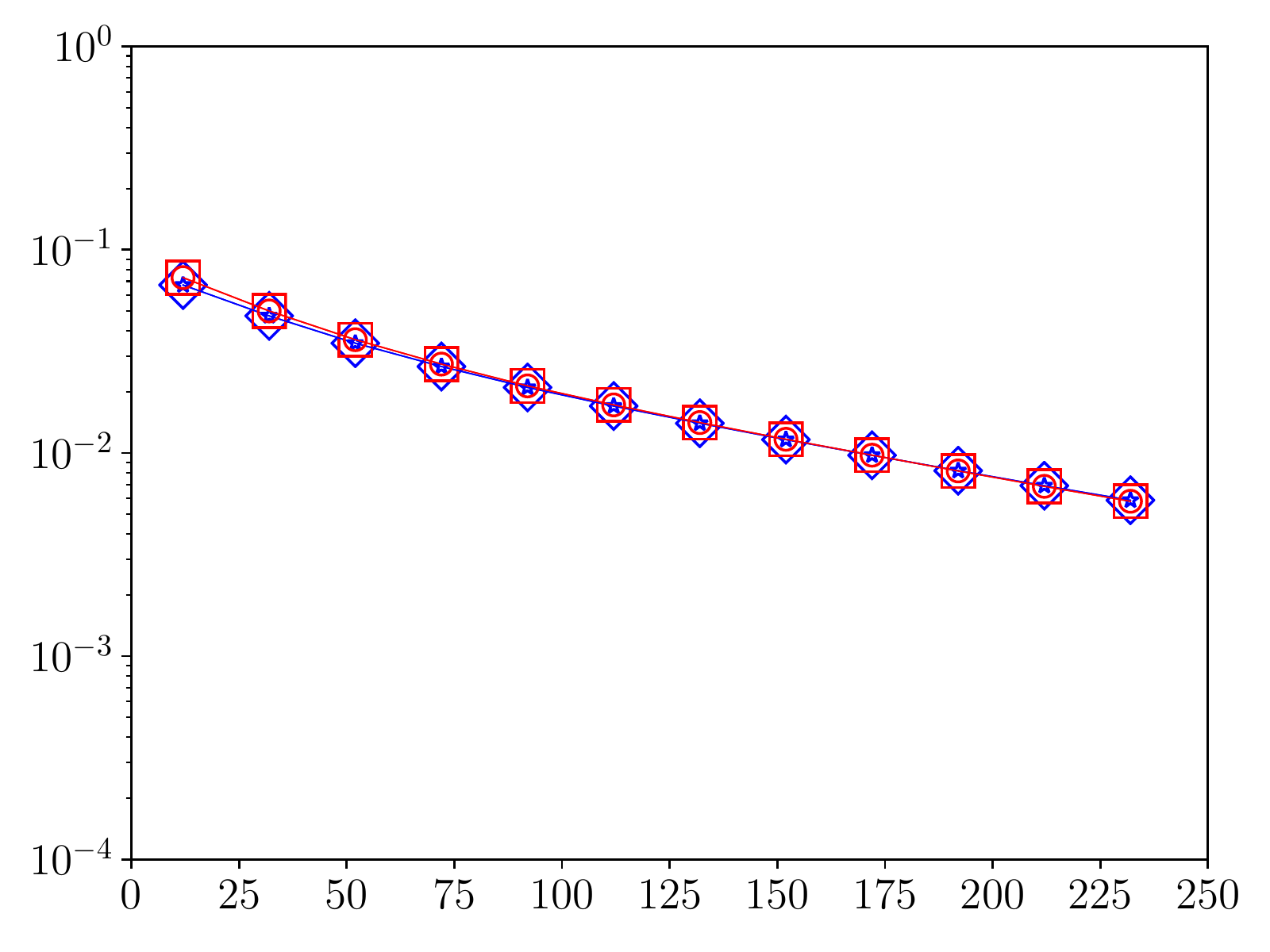}
\end{subfigure}
%
\begin{subfigure}[t]{0.33\textwidth}
\includegraphics[trim={0 1.1cm 0 0},clip,width=1.0\linewidth]{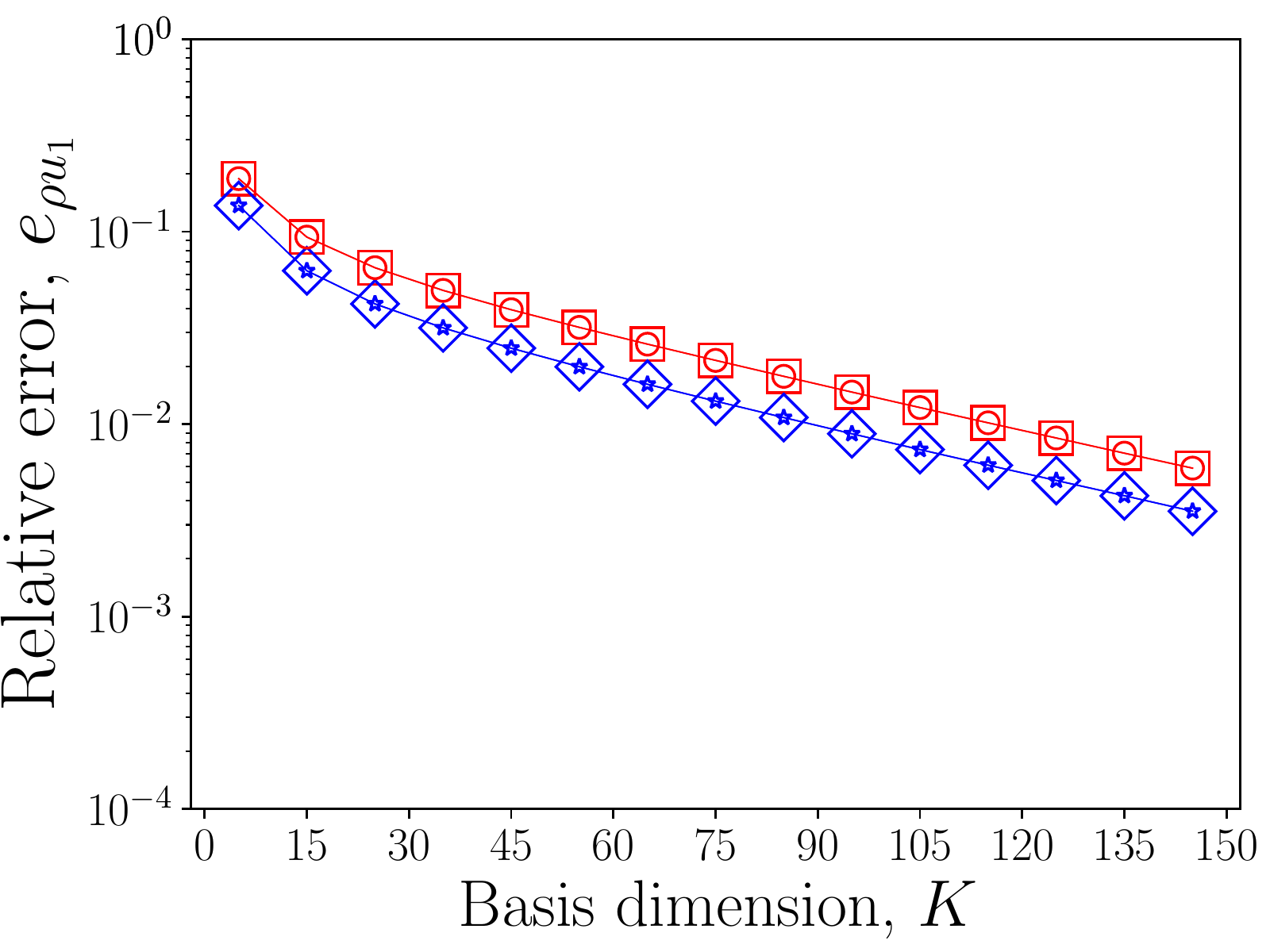}
\end{subfigure}
\begin{subfigure}[t]{0.31\textwidth}
\includegraphics[trim={0cm 0.5cm 0 0},clip,width=1.0\linewidth]{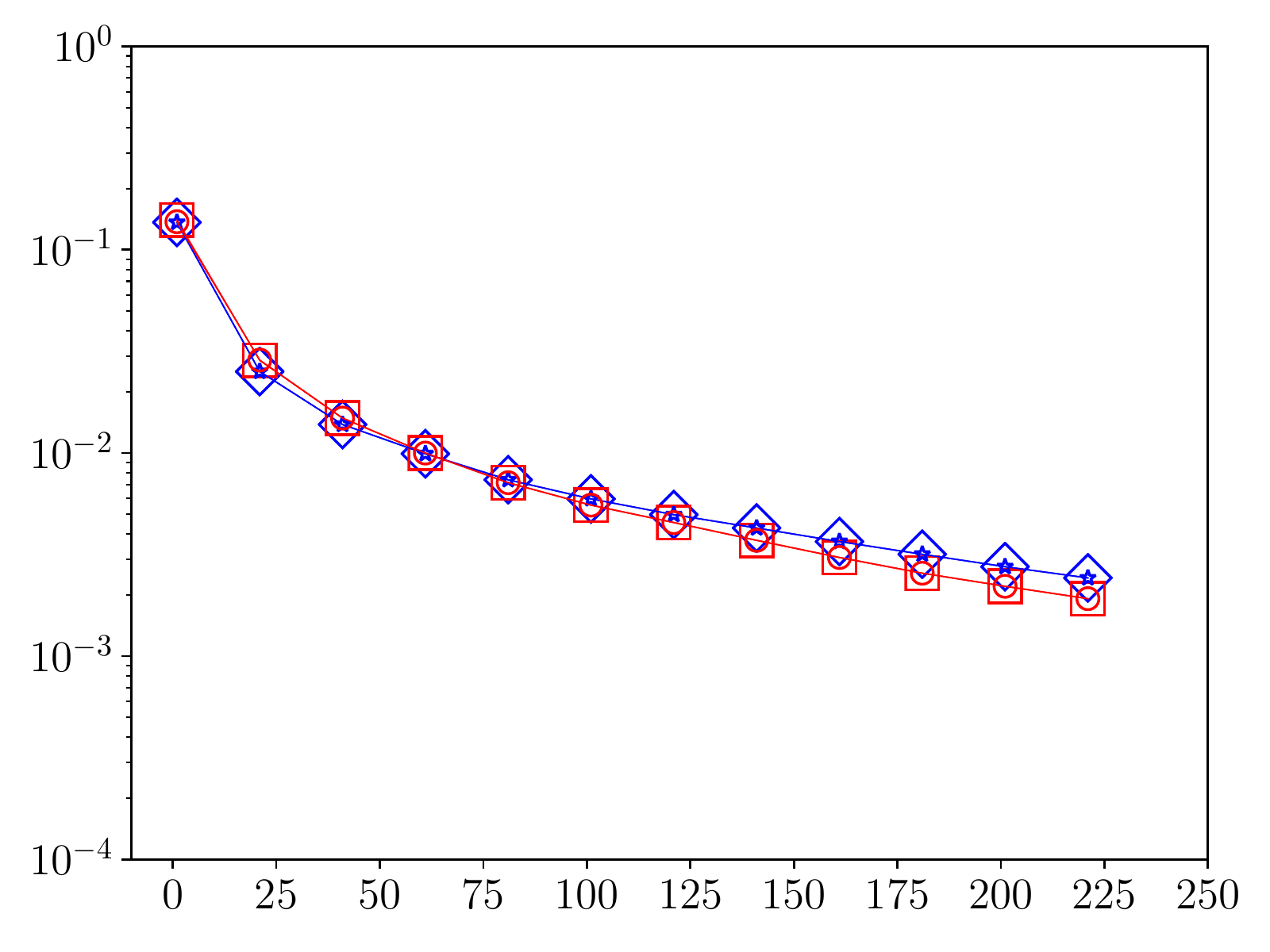}
\end{subfigure}
\begin{subfigure}[t]{0.31\textwidth}
\includegraphics[trim={0cm 0.5cm 0 0},clip,width=1.0\linewidth]{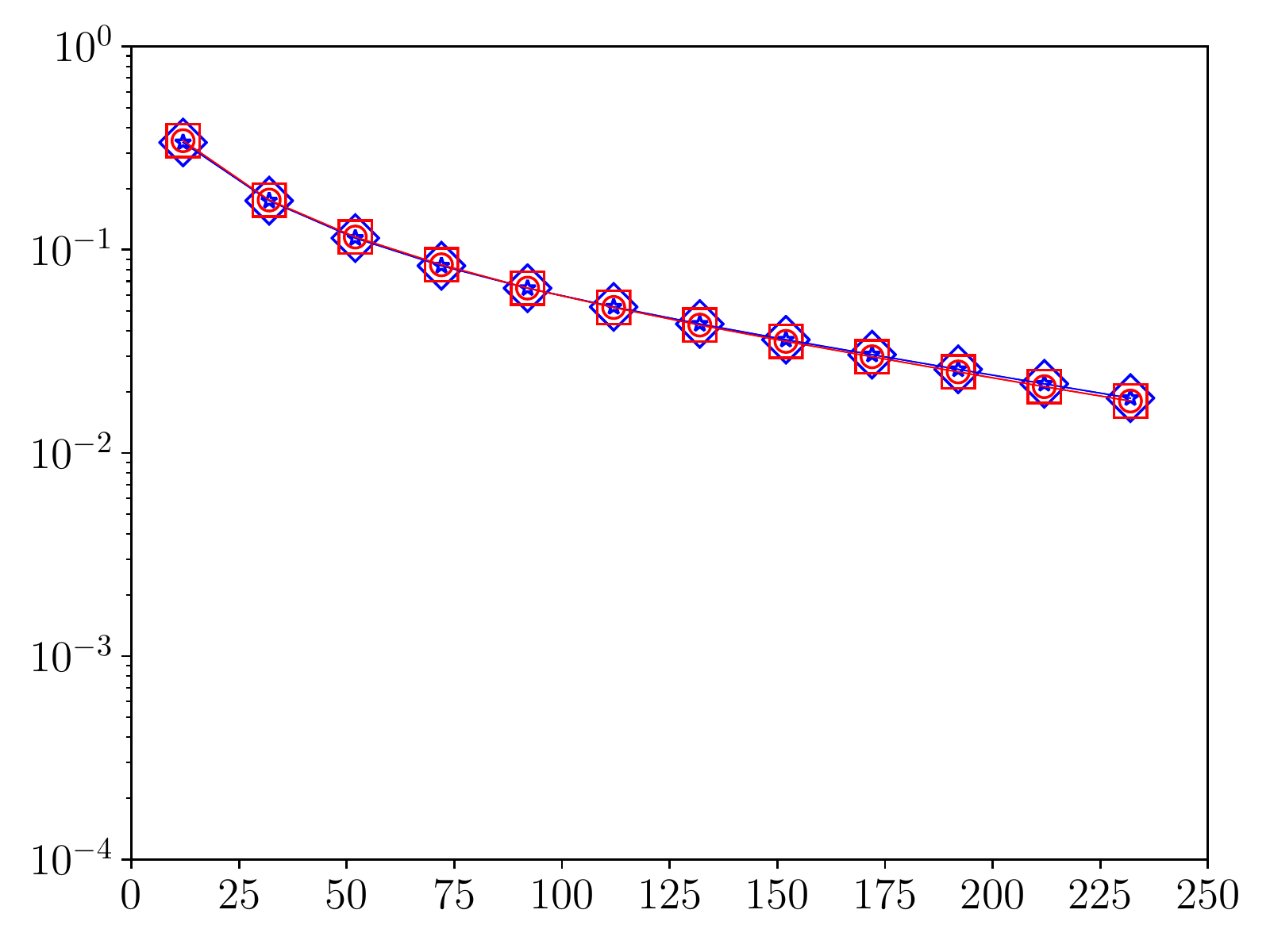}
\end{subfigure}
%
\begin{subfigure}[t]{0.33\textwidth}
\includegraphics[trim={0 1.1cm 0 0},clip,width=1.0\linewidth]{figures/sod1d/fig5_dummy.png}
\end{subfigure}
\begin{subfigure}[t]{0.31\textwidth}
\includegraphics[trim={0cm 0.5cm 0 0},clip,width=1.0\linewidth]{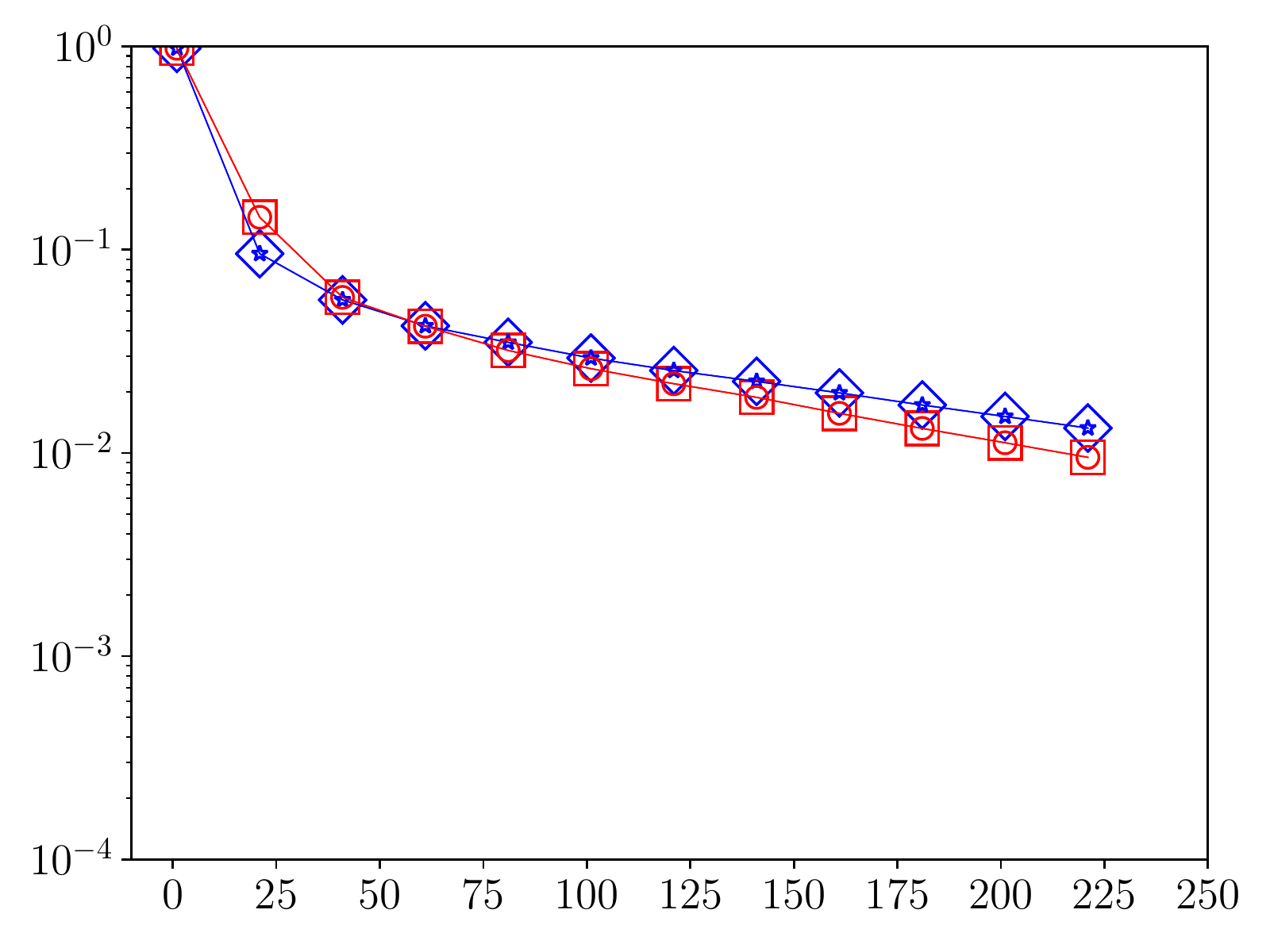}
\end{subfigure}
\begin{subfigure}[t]{0.31\textwidth}
\includegraphics[trim={0cm 0.5cm 0 0},clip,width=1.0\linewidth]{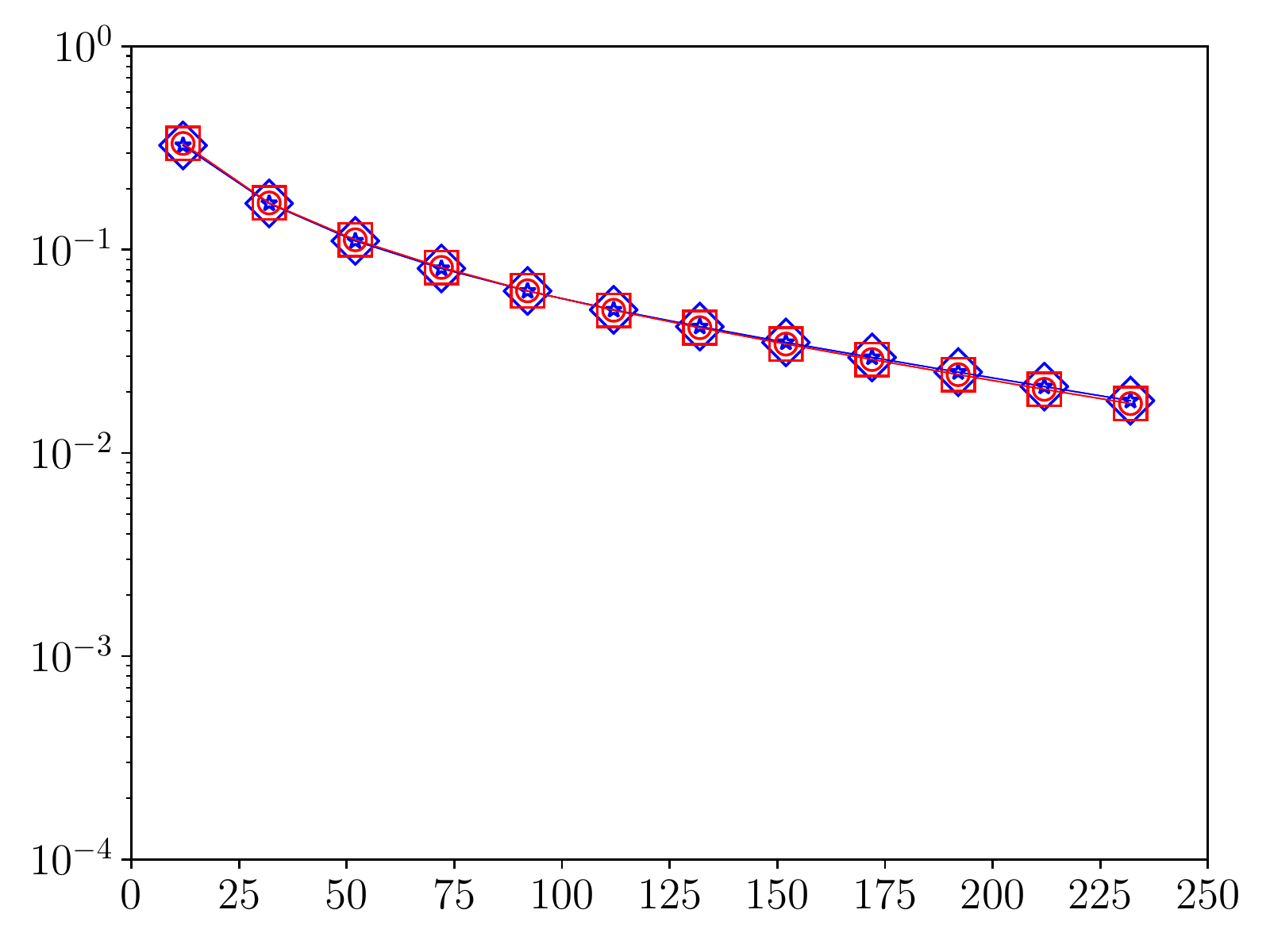}
\end{subfigure}
%
\begin{subfigure}[t]{0.33\textwidth}
\includegraphics[trim={0 0 0 -0.1},clip,width=1.0\linewidth]{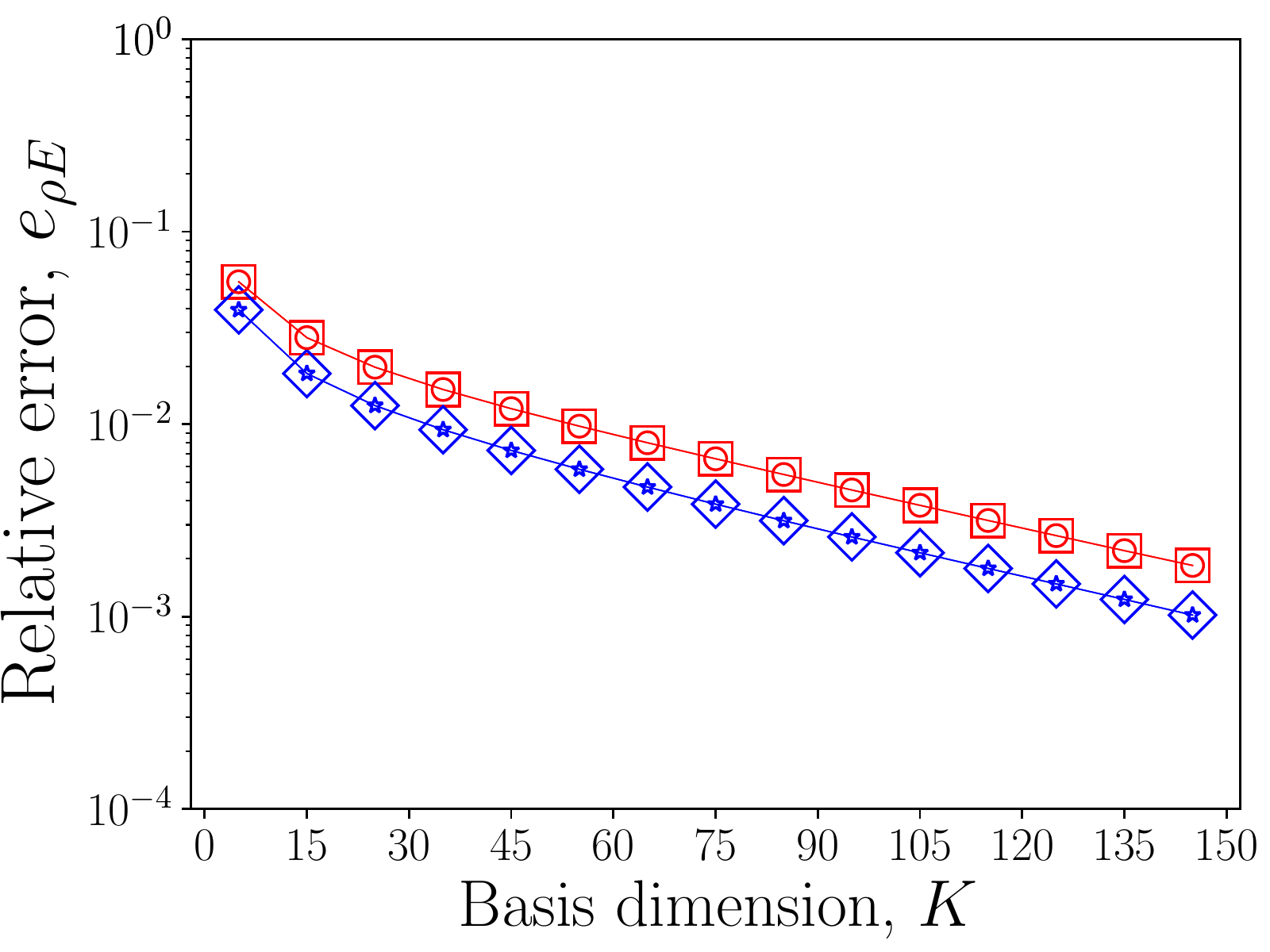}
\end{subfigure}
\begin{subfigure}[t]{0.31\textwidth}
\includegraphics[trim={0 0 0.4cm 0},clip,width=1.0\linewidth]{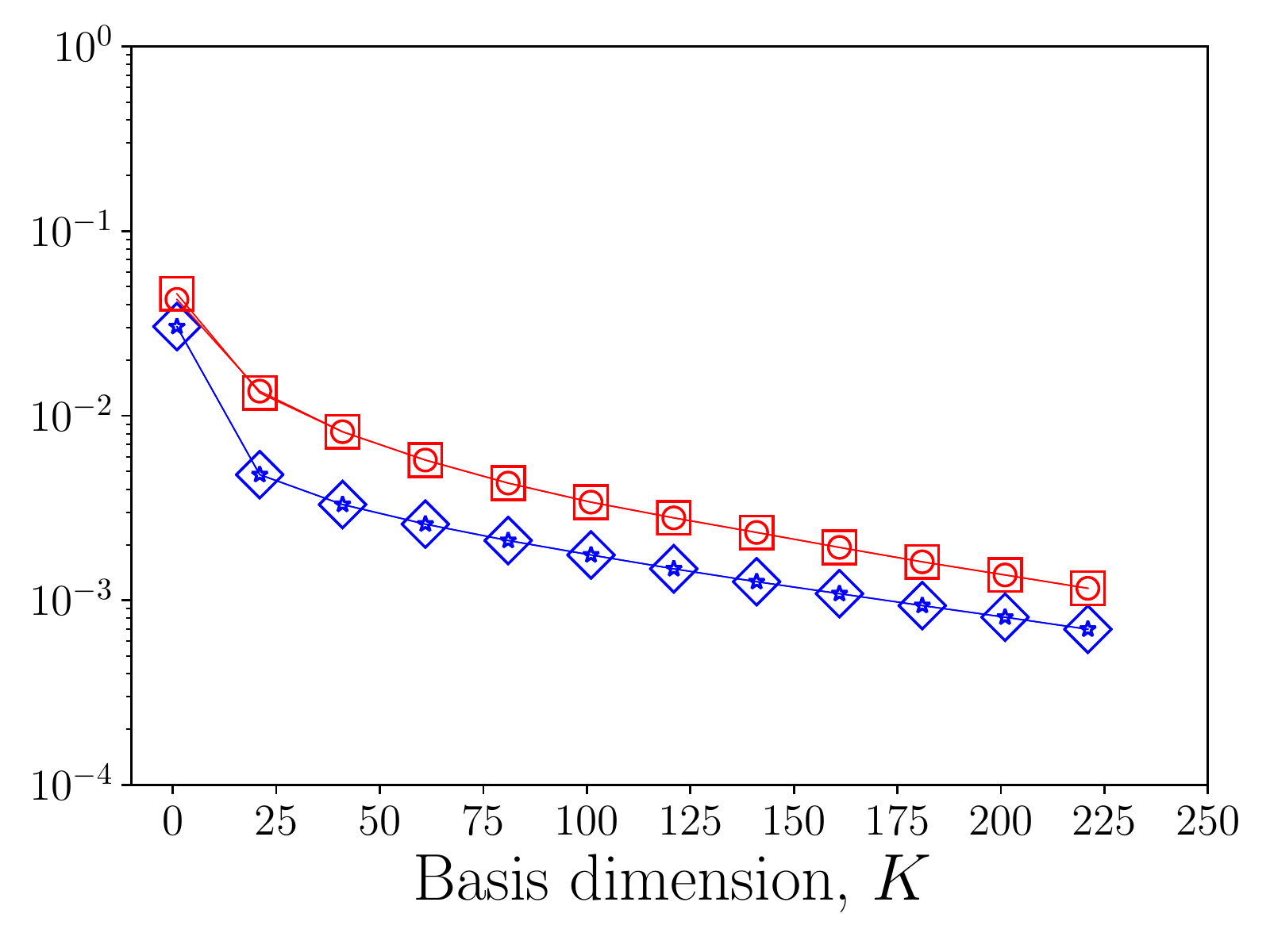}
\end{subfigure}
\begin{subfigure}[t]{0.31\textwidth}
\includegraphics[trim={0 0 0.4cm 0},clip,width=1.0\linewidth]{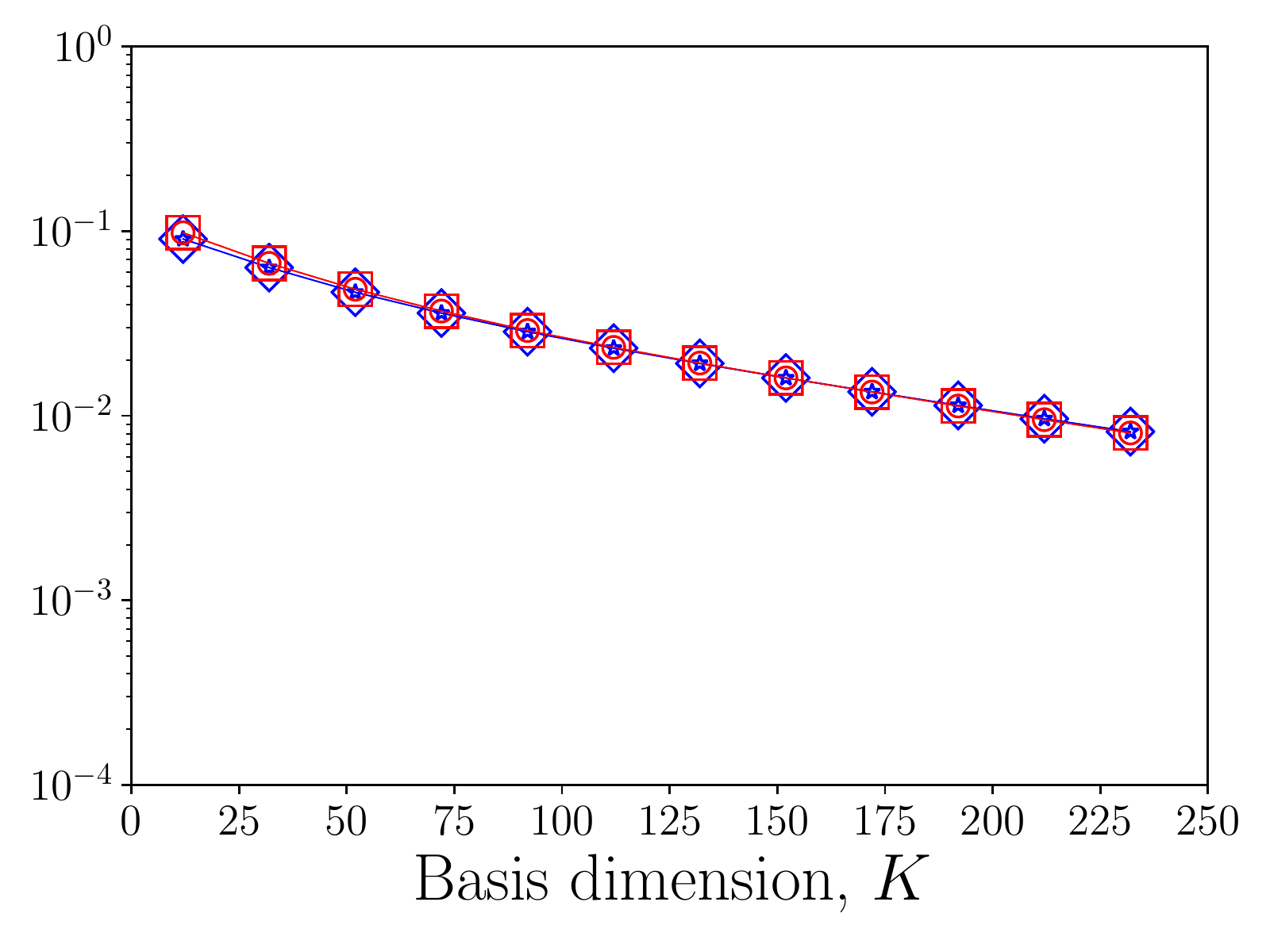}
\end{subfigure}
\caption{Relative errors associated with projection of the snapshot data onto the trial subspace in the $\EntropySymbol$ and $\NonDimensionalLTwoSymbol$ inner products
  for the Sod 1D problem (left column), Kelvin Helmholtz problem (center column),
  and homogeneous turbulence problem (right column).
  Results are shown for the mass (top row),
$x_1$ momentum (second row), $x_2$ momentum (where applicable, third row),
and energy (last row) state variables as a function of basis dimension.
Results are shown for both the dimensional and non dimensional data.
}
\label{fig:l2star_vs_entropy_pod_fig}
\end{center}
\end{figure}

We again first study POD in these different inner products. Figure~\ref{fig:l2star_vs_entropy_pod_fig} presents the projection error\footnote{Each projection is performed in the relevant inner product.} of the FOM data onto the trial subspace as a function of the basis dimension for the non-dimensional $\NonDimensionalLTwoSymbol$ and entropy $\EntropySymbol$ inner products examined in Section~\ref{sec:podInnerProduct}; we note that we don't consider bases constructed from the $\EntropyConservativeSymbol$ inner product. We observe that all bases have equivalent approximation power for both the dimensional and non-dimensional configurations; this result is expected as all methods are dimensionally consistent. In addition, we observe that the trial subspaces obtained from the entropy inner product have a slightly different approximation power than the bases obtained from the non-dimensional $\NonDimensionalLTwoSymbol$ inner product, but the same trends are present. We observe the same convergence in the relative error for both methods and, in general, no method appears significantly ``better" than the other.

\begin{figure}
  \centering
  \begin{subfigure}{0.42\linewidth}
    \centering
    \begin{subfigure}[t]{0.9\textwidth}
    \hspace{0.24cm}
    \includegraphics[clip,trim={1cm 8.006cm 6.051cm 2cm},width=0.98\linewidth]{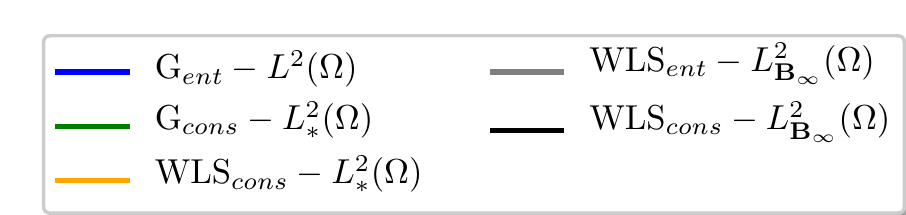}
    \end{subfigure}\\
    \begin{subfigure}[t]{1\textwidth}
    \includegraphics[clip,width=1.0\linewidth]{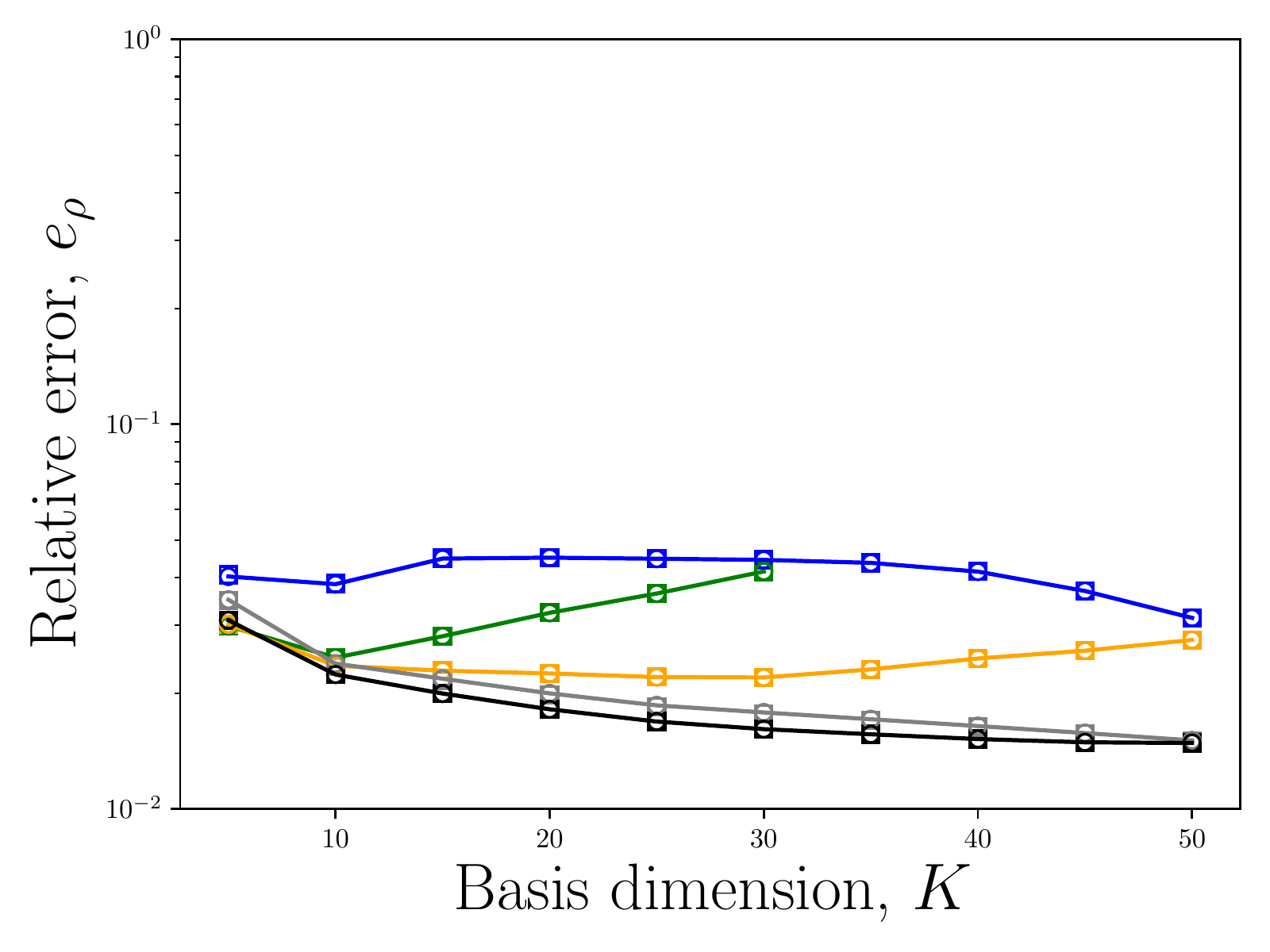}
    \end{subfigure}\\
    \begin{subfigure}[t]{1\textwidth}
    \includegraphics[clip,width=1.0\linewidth]{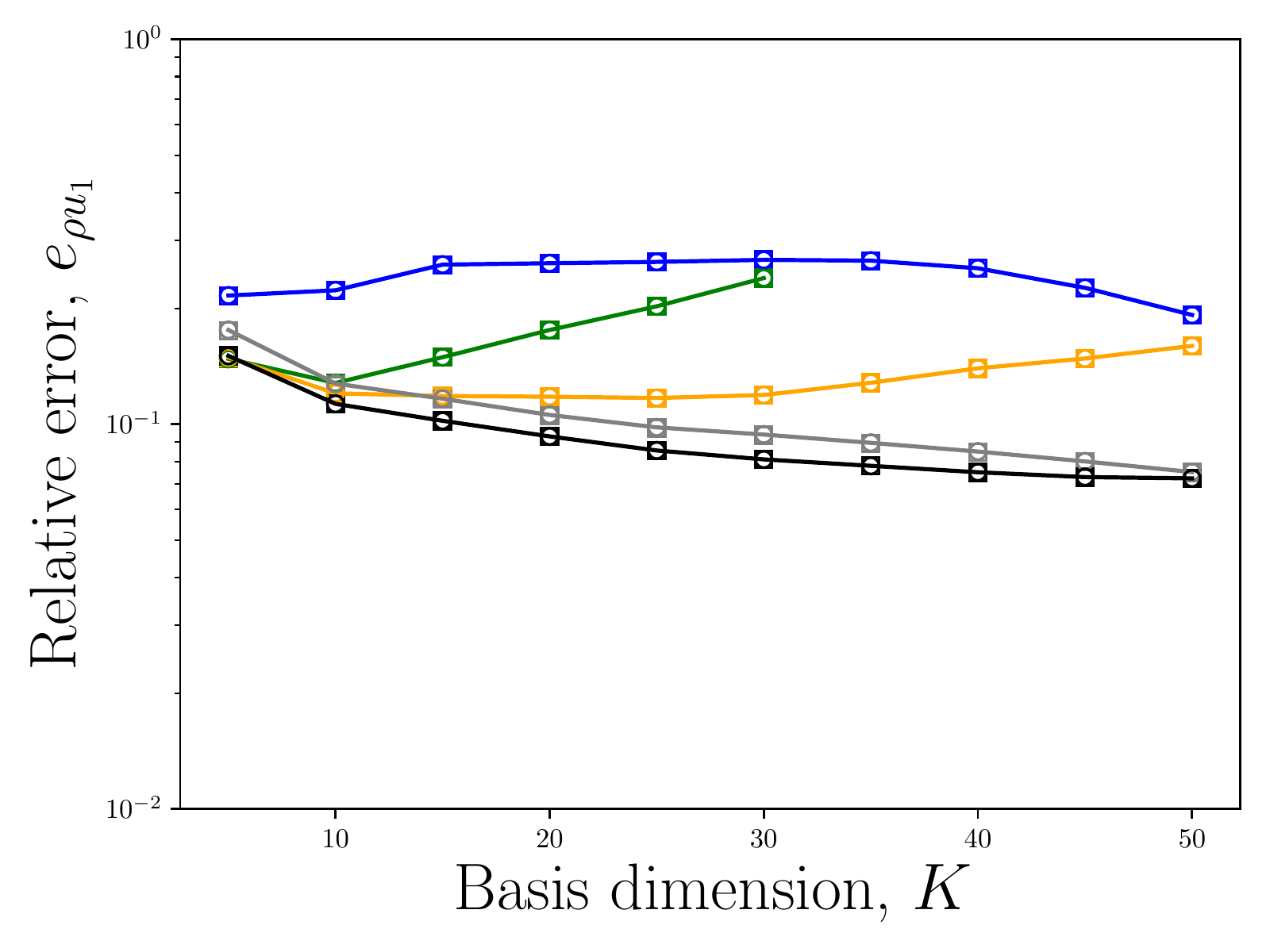}
    \end{subfigure}\\
    \begin{subfigure}[t]{1\textwidth}
    \includegraphics[clip,width=1.0\linewidth]{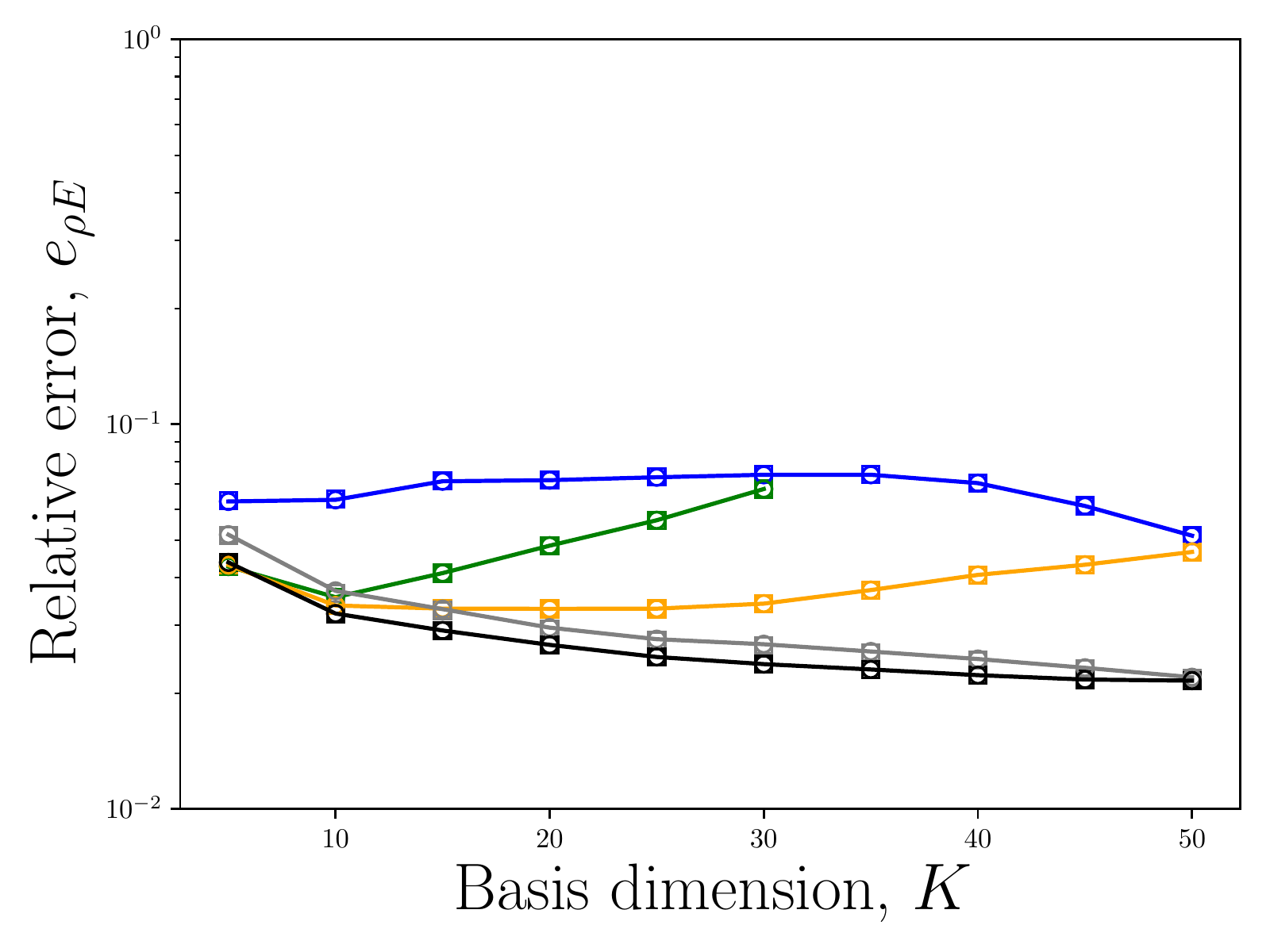}
    \end{subfigure}
  \end{subfigure}
  \begin{subfigure}{0.4\linewidth}
    \centering
    \begin{subfigure}[t]{0.88\textwidth}
    \includegraphics[clip,trim={0cm 1.4cm 0 0cm},width=1.0\linewidth]{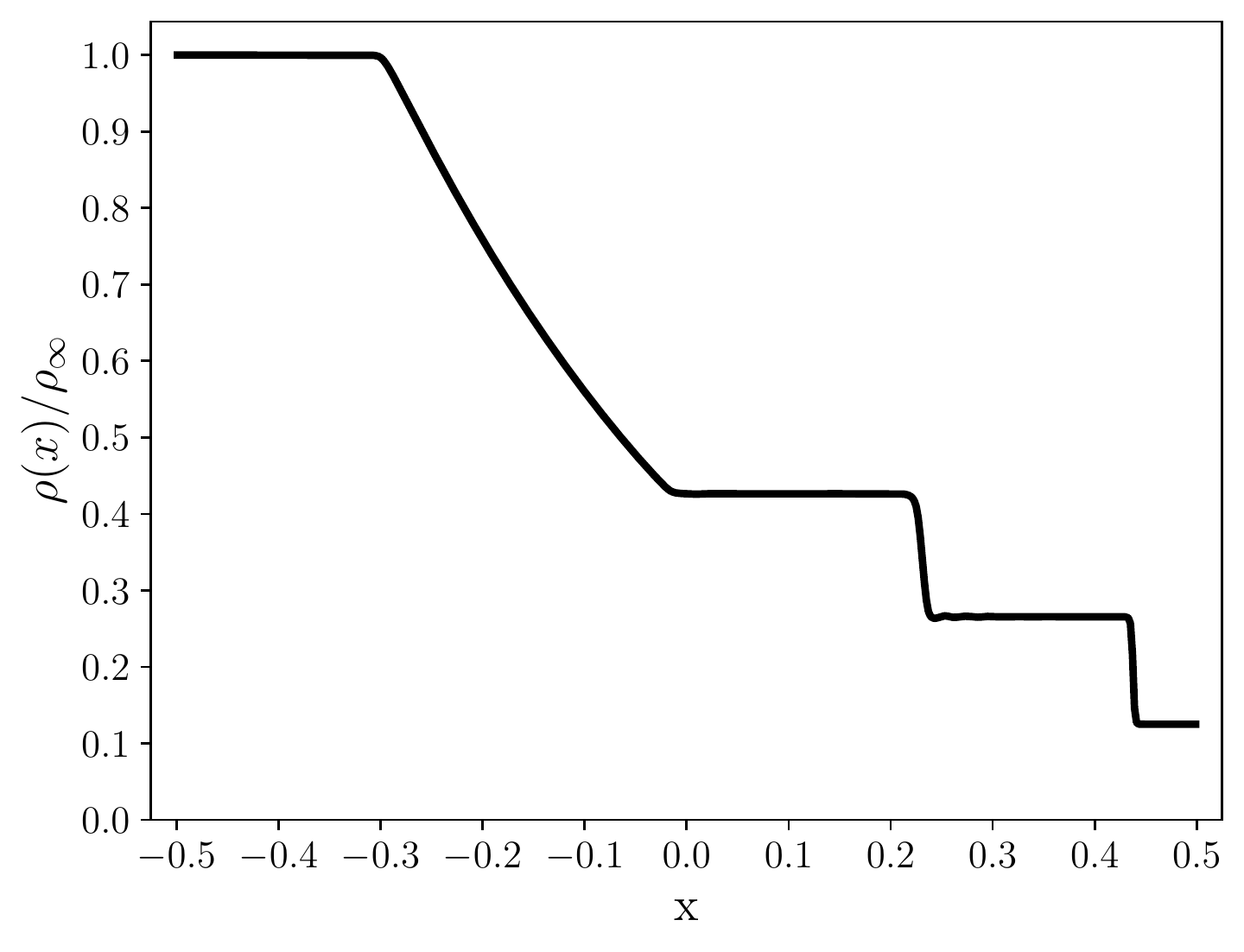}
    \caption{FOM}
    \end{subfigure}\\
    \begin{subfigure}[t]{0.88\textwidth}
    \includegraphics[clip,trim={0cm 1.4cm 0 0cm},width=1.0\linewidth]{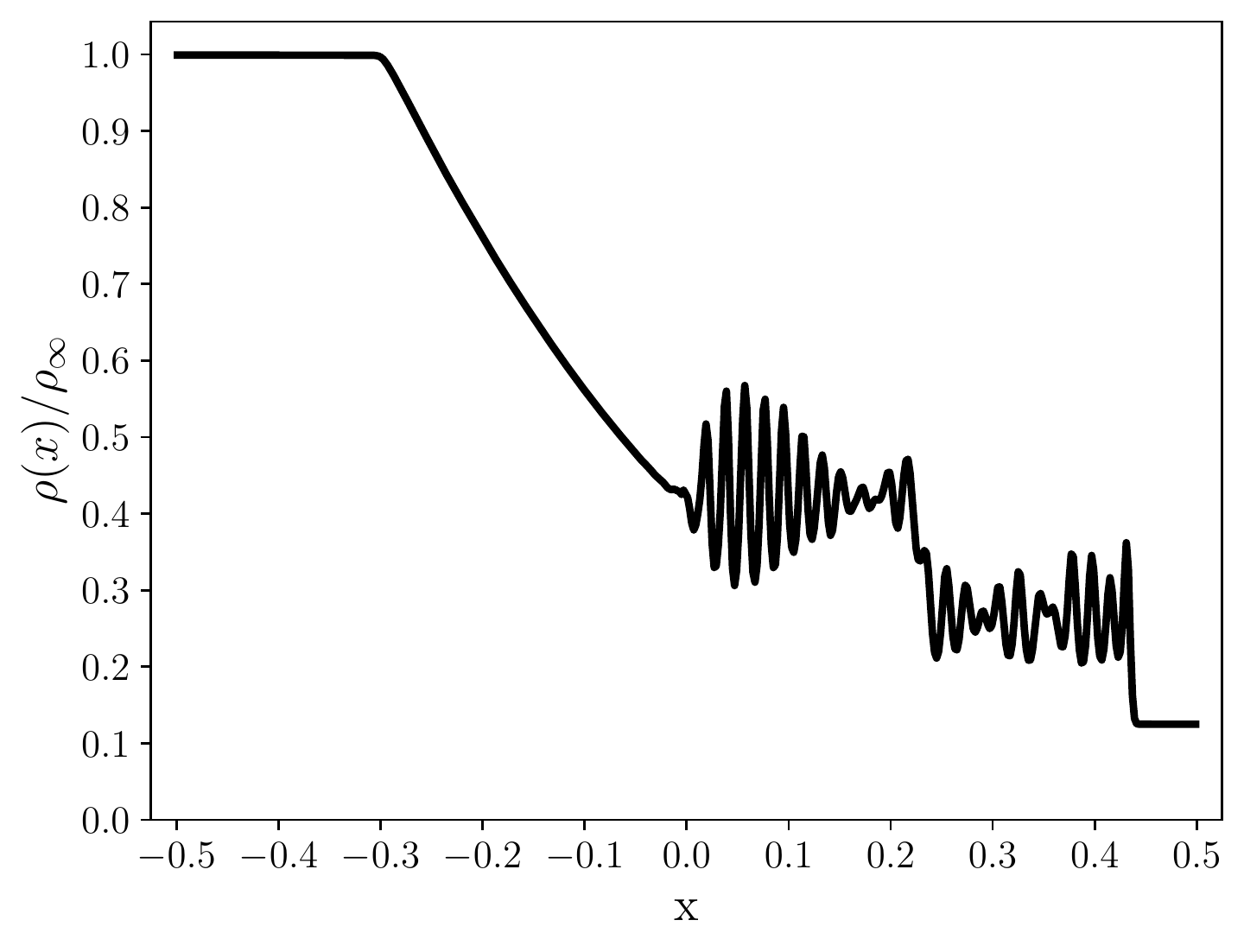}
    \caption{\GalerkinLTwoEntropyRomName}
    \end{subfigure}\\
    \begin{subfigure}[t]{0.88\textwidth}
    \includegraphics[clip,trim={0cm 1.4cm 0 0cm},width=1.0\linewidth]{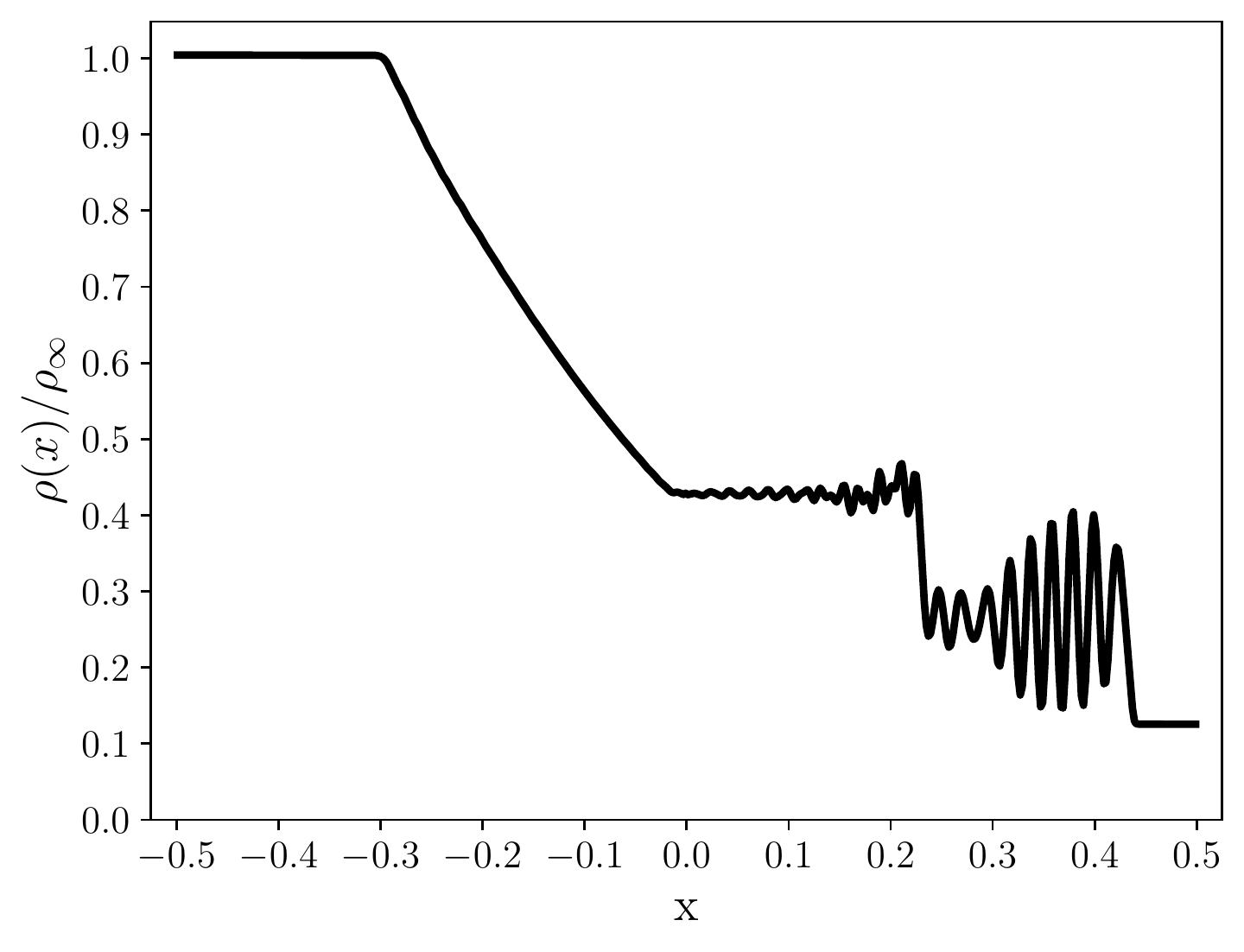}
    \caption{\WLSNonDimensionalLTwoConservedRomName}
    \end{subfigure}\\
    \begin{subfigure}[t]{0.88\textwidth}
    \includegraphics[clip,trim={0cm 1.4cm 0 0cm},width=1.0\linewidth]{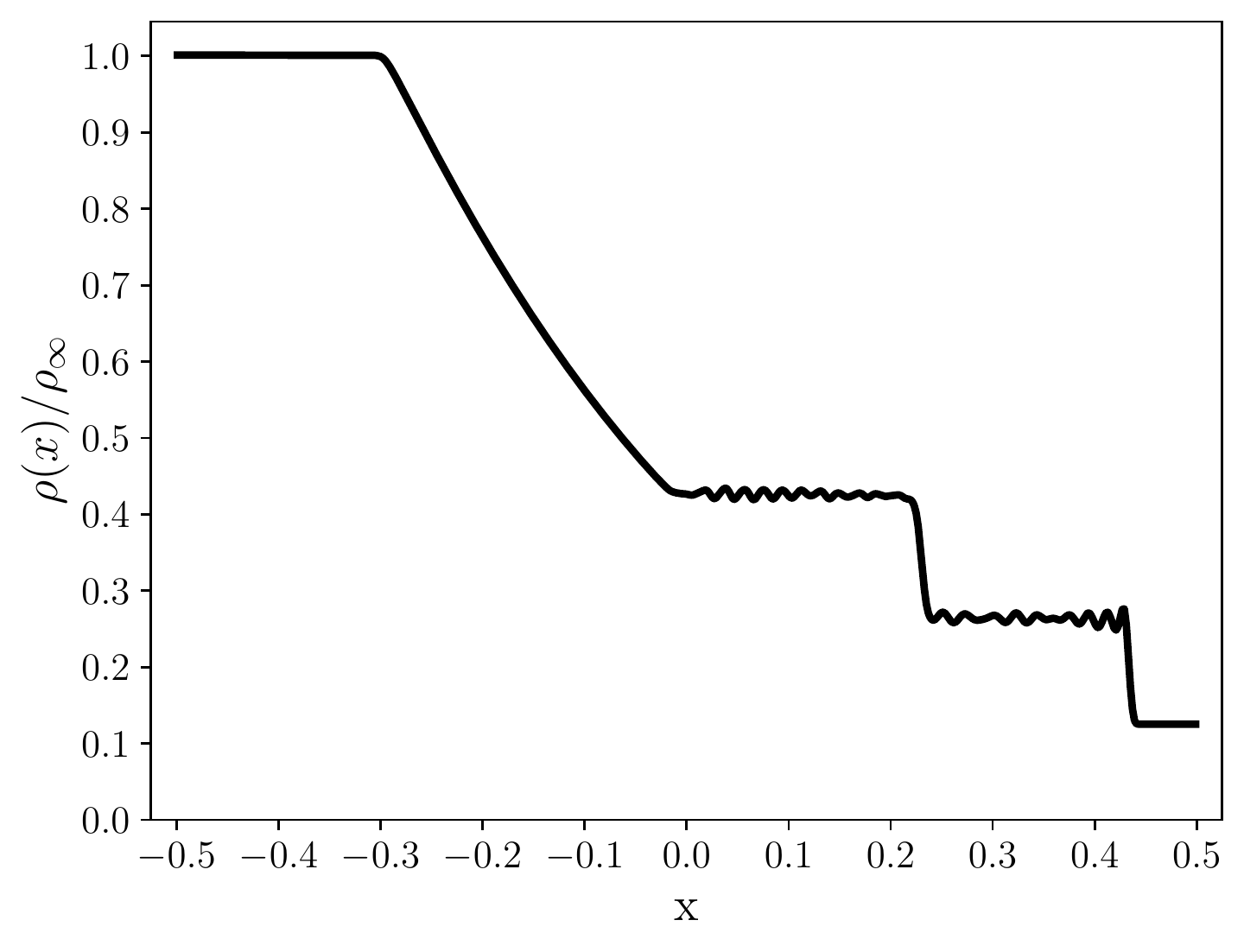}
    \caption{\WLSEntropyLTwoEntropyRomName}
    \end{subfigure}\\
    \begin{subfigure}[t]{0.88\textwidth}
    \includegraphics[clip,width=1.0\linewidth]{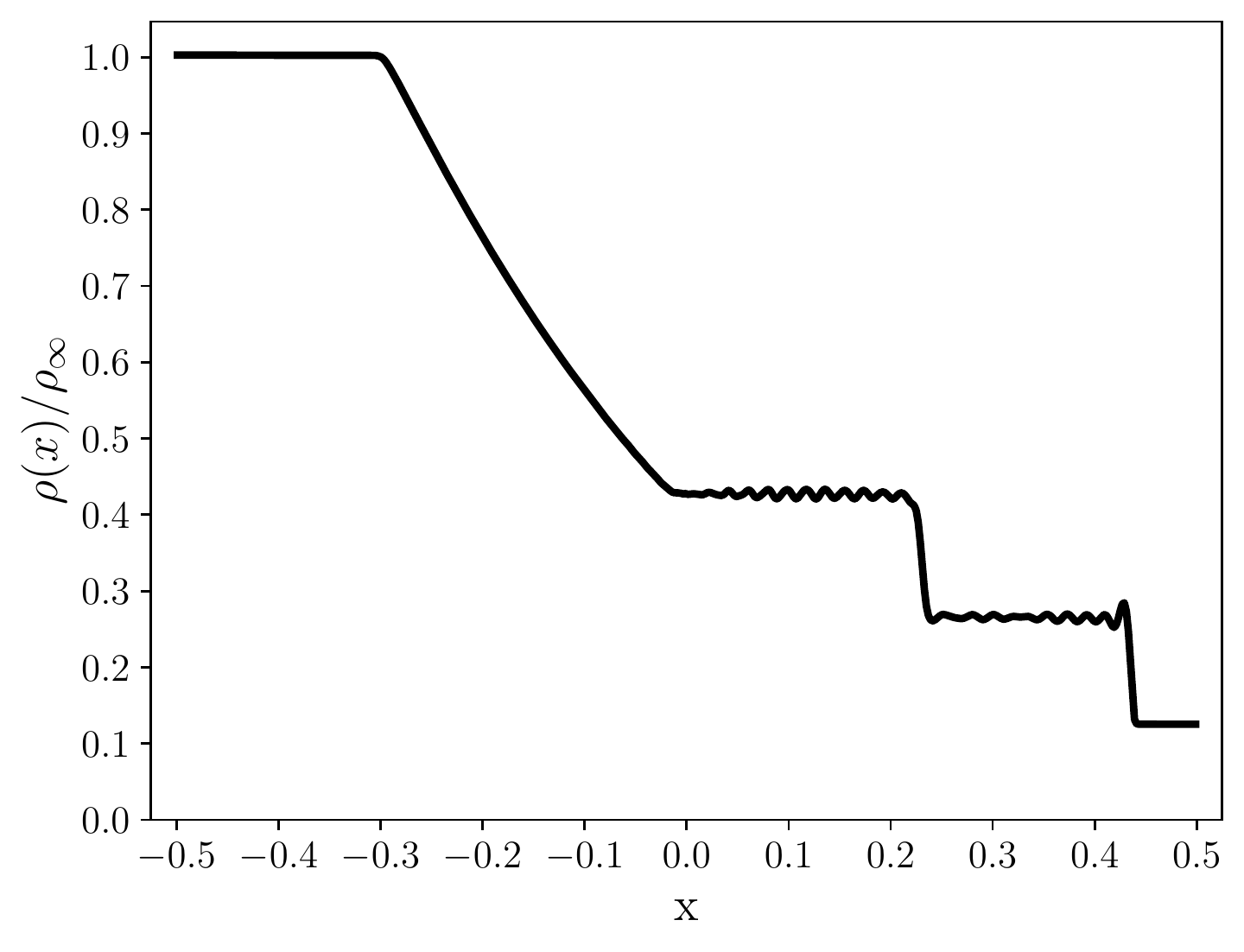}
    \caption{\WLSEntropyLTwoConservedRomName}
    \end{subfigure}
  \end{subfigure}
  \caption{Sod 1D problem. Convergence of ROM errors as a function of basis dimension (left) and physical space solutions for the density field at a non-dimensional time of $t=0.25$ (right) for a ROM dimension of $\romDim = 50$.}
  \label{fig:sod1dfig}
\end{figure}

\begin{figure}
  \centering
  \begin{subfigure}{0.42\linewidth}
    \centering
    \begin{subfigure}[t]{0.9\textwidth}
      \hspace{0.39cm}
    \includegraphics[clip,trim={1cm 8.006cm 6.051cm 2cm},width=0.945\linewidth]{figures/kelvin_helmholtz/legend.pdf}
    \end{subfigure}\\
    \begin{subfigure}[t]{1\textwidth}
    \includegraphics[clip,width=1.0\linewidth]{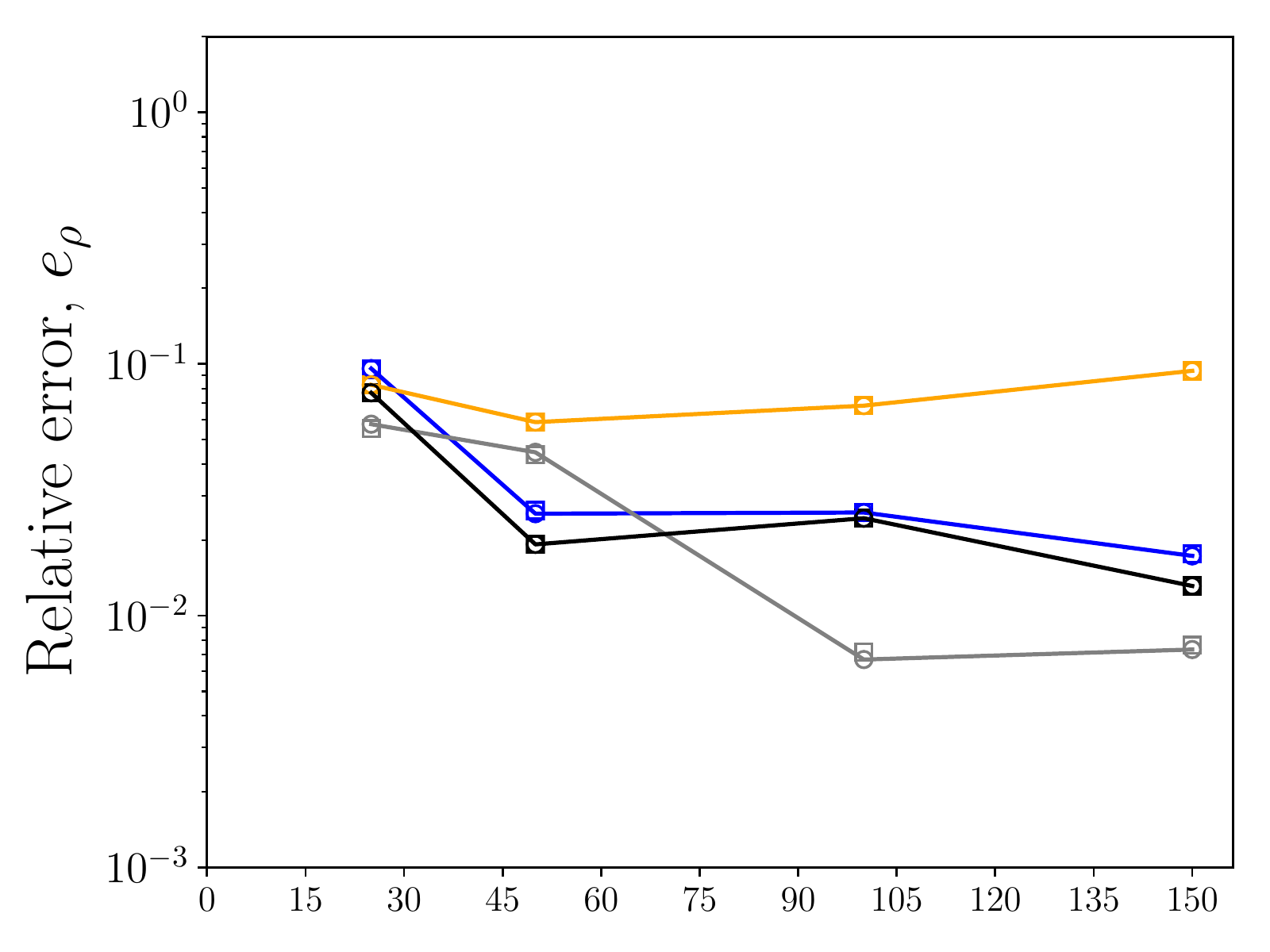}
    \end{subfigure}\\
    \begin{subfigure}[t]{1\textwidth}
    \includegraphics[clip,width=1.0\linewidth]{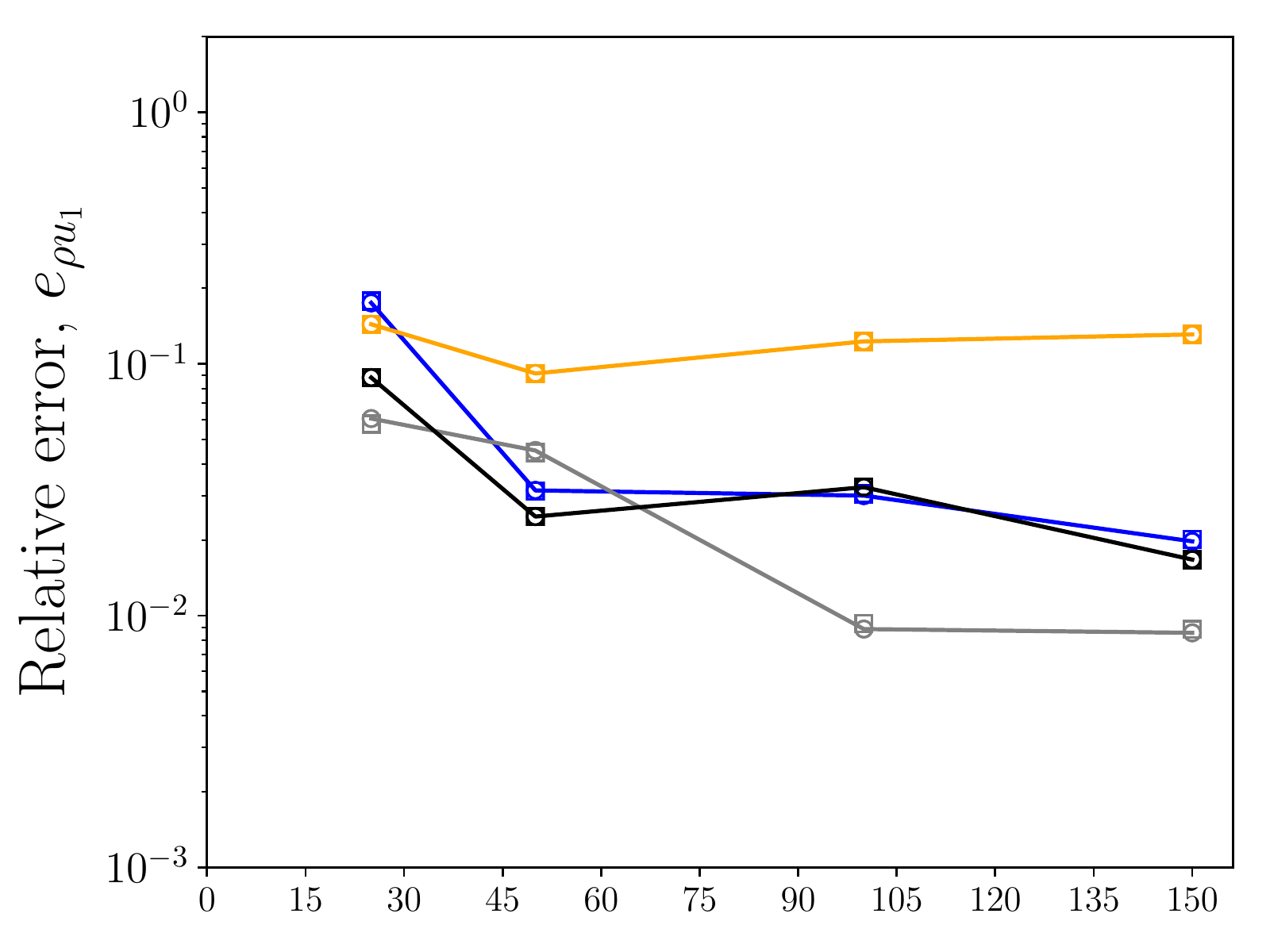}
    \end{subfigure}\\
    \begin{subfigure}[t]{1\textwidth}
    \includegraphics[clip,width=1.0\linewidth]{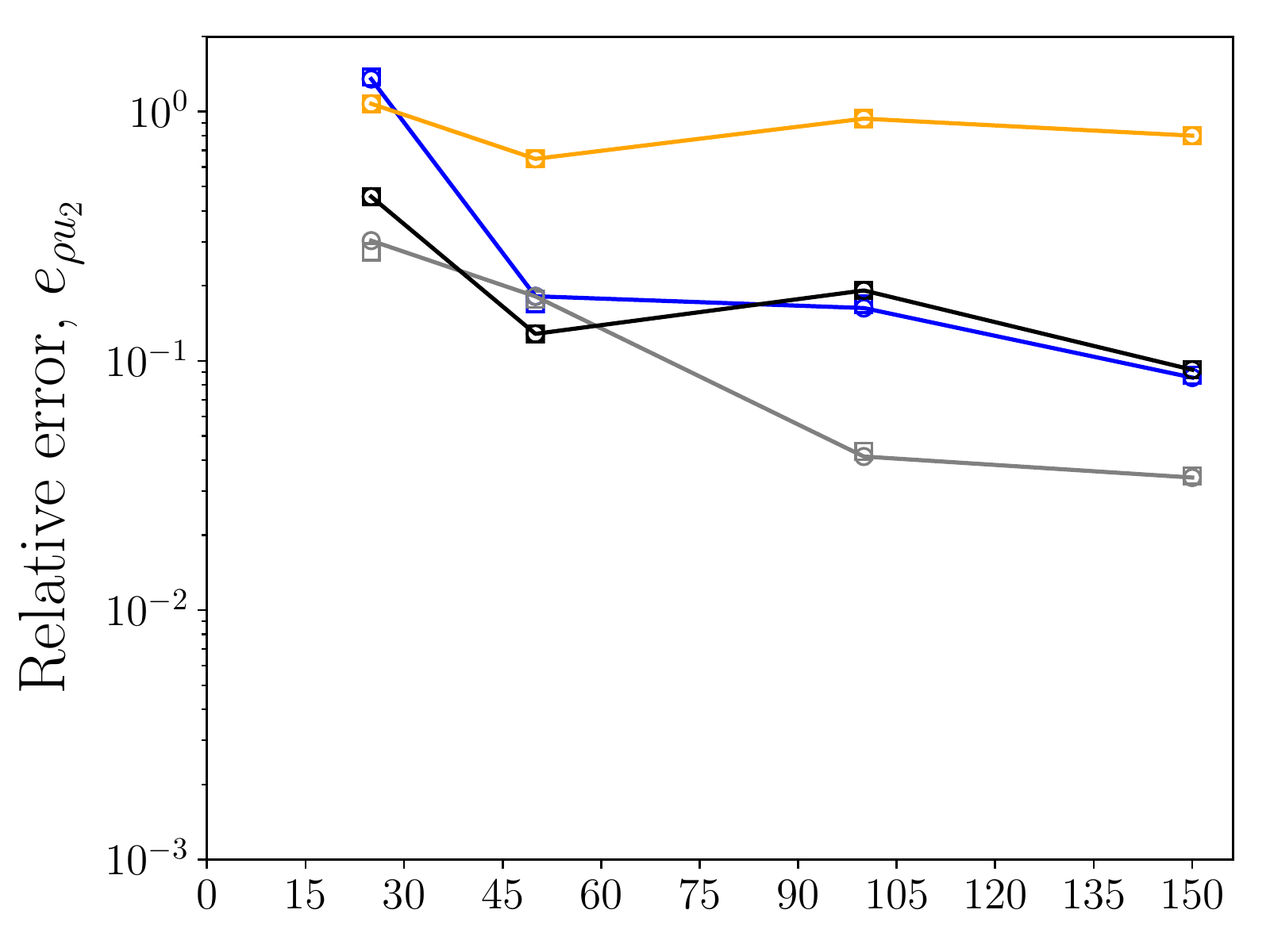}
    \end{subfigure}\\
    \begin{subfigure}[t]{1\textwidth}
    \includegraphics[clip,width=1.0\linewidth]{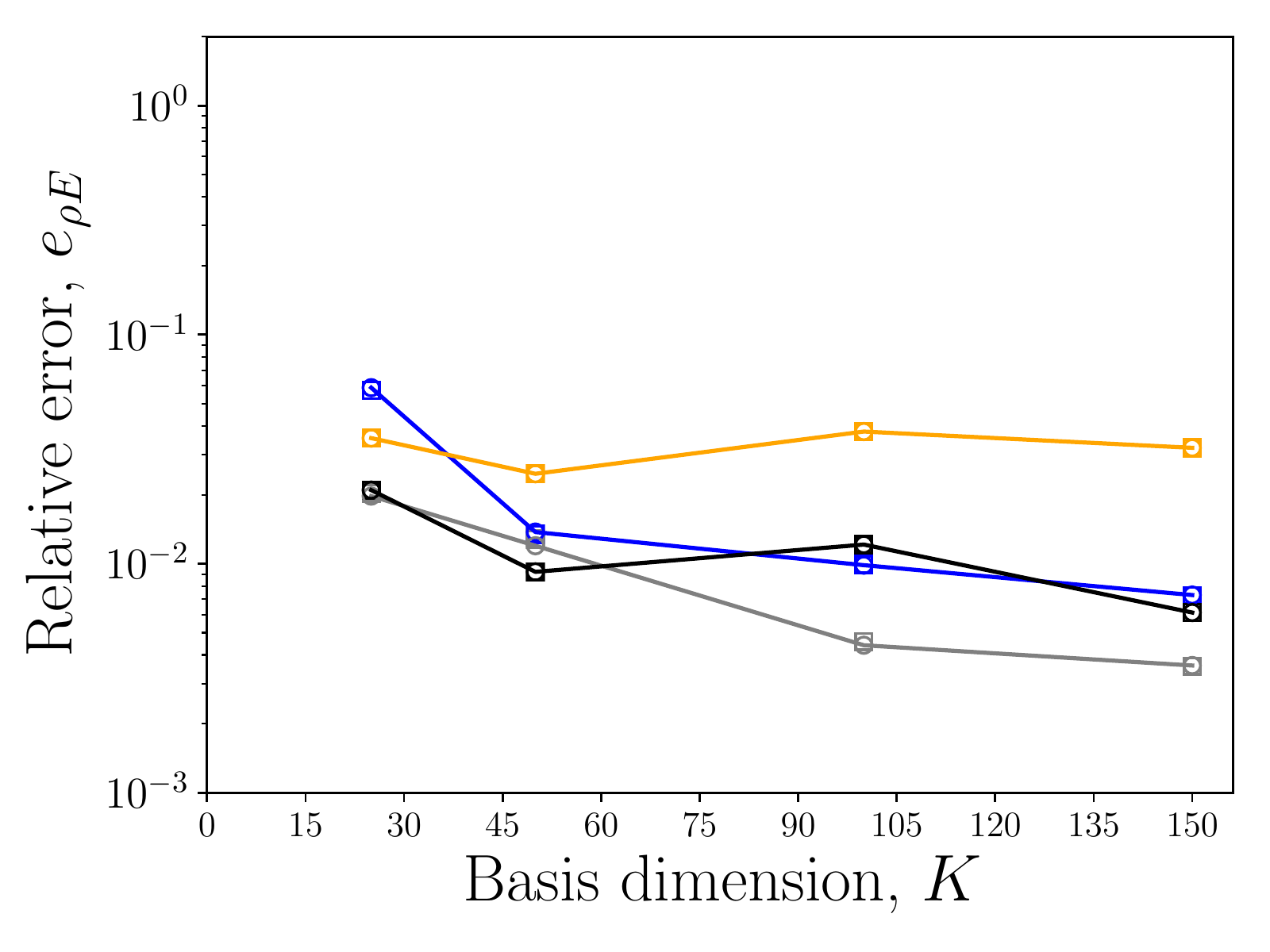}
    \end{subfigure}
  \end{subfigure}
  \begin{subfigure}{0.4\linewidth}
    \centering
    \begin{subfigure}[t]{0.9\textwidth}
    \includegraphics[clip,trim={0cm 1.32cm 0 0cm},width=1.0\linewidth]{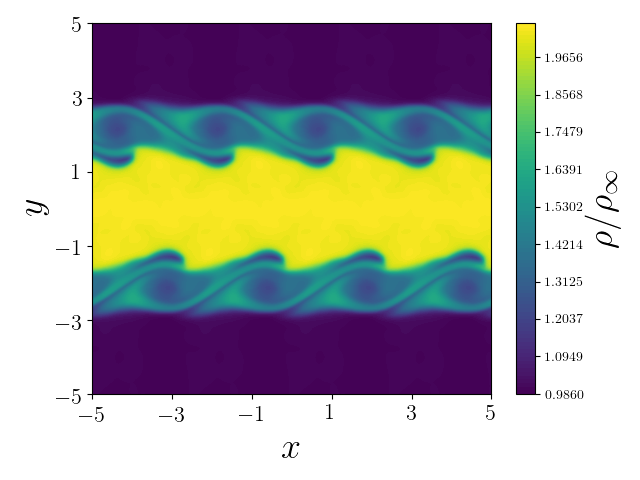}
    \caption{FOM}
    \end{subfigure}\\
    \begin{subfigure}[t]{0.9\textwidth}
    \includegraphics[clip,trim={0cm 1.32cm 0 0cm},width=1.0\linewidth]{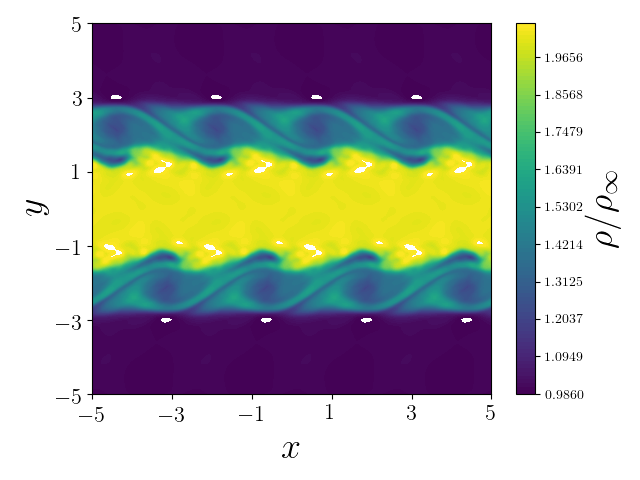}
    \caption{\GalerkinLTwoEntropyRomName}
    \end{subfigure}\\
    \begin{subfigure}[t]{0.9\textwidth}
    \includegraphics[clip,trim={0cm 1.32cm 0 0cm},width=1.0\linewidth]{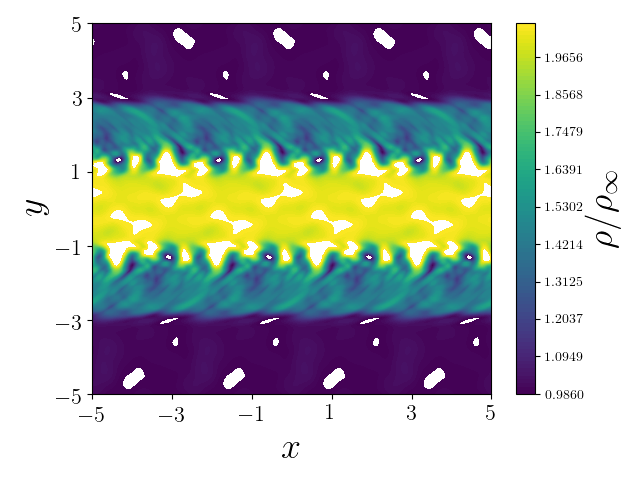}
    \caption{\WLSNonDimensionalLTwoConservedRomName}
    \end{subfigure}\\
    \begin{subfigure}[t]{0.9\textwidth}
    \includegraphics[clip,trim={0cm 1.32cm 0 0cm},width=1.0\linewidth]{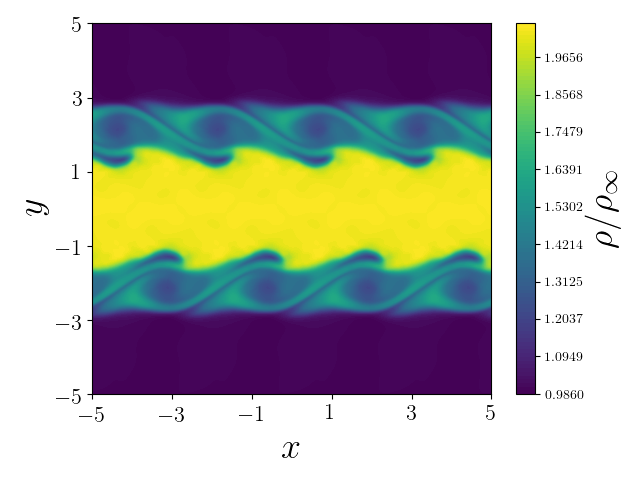}
    \caption{\WLSEntropyLTwoEntropyRomName}
    \end{subfigure}\\
    \begin{subfigure}[t]{0.9\textwidth}
    \includegraphics[clip,trim={0cm 0cm 0 0cm},width=1.0\linewidth]{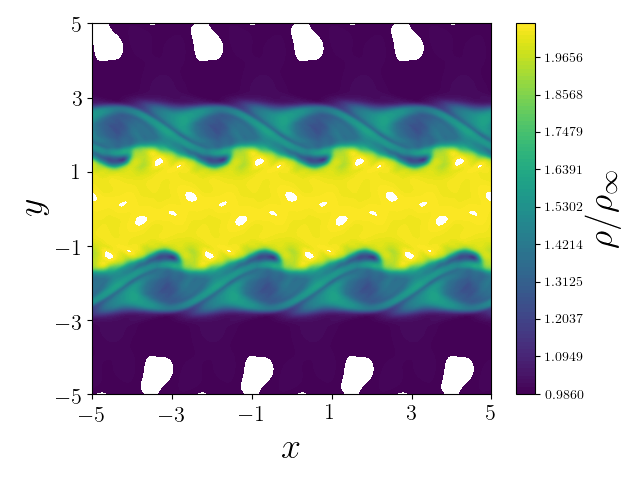}
    \caption{\WLSEntropyLTwoConservedRomName}
    \end{subfigure}
  \end{subfigure}
  \caption{Kelvin Helmholtz problem. Convergence of ROM errors as a function of basis dimension (left) and physical space solutions for the density field at a non-dimensional time of  $t=19.5$ (right) with a ROM dimension of $\romDim = 150$.}
  \label{fig:kh2dfig}
\end{figure}

\begin{figure}
  \centering
  \begin{subfigure}{0.42\linewidth}
    \centering
    \begin{subfigure}[t]{0.9\textwidth}
    \hspace{0.39cm}
    \includegraphics[clip,trim={1cm 8.006cm 6.051cm 2cm},width=0.925\linewidth]{figures/kelvin_helmholtz/legend.pdf}
    \end{subfigure}\\
    \begin{subfigure}[t]{1\textwidth}
    \includegraphics[clip,width=1.0\linewidth]{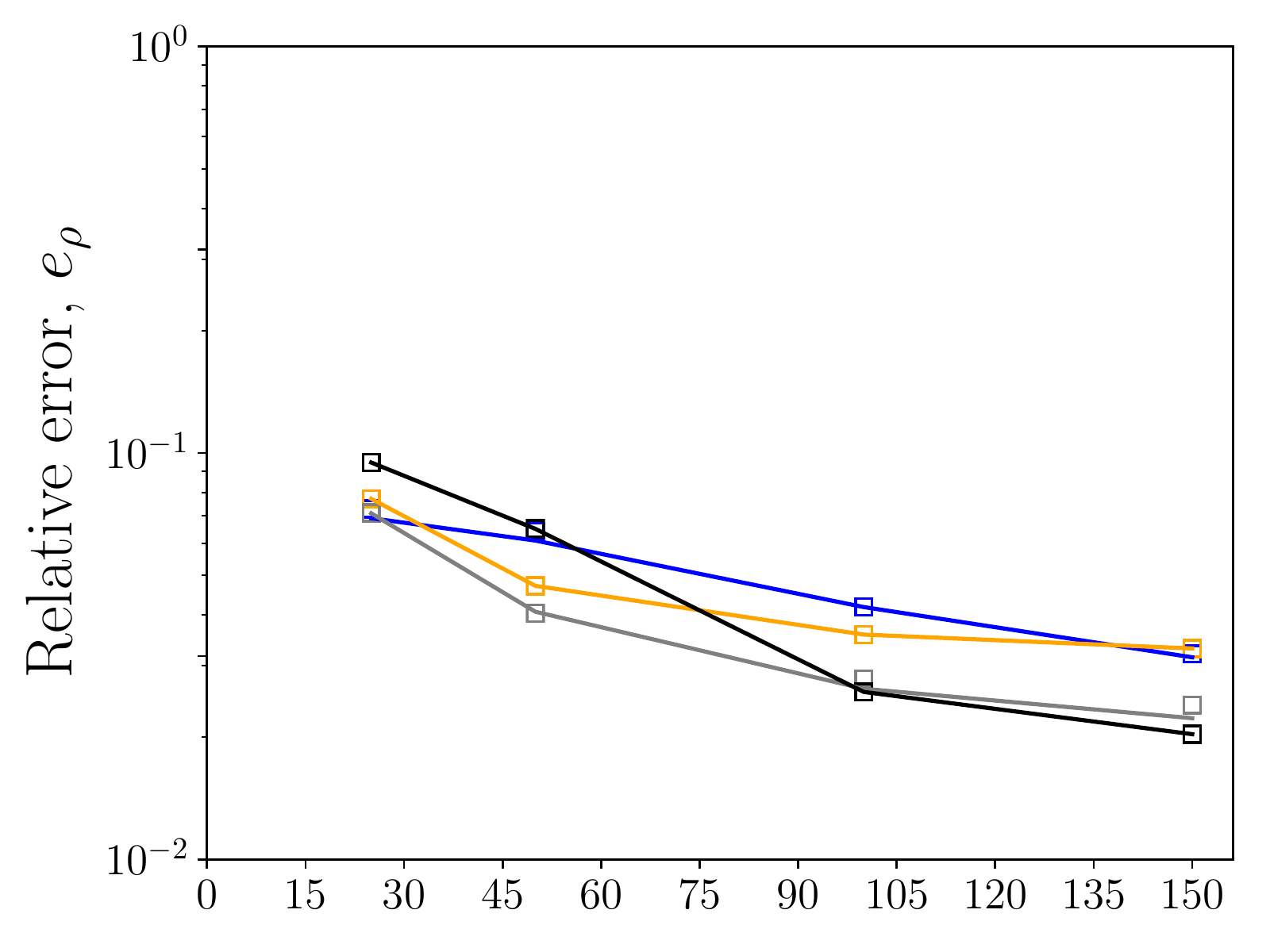}
    \end{subfigure}\\
    \begin{subfigure}[t]{1\textwidth}
    \includegraphics[clip,width=1.0\linewidth]{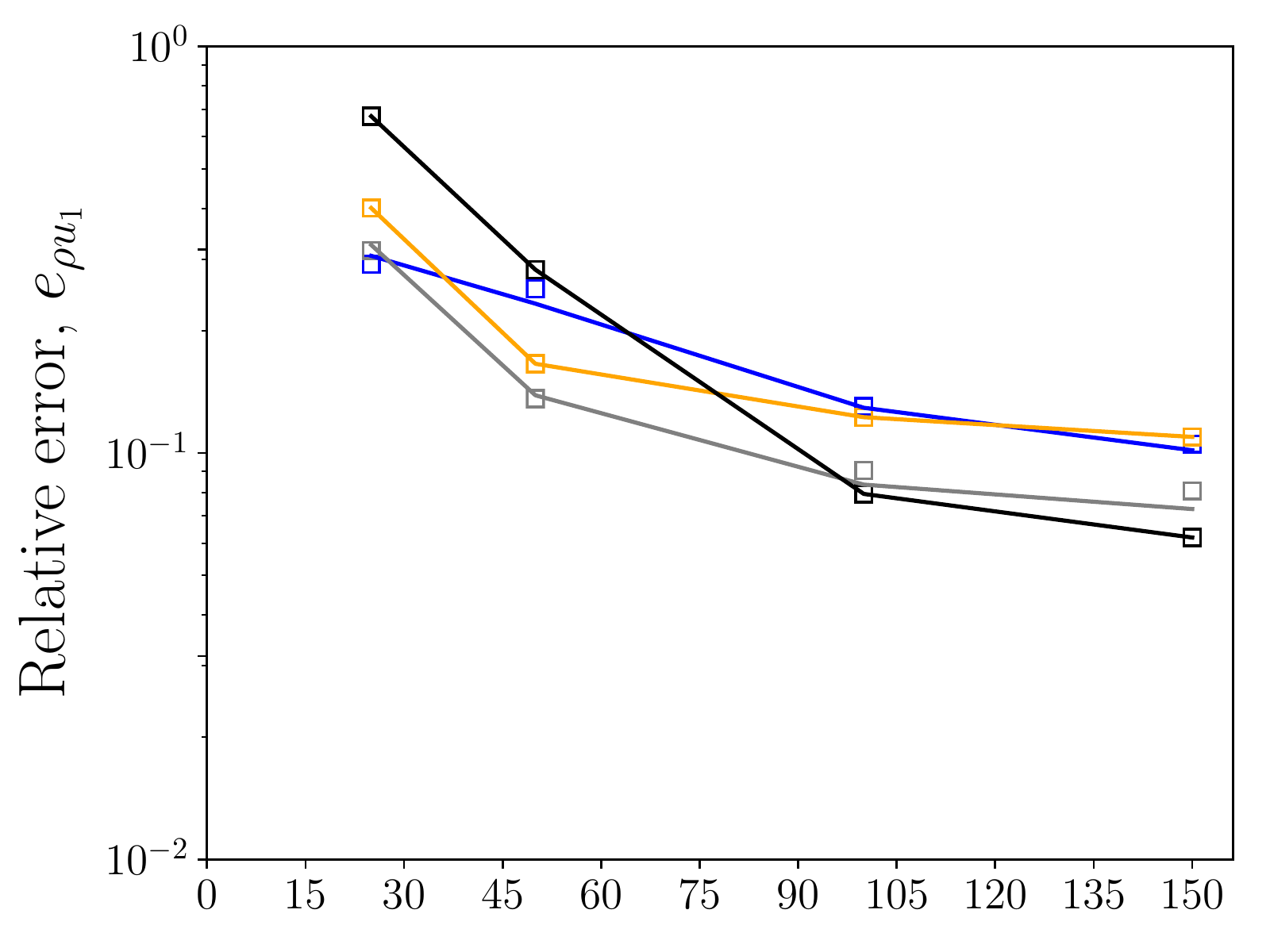}
    \end{subfigure}\\
    \begin{subfigure}[t]{1\textwidth}
    \includegraphics[clip,width=1.0\linewidth]{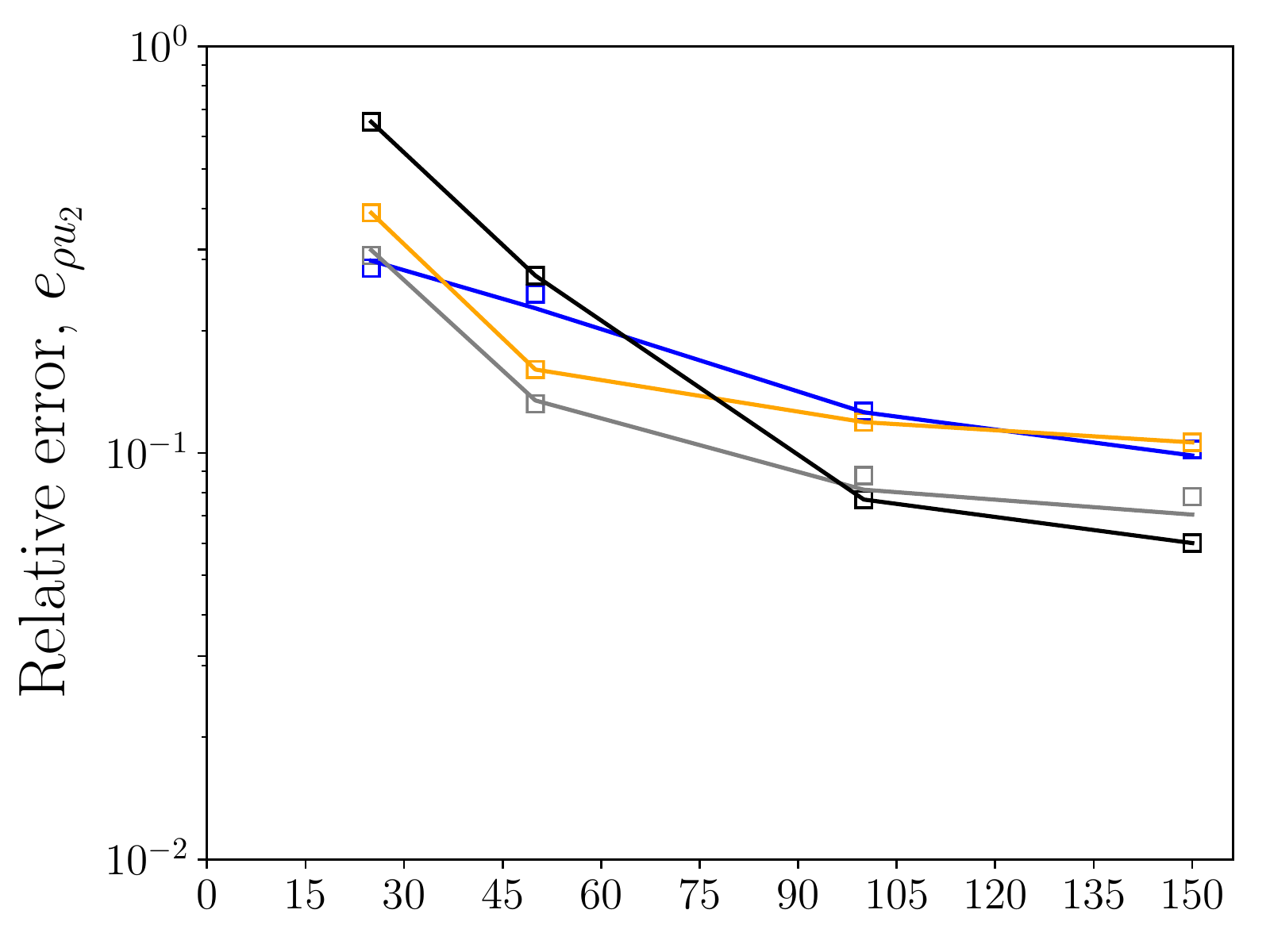}
    \end{subfigure}\\
    \begin{subfigure}[t]{1\textwidth}
    \includegraphics[clip,width=1.0\linewidth]{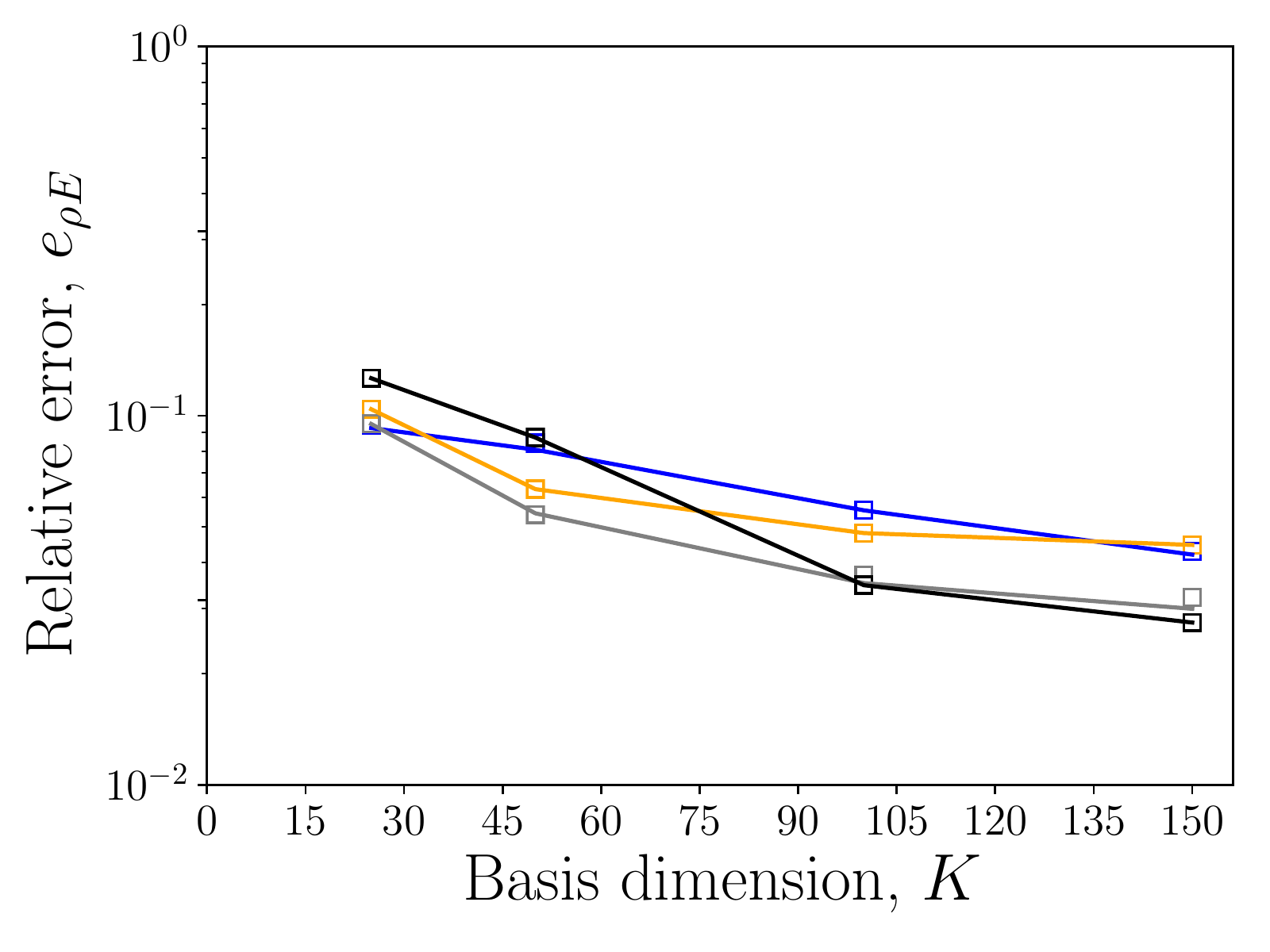}
    \end{subfigure}
  \end{subfigure}
  \begin{subfigure}{0.4\linewidth}
    \centering
    \begin{subfigure}[t]{0.9\textwidth}
    \includegraphics[clip,trim={0cm 1.32cm 0 0cm},width=1.0\linewidth]{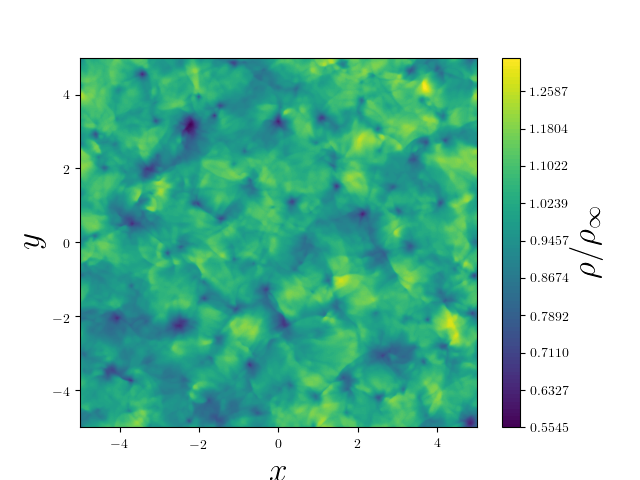}
    \caption{FOM}
    \end{subfigure}\\
    \begin{subfigure}[t]{0.9\textwidth}
    \includegraphics[clip,trim={0cm 1.32cm 0 0cm},width=1.0\linewidth]{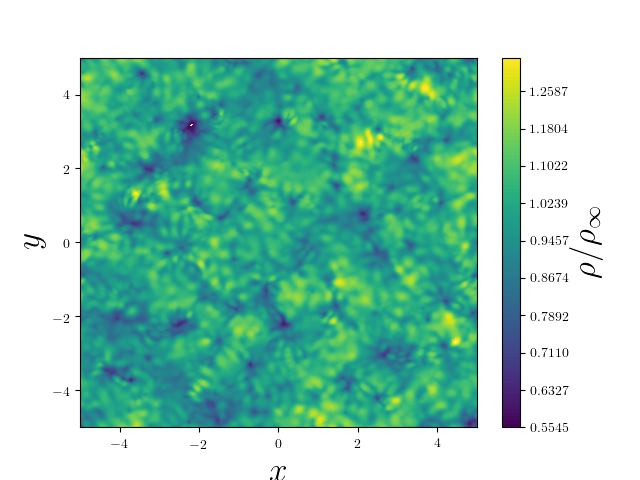}
    \caption{\GalerkinLTwoEntropyRomName}
    \end{subfigure}\\
    \begin{subfigure}[t]{0.9\textwidth}
    \includegraphics[clip,trim={0cm 1.32cm 0 0cm},width=1.0\linewidth]{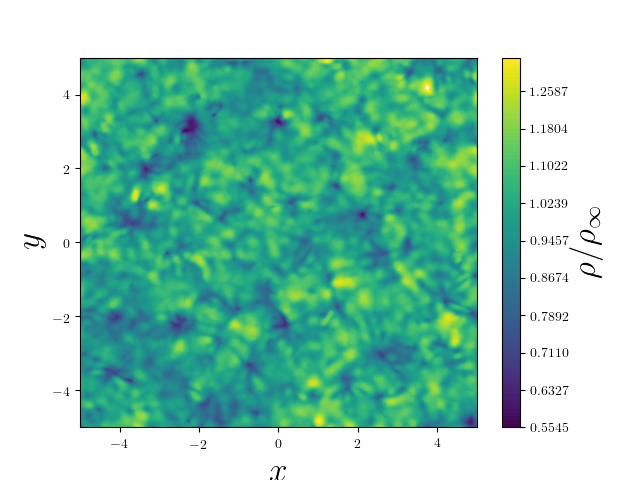}
    \caption{\WLSNonDimensionalLTwoConservedRomName}
    \end{subfigure}\\
    \begin{subfigure}[t]{0.9\textwidth}
    \includegraphics[clip,trim={0cm 1.32cm 0 0cm},width=1.0\linewidth]{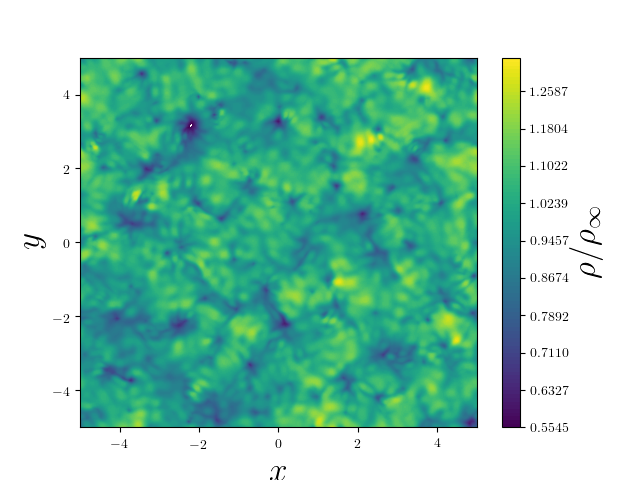}
    \caption{\WLSEntropyLTwoEntropyRomName}
    \end{subfigure}\\
    \begin{subfigure}[t]{0.9\textwidth}
    \includegraphics[clip,trim={0cm 0cm 0 0cm},width=1.0\linewidth]{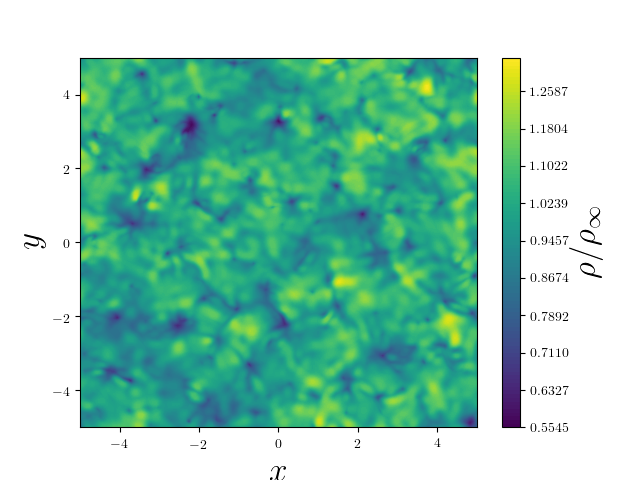}
    \caption{\WLSEntropyLTwoConservedRomName}
    \end{subfigure}
  \end{subfigure}
  \caption{Homogeneous turbulence problem. Convergence of ROM errors as a function of basis dimension (left) and physical space solutions for the density field at a non-dimensional time of  $t=14.7$ (right) with a ROM dimension of $\romDim = 100$.}
  \label{fig:hit2dfig}
\end{figure}

We next consider the various ROM formulations described above applied to the dimensional and non-dimensional FOMs. Figures~\ref{fig:sod1dfig}-\ref{fig:hit2dfig} present results for the various ROMs applied to the 1D Sod problem, Kelvin Helmholtz problem, and 2D homogeneous turbulence problem, respectively. The left columns of Figures~\ref{fig:sod1dfig}-\ref{fig:hit2dfig} present the convergence of the relative errors for mass, momentum, and energy as a function of the basis dimension, while the right columns present physical-space solutions. Results are presented for a basis dimension of $\romDim = 50$, $150$, and $150$ for the Sod 1D, Kelvin Helmholtz, and homogeneous turbulence problems, respectively. Across all solutions, we make the following observations:
\begin{itemize}
\item (Figures~\ref{fig:sod1dfig}-\ref{fig:hit2dfig}.) WLS ROMs based on the entropy inner product consistently out-perform the standard WLS ROMs based on the $\NonDimensionalLTwoSymbol$ inner product. We observe the entropy-based  WLS ROMs to (1) yield solutions that are qualitatively more accurate than the \WLSNonDimensionalLTwoConservedRomName\ ROM and (2) display better convergence properties. This is particularly evident in the field solutions for the Sod and Kelvin Helmholtz problems: the \WLSNonDimensionalLTwoConservedRomName\ ROM results in poor solutions in both cases, and simply changing the ROM formulation to measure the residual in the $\EntropyConservativeSymbol$ inner product results in a substantial increase in accuracy. We emphasize that the \WLSEntropyLTwoConservedRomName\ ROM employs the same basis as the \WLSNonDimensionalLTwoConservedRomName\ ROM, and thus the improvements are solely due to the ROM formulation.

\item (Figures~\ref{fig:sod1dfig}-\ref{fig:hit2dfig}.)  The \GalerkinLTwoEntropyRomName\ ROM based on the entropy inner product is stable for all ROM dimensions, while the \GalerkinNonDimensionalLTwoConservedRomName\ ROM based on the non-dimensional $L^2_*(\Omega)$ inner product is rarely stable. In particular, we observe that the \GalerkinNonDimensionalLTwoConservedRomName\ ROM often becomes less stable as the basis dimension grows, while convergence is, in general, observed for the \GalerkinLTwoEntropyRomName\ ROM.
\item (Figures~\ref{fig:sod1dfig}-\ref{fig:hit2dfig}.)  All methods yield the same results on the dimensional and non-dimensional configurations; this is again expected as all methods are dimensionally consistent.
\end{itemize}

\section{Summary of numerical results}\label{sec:numerical_summary}
\begin{figure}
\begin{center}
\begin{subfigure}[t]{0.44\textwidth}
\includegraphics[clip,width=1.0\linewidth]{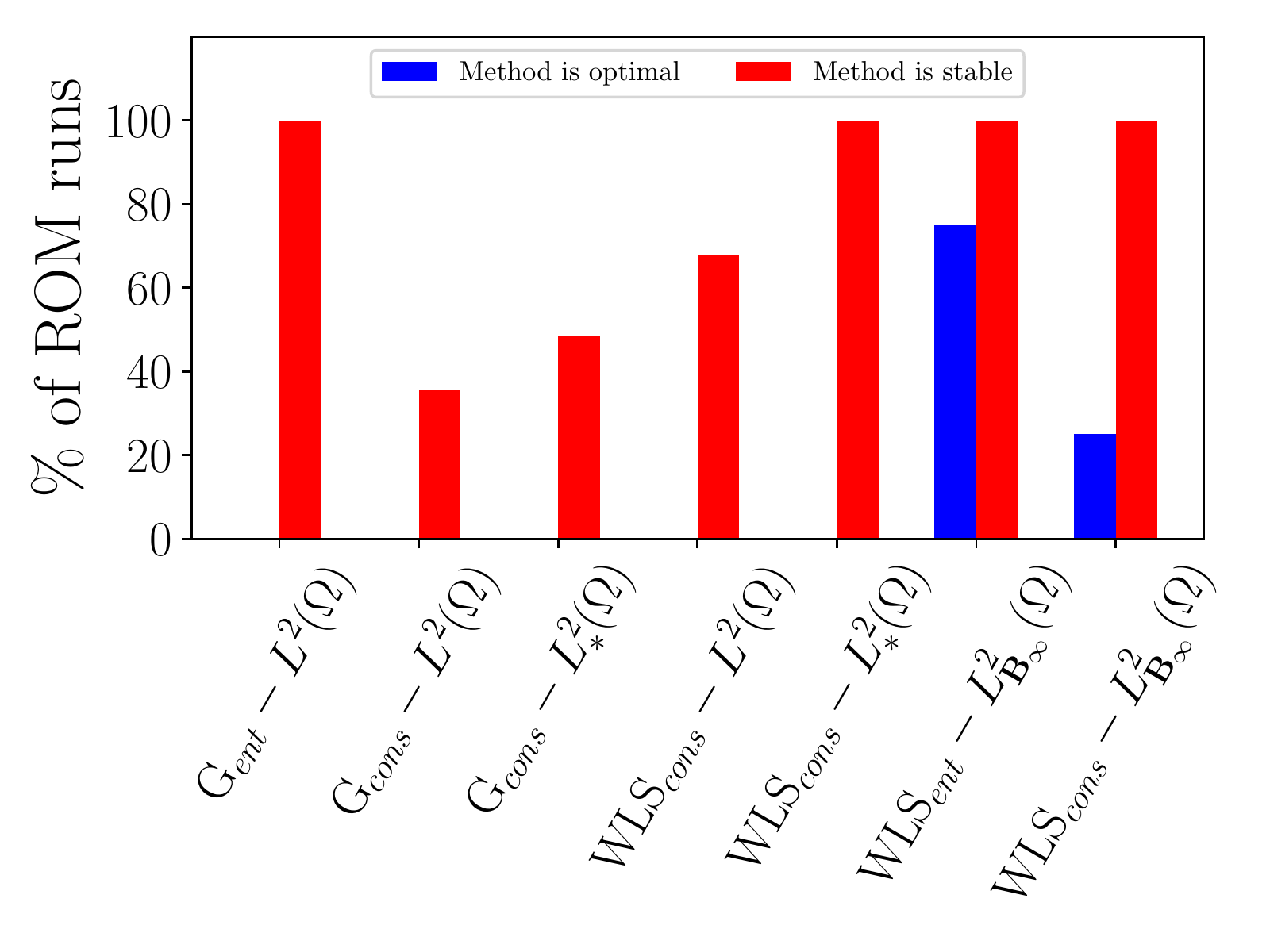}
\caption{Dimensional}
\end{subfigure}
\begin{subfigure}[t]{0.44\textwidth}
\includegraphics[clip,width=1.0\linewidth]{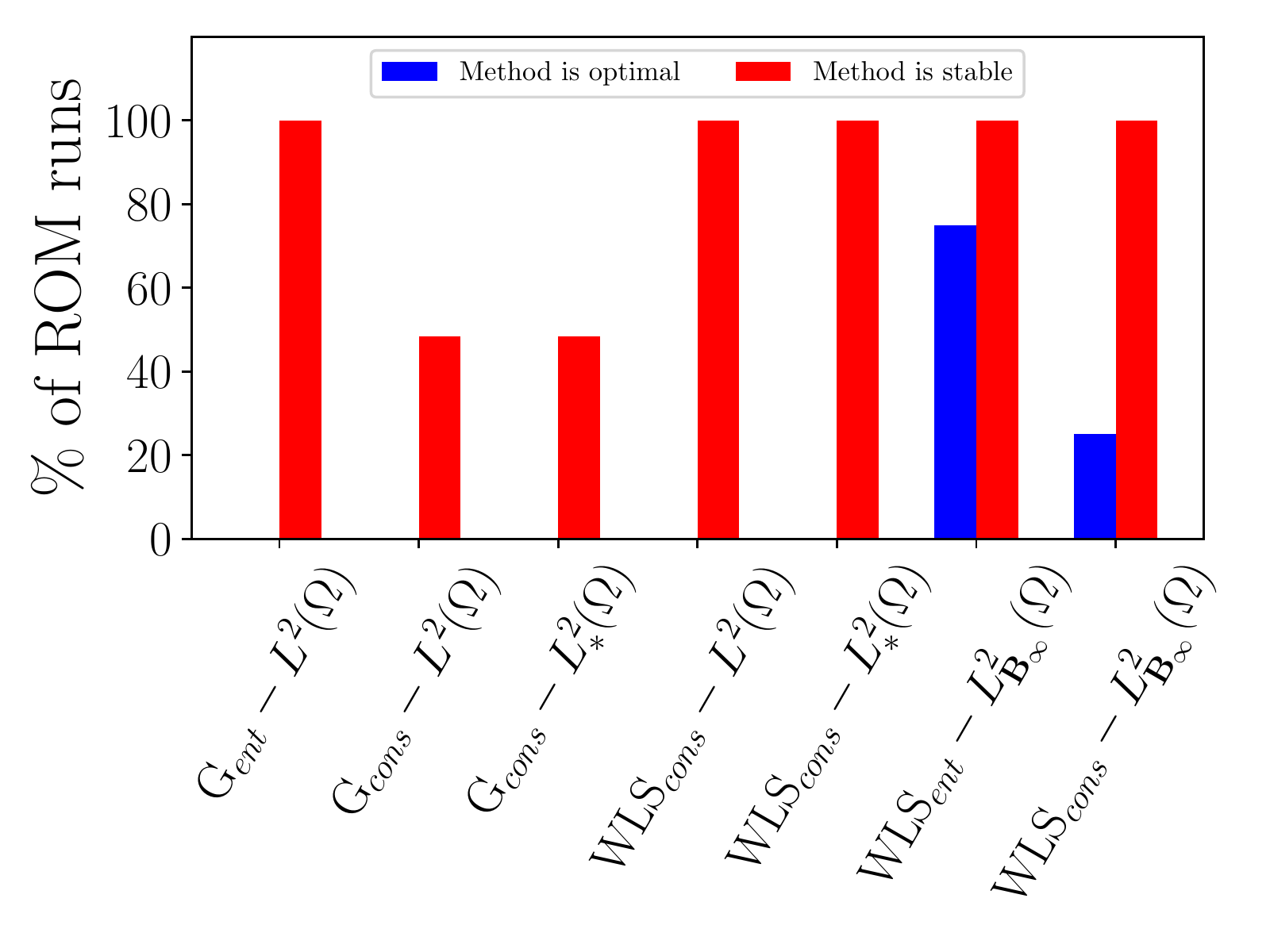}
\caption{Non-dimensional}
\end{subfigure}
\caption{Percentage of the time a ROM formulation gave a solution with the lowest mean-squared-error and percentage of time a given ROM formulation was stable. Results are compiled over all bases dimensions and state variables. The dimensional configuration is shown on the left while the non-dimensional configuration is shown on the right.}
\label{fig:romSummary}
\end{center}
\end{figure}

To summarize our numerical results, Figure~\ref{fig:romSummary} presents the percentage of time a given ROM formulation resulted in the most accurate solution for a state variable and the number of times that a given ROM formulation was stable. Results are presented for both the dimensional configuration (left) and non-dimensional configuration (right). The results are aggregated across all ROM dimensions and all state variables. We make the following remarks:
\begin{itemize}
\item ROMs based on the entropy inner products are always stable for both the dimensional and non-dimensional configurations.
\item The WLS-Entropy ROMs using either a conserved basis or entropy basis always result in the most accurate predictions.
\item The Galerkin ROMs are rarely stable. Employing a non-dimensional $\NonDimensionalLTwoSymbol$ inner product helps slightly, but performance is still poor.
\item The non-dimensional \WLSNonDimensionalLTwoConservedRomName\ ROM is always stable while the dimensional \WLSLTwoConservedRomName\ is not.
\item The WLS-Entropy ROM with an entropy basis results in the most accurate prediction for a state variable $80\%$ of the time.
\end{itemize}
Lastly, compiling results across our examples, we can rank the performance of our dimensionally-consistent ROMs in terms of solution accuracy from best-to-worst as:
\begin{enumerate}
\item \WLSEntropyLTwoEntropyRomName,
\item \WLSEntropyLTwoConservedRomName,
\item \GalerkinLTwoEntropyRomName,
\item \WLSNonDimensionalLTwoConservedRomName,
\item \GalerkinNonDimensionalLTwoConservedRomName,
\end{enumerate}
where we note it is a toss-up between (3) and (4).
\section{Conclusions}\label{sec:summary}
This work investigated the impact of the choice of inner product on the performance of POD Galerkin and POD least-squares reduced-order models within the context of the compressible Euler equations. We demonstrated that employing the classic $\LTwoSymbol$ inner product is inappropriate for POD, Galerkin, and least-squares ROMs discretized with conserved variables. This inner product resulted in ROMs that were dimensionally inconsistent and resulted in different solutions when applied to non-dimensional and dimensional versions of the same problem. Further, these formulations were often unstable. We showed that this issued can be alleviated by employing a non-dimensional $\NonDimensionalLTwoSymbol$ inner product, in which case a dimensionally-consistent POD process and Galerkin/least-squares ROM can be formulated. Our numerical experiments demonstrated that employing this non-dimensional inner product led to ROMs that were consistently more accurate and more stable. However, even when employing this non-dimensional inner product, the Galerkin ROM was still often unstable and least-squares ROMs, while stable, were not particularly accurate. We demonstrated that the use of a physics-based ``entropy" inner product addressed this issue and vastly improved the ROM performance. We demonstrated that employing an entropy inner product in the POD process resulted in a dimensionally-consistent POD formulation that provided bases with a similar approximation power to a non-dimensional $\NonDimensionalLTwoSymbol$ inner product. Next, we demonstrated that the Galerkin ROM built off of entropy variables led to significant improvements over the standard Galerkin ROM. Lastly, we demonstrated that least-squares ROMs which minimize the residual in an entropy-based norm were significantly more accurate than least-squares ROMs that minimize the residual in a non-dimensional $\NonDimensionalLTwoSymbol$ norm. Overall, the most accurate ROM formulations were the least-squares ROMs formulated in the entropy inner product. These methods provided solutions that were qualitatively accurate and associated with the lowest errors across every exemplar.

The results of this work demonstrated the promise of ROMs based on physics-based inner products, and in particular demonstrated the efficacy of entropy inner products for the compressible Euler equations. While we believe the presented results are promising, we emphasize that important aspects of a practical ROM formulation, such as online efficiency through hyper-reduction and certification via error analysis, were outside the scope this study. As a result, future work will investigate theoretical stability properties of entropy ROMs, integrating them with hyper-reduction and/or other techniques to enhance their online efficiency, and application of physics-based inner products to other domain areas.

\section{Acknowledgments}
The authors thank Irina Tezaur, Chad Sockwell, Patrick Blonigan, and John Tencer for conversations from which this work benefited.
This work was sponsored by Sandia's Advanced Simulation and
Computing (ASC) Verification and Validation (V\&V) Project/Task
\#103723/05.30.02 and 103723/06.20.05.  This paper describes objective technical results and
analysis. Any subjective views or opinions that might be expressed in the
paper do not necessarily represent the views of the U.S.\ Department of Energy
or the United States Government.  Sandia National Laboratories is a
multimission laboratory managed and operated by National Technology and
Engineering Solutions of Sandia, LLC., a wholly owned subsidiary of Honeywell
International, Inc., for the U.S.\  Department of Energy's National Nuclear
Security Administration under contract DE-NA-0003525.

%
%
%

%% file: appendix.tex
\section{Entropy Jacobian}\label{sec:appendix}
The Jacobian of the conserved variables with respect to the entropy variables is
$$
\frac{\partial \qConserved}{\partial \qEntropy} \equiv \conservativeToEntropy = \begin{bmatrix}
\qConserved_1 & \qConserved_2 & \qConserved_3 & \qConserved_4 \\
 
\text{sym.} & \frac{\qConserved_2^2}{\qConserved_1} + p &  - \frac{ \qConserved_2 \qEntropy_3}{\qEntropy_4} & -\frac{\qConserved_2}{\qEntropy_4} - \frac{\qEntropy_2 \qConserved_4}{\qEntropy_4}  \\
\text{sym.} & \text{sym.} &  \frac{\qConserved_3^2}{\qConserved_0} + p & -\frac{\qConserved_3}{\qEntropy_4} - \frac{\qEntropy_3 \qConserved_4}{\qEntropy_4} \\ 
\text{sym.}  & \text{sym.}  & \text{sym.} & \qConserved_1 H^2 - \frac{a^2 p}{\gamma - 1} 
\end{bmatrix}
$$
where $H = \frac{\qConserved_4 + p}{\rho}$ and $a^2 = \frac{\gamma p}{\rho}$. 